\documentclass[reqno,a4paper,11pt]{amsart}
\usepackage[english]{babel}
\renewcommand{\baselinestretch}{1.1}
\usepackage{graphicx,mathrsfs,tikz,latexsym,ifthen,amsmath,amsfonts,amssymb,amsthm,stmaryrd,fancyhdr}
\usepackage{enumerate,enumitem,empheq,mathtools,times,tikz}%showkeys
\mathtoolsset{showonlyrefs,showmanualtags}
\usepackage[utf8x]{inputenc}
\usepackage[toc,page]{appendix}
\usepackage[pdftex]{hyperref}

\allowdisplaybreaks

\input xy
\xyoption{all}

\usepackage{color}
\definecolor{MyDarkBlue}{rgb}{0.15,0.25,0.45}
\usepackage{hyperref}
\hypersetup{
colorlinks=true,
citecolor=black,
linkcolor=black,
urlcolor=MyDarkBlue,
breaklinks=true
}

\newcommand{\vacuum}{\vert 0\rangle}

\newcommand{\hilb}[2]{\mathrm{Hilb}^{#1}(#2)}

\newcommand{\localizedi}{\otimes_{\mathbb{C}[\varepsilon_1^{(i)},\varepsilon_2^{(i)}]}\mathbb{C}\big(\varepsilon_1^{(i)},\varepsilon_2^{(i)}\big)}

\newcommand{\Acal}{{\mathcal A}}

\newcommand{\Ecal}{{\mathcal E}}

\newcommand{\Ncal}{{\mathcal N}}
\newcommand{\Mcal}{\mathcal{M}}

\newcommand{\Zcal}{{\mathcal Z}}
\newcommand{\Ocal}{{\mathcal O}}
\newcommand{\Gcal}{{\mathcal G}}
\newcommand{\Fcal}{{\mathcal F}}
\newcommand{\Ical}{{\mathcal I}}

\newcommand{\Rcal}{{\mathcal R}}
\newcommand{\Lcal}{{\mathcal L}}
\newcommand{\Vcal}{{\mathcal V}}

\newcommand{\Ucal}{\mathcal{U}}

\newcommand{\Wcal}{\mathcal{W}}

\newcommand{\hfrak}{\mathfrak{h}}
\newcommand{\cfrak}{\mathfrak{c}}
\newcommand{\tfrak}{\mathfrak{t}}
\newcommand{\Qfrak}{\mathfrak{Q}}
\newcommand{\Pfrak}{\mathfrak{P}}
\newcommand{\Pfrakhat}{\widehat{\mathfrak{P}}}
\newcommand{\pfrak}{\mathfrak{p}}
\newcommand{\qfrak}{\mathfrak{q}}
\newcommand{\slfrak}{\mathfrak{sl}}
\newcommand{\glfrakhat}{\widehat{\mathfrak{gl}}}
\newcommand{\slfrakhat}{\widehat{\mathfrak{sl}}}
\newcommand{\tfrakhat}{\widehat{\mathfrak{t}}}
\newcommand{\nfrakhat}{\widehat{\mathfrak{n}}}
\newcommand{\bfrakhat}{\widehat{\mathfrak{b}}}

\newcommand{\glfrak}{\mathfrak{gl}}
\newcommand{\virfrak}{\mathfrak{Vir}}
\newcommand{\Lfrak}{\mathfrak{L}}
\newcommand{\Ufrak}{\mathfrak{U}}

\newcommand{\R}{{\mathbb{R}}}
\newcommand{\N}{{\mathbb{N}}}
\newcommand{\C}{{\mathbb{C}}}
\newcommand{\Z}{{\mathbb{Z}}}
\newcommand{\Q}{{\mathbb{Q}}}

\newcommand{\PP}{{\mathbb{P}}}

\newcommand{\F}{{\mathbb{F}}}
\newcommand{\HH}{{\mathbb{H}}}
\newcommand{\W}{{\mathbb{W}}}

\newcommand{\Xscr}{{\mathscr{X}}}

\newcommand{\Dscr}{{\mathscr{D}}}
\newcommand{\Bscr}{{\mathscr{B}}}

\newcommand{\qsf}{\mathsf{q}}

\newcommand{\quiv}{{\tt Q}}
\newcommand{\source}{{\tt s}}
\newcommand{\tail}{{\tt t}}

\newcommand{\Ebf}{\boldsymbol{E}}

\newcommand{\ubf}{\boldsymbol{u}}
\newcommand{\Vbf}{\boldsymbol{V}}

\newcommand{\Zbf}{\boldsymbol{Z}}

\newcommand{\mubf}{\boldsymbol{\mu}}
\newcommand{\qbf}{\boldsymbol{\mathsf{q}}}
\newcommand{\nbf}{\boldsymbol{n}}
\newcommand{\lambdabf}{\boldsymbol{\lambda}}
\newcommand{\xibf}{\boldsymbol{\xi}}
\newcommand{\jbf}{\boldsymbol{j}}

\newcommand{\rk}{\operatorname{rk}}
\newcommand{\cc}{\operatorname{c}}
\newcommand{\Trr}{\operatorname{Tr}}
\newcommand{\id}{\operatorname{id}}

\newcommand{\Hom}{\operatorname{Hom}}

\newcommand{\crm}{\operatorname{c}}
\newcommand{\End}{\operatorname{End}}

\newcommand{\Pic}{\operatorname{Pic}}
\newcommand{\Ext}{\operatorname{Ext}}
\newcommand{\V}{\operatorname{V}}
\newcommand{\eu}{\operatorname{eu}}
\newcommand{\ch}{\operatorname{ch}}
\def\e{{\,\rm e}\,}
\def\ii{{\,{\rm i}\,}}

\newcommand{\triend}{\parbox{2mm}{\hfill} \hfill\mbox{\hspace{0.2mm}}\hfill$\triangle$}
\newcommand{\ocend}{\parbox{2mm}{\hfill} \hfill\mbox{\hspace{0.2mm}}\hfill$\oslash$}

\newtheorem{theorem}[equation]{Theorem}
\newtheorem{proposition}[equation]{Proposition}
\newtheorem{lemma}[equation]{Lemma}
\newtheorem{corollary}[equation]{Corollary}
\newtheorem*{corollary*}{Corollary}
\newtheorem*{theorem*}{Theorem}
\newtheorem*{proposition*}{Proposition}

\numberwithin{equation}{section}

\theoremstyle{remark}
\newtheorem{example}[equation]{Example}

\theoremstyle{remark}
\newtheorem{rem}[equation]{Remark}
\newenvironment{remark}{\begin{rem}}{\triend\end{rem}}

\theoremstyle{definition}
\newtheorem{defin}[equation]{Definition}
\newenvironment{definition}{\begin{defin}}{\ocend\end{defin}}

\begin{document}

\begin{flushright}
EMPG--14--10
\end{flushright}

\vskip 1cm

\title[AGT on ALE: The abelian case]{\large{AGT relations for abelian quiver gauge theories on ALE spaces}}

\vskip 1cm

\maketitle \thispagestyle{empty}

\begin{center}
{\large \sc Mattia Pedrini$^{\S\ddag}$,  Francesco Sala$^{\bullet\P\star}$} \ and \ {\large \sc Richard J. Szabo$^{\P\star\circ}$} \\[8pt] 
$^\S$ Scuola Internazionale Superiore di Studi Avanzati {\sc (SISSA),}\\ Via Bonomea 265, 34136
Trieste, Italia; \\[3pt] $^\ddag$ Istituto Nazionale di Fisica Nucleare, Sezione di Trieste \\[3pt]
$^\bullet$ Department of Mathematics, The University of Western Ontario,\\
Middlesex College, London N6A 5B7, Ontario, Canada;\\[3pt]
$^\P$ Department of Mathematics, Heriot-Watt University,\\ Colin
Maclaurin Building, Riccarton, Edinburgh EH14 4AS, United Kingdom; \\[3pt]
$^\star$ Maxwell Institute for Mathematical Sciences, Edinburgh,
United Kingdom; \\[3pt] $^\circ$ The Tait Institute, Edinburgh, United Kingdom
 \end{center}

\par\vfill

\noindent \begin{quote} 
{\sc Abstract.} \small
We construct level one dominant representations of the affine Kac-Moody
algebra $\glfrakhat_k$ on the equivariant cohomology groups of moduli
spaces of rank one framed sheaves on the orbifold compactification of
the minimal resolution $X_k$ of the $A_{k-1}$ toric singularity $\C^2/\Z_k$. We show that the direct sum of the fundamental classes of these moduli spaces is a Whittaker vector for $\glfrakhat_k$, which proves the AGT correspondence for pure $\Ncal=2$ $U(1)$ gauge theory on $X_k$. We consider Carlsson-Okounkov type Ext-bundles over products of the moduli spaces and use their Euler classes to define vertex operators. Under the decomposition $\glfrakhat_k\simeq \hfrak\oplus \slfrakhat_k$, these vertex operators decompose as products of bosonic exponentials associated to the Heisenberg algebra $\hfrak$ and primary fields of $\slfrakhat_k$. We use these operators to prove the AGT correspondence for $\Ncal=2$ superconformal abelian quiver gauge theories on $X_k$. 
\end{quote}

\par
\vfill
\parbox{.95\textwidth}{\hrulefill}\par
\noindent \begin{minipage}[c]{\textwidth}\parindent=0pt \renewcommand{\baselinestretch}{1.2}
\small
\emph{Date:}   November 2015  \par 
 \par
\emph{E-Mail:} \texttt{mattia.pedro@gmail.com,  salafra83@gmail.com, R.J.Szabo@hw.ac.uk} \par
\end{minipage}

\newpage

\setlength{\parskip}{0.2ex}

{\baselineskip=14pt
\tableofcontents
}

\setlength{\parskip}{0.6ex plus 0.3ex minus 0.2ex}

\newpage

\section{Introduction and summary}

\subsection{AGT relations and ALE spaces}

In this paper we study a new occurrence of the deep relations between the moduli theory of  sheaves and the representation theory of affine/vertex algebras. 

We are particularly interested in the kind of relations which come
from gauge theory considerations. An important example of these
relations is the AGT correspondence for gauge theories on $\R^4$: in
\cite{art:aldaygaiottotachikawa2010} Alday, Gaiotto and Tachikawa
conjectured a relation between the instanton partition functions of
$\Ncal=2$ supersymmetric quiver gauge theories on $\R^4$ and the
conformal blocks of two-dimensional $A_{r-1}$ Toda conformal field theories (see also
\cite{art:wyllard2009,art:aldaytachikawa2010}); this conjecture has
been explicitly confirmed in some special cases, see
e.g. \cite{art:mironovmorozovshakirov2011,
  art:albafateevlitvinovtarnopolskiy2011, art:tan2013,
  art:aganagichaouzishakirov2014}. From a mathematical perspective,
this correspondence implies: (1) the existence of a representation of the
W-algebra $\Wcal(\glfrak_r)$ on the equivariant cohomology of the
moduli spaces $\Mcal(r,n)$ of framed sheaves on the projective plane $\PP^2$ of rank $r$
and second Chern class $n$ such that the latter is isomorphic to a Verma module of $\Wcal(\glfrak_r)$; (2) the fundamental classes of $\Mcal(r,n)$ give a Whittaker vector of
$\Wcal(\glfrak_r)$ (pure gauge theory); (3) the Ext vertex operator is related to a certain ``intertwiner" of $\Wcal(\glfrak_r)$ under the isomorphism stated in (1) (quiver gauge theory). The instances (1) and (2) were proved by
Schiffmann and Vasserot \cite{art:schiffmannvasserot2013}, and
independently by Maulik and Okounkov
\cite{art:maulikokounkov2012}. For $r=1$, the moduli space
$\Mcal(1,n)$ is isomorphic to the Hilbert scheme of $n$ points on $\C^2$ and $\Wcal(\glfrak_1)$ is the W-algebra associated with an infinite-dimensional Heisenberg algebra; the AGT correspondence for pure $U(1)$ gauge theory reduces to the famous result of Nakajima \cite{art:nakajima1997, book:nakajima1999} in the equivariant case \cite{art:vasserot2001, art:liqinwang2004, art:nakajima2014}. Presently, (3) has been proved only in the rank one case \cite{art:carlssonokounkov2012} and in the rank two case \cite{art:carlsson2015,art:negut2015}.

In this paper we are interested in the AGT correspondence for $\Ncal=2$ quiver gauge theories on ALE spaces associated with the Dynkin diagram of type $A_{k-1}$ for $k\geq 2$. The corresponding instanton partition functions are defined in terms of equivariant cohomology classes over Nakajima quiver varieties of type the affine Dynkin diagram $\widehat{A}_{k-1}$. These quiver varieties depend on a real stability parameter $\xi_\R$, which lives in an open subset of $\R^k$ having a ``chambers" decomposition: if two real stability parameters belong to the same chamber, the corresponding quiver varieties are (equivariantly) isomorphic; otherwise, the corresponding quiver varieties are only $\C^\ast$-diffeomorphic. Therefore, the pure gauge theories partition functions should be all nontrivially equivalent, while the partition functions for quiver gauge theories should satisfy ``wall-crossing" formulas (cf.\ \cite{art:belavinbershteinfeiginlitvinovtarnopolsky2011,art:itomaruyoshiokuda2013}).

By looking at instanton partition functions of pure gauge theories associated with moduli spaces of $\Z_k$-equivariant framed sheaves on $\PP^2$ (which are quiver varieties depending on a so-called ``level zero chamber"), the authors of \cite{art:belavinfeigin2011,art:nishiokatachikawa2011, art:belavinbelavinbershtein2011} conjectured an extension of the AGT correspondence in the A-type ALE case as a relation between instanton partition functions of  $\Ncal=2$ quiver gauge theories and conformal blocks of Toda-like conformal field theories with $\Z_k$ parafermionic symmetry. In particular, the pertinent algebra to consider in this case is the coset
\begin{equation}
\Acal(r,k):= \frac{\glfrakhat_N}{\glfrakhat_{N-k}}
\end{equation}
acting at level $r$, where $N$ is related to the equivariant parameters. For $r=1$ the algebra $\Acal(1,k)$ is simply $\glfrakhat_k$ acting at level one. In general, $\Acal(r,k)$ is isomorphic to the direct sum of the affine Lie algebra $\glfrakhat_k$ acting at level $r$ and the $\Z_k$-parafermionic $\Wcal(\glfrak_r)$-algebra. Checks of the conjecture has been done \cite{art:wyllard2011, art:ito2011} by using partition functions of pure gauge theories associated with moduli spaces of $\Z_k$-equivariant framed sheaves on $\PP^2$. In \cite{art:bonellimaruyoshitanzini2011, art:bonellimaruyoshitanzini2012} the authors studied in details $\Ncal=2$ quiver gauge theories on the minimal resolution $X_2$ of the Kleinian singularity $\C^2/\Z_2$ and provided evidences for the conjecture: in this case, the quiver variety depends on a so-called ``level infinity chamber" and corresponds to moduli spaces of framed sheaves on a suitable stacky compactification of $X_2$. In the $k=2$ case, a comparison of these approaches using different stability chambers is done in \cite{art:alfimovbelavintarnopolsky2013}; further speculations in the arbitrary $k$ case are in \cite{art:bonellimaruyoshitanziniyagi2012}.

Mathematically, this correspondence should imply: (1) the existence of a representation of the
coset $\Acal(r,k)$ on the equivariant cohomology of Nakajima quiver varieties associated with the affine A-type Dynkin diagram such that the latter is isomorphic to a Verma module of $\Acal(r,k)$; (2) the fundamental classes of the quiver varieties give a Whittaker vector of
$\Acal(r,k)$ (pure gauge theory); (3) the Ext vertex operator is related to a certain ``intertwiner" of $\Acal(r,k)$ under the isomorphism stated in (1) (quiver gauge theory). As pointed out in \cite{art:alfimovbelavintarnopolsky2013}, different chambers should provide different realizations of the action conjectured in (1). On the other hand, the conjectural wall-crossing behavior of the instanton partition functions for quiver gauge theories \cite{art:itomaruyoshiokuda2013} should be related by a similar behavior of the Ext vertex operators by varying of the stability chambers.

The ALE space we consider in this paper is the minimal resolution $X_k$ of the simple Kleinian singularity $\C^2/\Z_k$. In \cite{art:bruzzopedrinisalaszabo2013} an orbifold compactification $\Xscr_k$ of $X_k$ is constructed by adding a smooth divisor $\Dscr_\infty$, which lays the foundations for a new sheaf theory approach to the study of $U(r)$ instantons on $X_k$ (cf.\ \cite{art:eyssidieuxsala2013}). Moduli spaces of sheaves on $\Xscr_k$ framed along $\Dscr_\infty$ are also constructed in \cite{art:bruzzopedrinisalaszabo2013}; by using these moduli spaces we have a new sheaf theory approach to the study of Nakajima quiver varieties with the stability parameter of $X_k$ and, consequently, of $U(r)$ gauge theories on ALE spaces of type $A_{k-1}$ which are isomorphic to $X_k$. In the present paper we use this new approach to study the AGT correspondence for abelian quiver gauge theories on $X_k$: from a physics point of view we prove the relations between instanton partition functions and conformal blocks and from a mathematical point of view we prove (1), (2) and (3).

\subsection{Summary of results}

Let us now summarize our main results. Recall that the compactification $\Xscr_k$ is a two-dimensional projective toric orbifold with Deligne-Mumford torus $T:=\C^\ast\times \C^\ast$; the complement $\Xscr_k\setminus X_k$ is a smooth Cartier divisor $\Dscr_\infty$ endowed with the structure of a $\Z_k$-gerbe. There exist line bundles $\Ocal_{\Dscr_\infty}(j)$ on $\Dscr_\infty$, for $j=0, 1, \ldots, k-1$, endowed with unitary flat connections associated with the irreducible unitary representations of $\Z_k$. Hence by \cite[Theorem 6.9]{art:eyssidieuxsala2013} locally free sheaves on $\Xscr_k$ which are isomorphic along $\Dscr_\infty$ to $\Ocal_{\Dscr_\infty}(j)$ correspond to $U(1)$ instantons on $X_k$ with holonomy at infinity given by the $j$-th irreducible unitary representation of $\Z_k$, for $j=0,1, \ldots, k-1$.

Fix $j=0, 1, \ldots, k-1$. A rank one $(\Dscr_\infty, \Ocal_{\Dscr_\infty}(j))$-framed sheaf on $\Xscr_k$ is a pair $(\Ecal, \phi_{\Ecal})$, where $\Ecal$ is a rank one torsion free sheaf on $\Xscr_k$, locally free in a neighbourhood of $\Dscr_\infty$, and $\phi_{\Ecal}\colon \Ecal\big\vert_{\Dscr_\infty} \xrightarrow{\sim} \Ocal_{\Dscr_\infty}(j)$ is an isomorphism. Let $\Mcal(\vec{u},n,j)$ be the fine moduli space parameterizing isomorphism classes of rank one $(\Dscr_\infty, \Ocal_{\Dscr_\infty}(j))$-framed sheaves on $\Xscr_k$, with first Chern class given by $\vec{u}\in \Z^{k-1}$ and second Chern class $n$. As explained in Remark \ref{rem:firstchern}, the vector $\vec{u}$ is canonically associated with an element $\gamma_{\vec{u}}+\omega_j\in \Qfrak+\omega_j$, where $\Qfrak$ is the root lattice of the Dynkin diagram of type $A_{k-1}$ and $\omega_j$ is the $j$-th fundamental weight of type $A_{k-1}$. We denote by $\Ufrak_j$ the set of vectors $\vec{u}$ associated with $\gamma+\omega_j$ for some $\gamma\in \Qfrak$.

The moduli space $\Mcal(\vec{u},n,j)$ is a smooth quasi-projective variety of dimension $2n$. On $\Mcal(\vec{u},n;j)$ there is a natural $T$-action induced by the toric structure of $\Xscr_k$. Let $\varepsilon_1, \varepsilon_2$ be the generators of the $T$-equivariant cohomology of a point and consider the localized equivariant cohomology
\begin{equation}
\W_{\vec{u},j}:=\bigoplus_{n\geq0}\,  H^\ast_T\big(\Mcal(\vec{u},n,j)\big)\otimes_{\C[\varepsilon_1, \varepsilon_2]}\C(\varepsilon_1, \varepsilon_2)\ .
\end{equation}
Define also the total localized equivariant cohomology by summing over all vectors $\vec{u}\in\Ufrak_j$:
\begin{equation}
\W_j:=\bigoplus_{\vec{u}\in \Ufrak_j}\, \W_{\vec{u},j}\ .
\end{equation}
The affine Lie algebra $\glfrakhat_k$ acts on $\W_j$ as follows (see Proposition \ref{prop:representation} and Proposition \ref{prop:Wkweight}).
\begin{proposition*}
There exists a $\glfrakhat_k$-action on $\W_j$ under which it is the $j$-th dominant representation of $\glfrakhat_k$ at level one, i.e., the highest weight representation of $\glfrakhat_k$ with fundamental weight $\widehat{\omega}_j$ of type $\widehat{A}_{k-1}$. Moreover, the weight spaces of $\W_j$ with respect to the $\glfrakhat_k$-action are the $\W_{\vec{u},j}$ with weights $\gamma_{\vec{u}}+\omega_j$.
\end{proposition*}

The vector spaces $\W_{\vec{u},j}$ also have a representation theoretic intepretation.
\begin{corollary*}[Corollary \ref{cor:hwrep}]
$\W_{\vec{u},j}$ is a highest weight representation of the Virasoro algebra associated with $\glfrakhat_k$ of conformal dimension $\Delta_{\vec{u}}:=\frac{1}{2}\, \vec{u}\cdot C^{-1}\vec{u}$, where $C$ is the Cartan matrix of type $A_{k-1}$.
\end{corollary*}

The representation is constructed by using a vertex algebra approach via the Frenkel-Kac construction. A similar construction for the cohomology groups of moduli spaces of rank one torsion free sheaves over smooth projective surfaces is outlined in \cite[Chapter 9]{book:nakajima1999}. In \cite{art:nagao2009}, Nagao analysed vertex algebra realizations of representations of $\slfrakhat_k$ on the equivariant cohomology groups of Nakajima quiver varieties associated with the affine Dynkin diagram $\widehat{A}_{k-1}$, for an integer $k\geq 2$, with dimension vector corresponding to the trivial holonomy at infinity $j=0$; in this case the pertinent representation is the basic representation of $\slfrakhat_k$.

In the following we describe our AGT relations, which connect together $\W_j$ for $j=0,1, \ldots, k-1$, the action of $\glfrakhat_k$ on $\W_j$ and abelian quiver gauge theories on $X_k$. The first relation we obtain concerns the pure gauge theory. Let $\Zcal_{X_k}\big(\varepsilon_1,\varepsilon_2; \qsf, \vec{\xi} \ \big)_j$ be the instanton partition function for the pure $\Ncal=2$ $U(1)$ gauge
theory on the ALE space $X_k$ with fixed holonomy at infinity given by the $j$-th irreducible representation of $\Z_k$ (see Section \ref{sec:N=2Xk}). It has the following representation theoretic characterization.
\begin{theorem*}[AGT relation for pure $\Ncal=2$ $U(1)$ gauge theory]
The Gaiotto state 
\begin{equation}
G_j:= \sum_{\vec{u}\in\Ufrak_j} \ \sum_{n\geq 0} \,
\big[\Mcal(\vec{u},n,j) \big]_T
\end{equation}
is a Whittaker vector for $\glfrakhat_k$. Moreover, the weighted norm of the weighted Gaiotto state
\begin{equation}
G_j\big(\qsf,\vec{\xi} \ \big) := \sum_{\vec{u}\in\Ufrak_j} \ \sum_{
  n\geq 0} \, \qsf^{n+\frac{1}{2}\, \vec{u} \cdot C^{-1}\vec{u}}\ 
\vec{\xi}^{\ C^{-1}\vec{u}}\ \big[\Mcal(\vec{u},n,j) \big]_T
\end{equation}
is exactly $\Zcal_{X_k}\big(\varepsilon_1,\varepsilon_2; \qsf, \vec{\xi} \ \big)_j$.
\end{theorem*}

We also consider $\Ncal=2$ superconformal quiver gauge theories with gauge group $U(1)^{r+1}$ for some $r\geq0$. By the ADE classification in \cite[Chapter 3]{art:nekrasovpestun2012} the admissible quivers in this case are the linear quivers of the finite-dimensional $A_r$-type Dynkin diagram and the cyclic quivers of the affine $\widehat{A}_{r}$-type extended Dynkin
diagram. In order to state AGT relations in these cases, we introduce Ext vertex operators \cite{art:carlssonokounkov2012,art:carlsson2015, art:negut2015}. Consider the element $\Ebf_\mu\in K\big(\Mcal(\vec{u}_1,n_1,j_1)\times\Mcal(\vec{u}_2,n_2,j_2)\big)$ whose fibre over a point $\big([(\Ecal,\phi_\Ecal)]$, $[(\Ecal',\phi_{\Ecal'})]\big)$ is
\begin{equation}
\big(\Ebf_\mu\big)_{\left([(\Ecal,\phi_\Ecal)]\,,\,[(\Ecal',\phi_{\Ecal'})]\right)}=\Ext^1\big(\Ecal,\Ecal'\otimes
\Ocal_{\Xscr_k}(\mu)\otimes\Ocal_{\Xscr_k}(-\Dscr_\infty) \big)\ ,
\end{equation}
where $\Ocal_{\Xscr_k}(\mu)$ is the trivial line bundle on $\Xscr_k$ on which the torus $T_\mu=\C^\ast$ acts by scaling the fibres with $H_{T_\mu}^\ast({\rm pt};\C)=\C[\mu]$. By using the Euler class of $\Ebf_\mu$ we define a vertex operator $\V_\mu(\vec x,z)\in\End\big(\bigoplus_{j=0}^{k-1}\,\W_j\, \big)[[z^{\pm\, 1},$ $x_1^{\pm\,1},$ $\dots,x_{k-1}^{\pm\,1}]]$ (see Section \ref{se:bifundXk}). Under the decomposition $\glfrakhat_k=\hfrak\oplus \slfrakhat_k$, we have the following characterization of $\V_\mu(\vec x,z)$ in terms of vertex operators depending respectively on $\hfrak$ and $\slfrakhat_k$.
\begin{theorem*}[Theorem \ref{thm:virasoroprimary}]
The vertex operator $\V_\mu(\vec x, z)$ can be expressed in the form
\begin{multline}
\V_\mu(\vec x, z)=\V_{-\frac{\mu}{\sqrt{-k\, \varepsilon_1\, \varepsilon
_2}},\frac{\mu+\varepsilon_1+\varepsilon_2}{\sqrt{-k\, \varepsilon_1\,
\varepsilon_2}}
}(z) \\ \otimes \ \sum_{j_1,j_2=0}^{k-1} \
\sum_{\vec{u}_1\in\Ufrak_{j_1},\vec{u}_2\in\Ufrak_{j_2}}\
\bar \V_\mu(\vec{v}_{21},\vec x, z)\, z^{\Delta_{\vec u_2}-\Delta_{\vec u_1}}\ \exp\big(\log z\ \cfrak-\gamma_{21}\,\big)\,\exp\big(\gamma_{21}\,\big)\big\vert_{ \W_{\vec{u}_1, j_1}}
\ ,
\end{multline}
where $\V_{\alpha,\beta}(z)$ denotes a generalized bosonic exponential associated with the Heisenberg algebra $\hfrak$ (see Definition \ref{def:bosonops}), $\exp\big(\log z\ \cfrak-\gamma_{21}\,\big)\,\exp\big(\gamma_{21}\,\big)$ is the vertex operator on $\W_{j_1}$ defined in Equation \eqref{eq:expopWj}, and $\bar \V_\mu(\vec{v}_{21},\vec x, z)$ is the primary field \eqref{eq:barVmudef} of the Virasoro algebra
associated with $\slfrakhat_k$ with conformal
 dimension $\Delta_{\vec u_2-\vec u_1}=\frac12\, \vec v_{21}\cdot C\vec v_{21}$, where $\vec v_{21}:= C^{-1}(\vec u_2-\vec u_1)$ for $j_1, j_2=0,1, \ldots, k-1$ and $\vec{u}_1\in\Ufrak_{j_1},\vec{u}_2\in\Ufrak_{j_2}$.
\end{theorem*}
For $j_1,j_2=0,1, \ldots, k-1$ denote by $\V_\mu^{j_1,j_2}(\vec x, z)$ the restriction of the vertex operator
$\V_\mu(\vec x, z)$ to $\Hom(\W_{j_1}, \W_{j_2})[[z^{\pm\,1},x_1^{\pm\,1},\dots,x_{k-1}^{\pm\,1}]]$.

Let $\Zcal_{X_k}^{\widehat{A}_r}\big(\varepsilon_1,\varepsilon_2,\mubf;\qbf, \vec\xibf
\ \big)_{\jbf}$ be the instanton partition function for the
$\Ncal=2$ superconformal $U(1)^{r+1}$ quiver gauge theory of type $\widehat{A}_{r}$ with holonomy at infinity associated with $\jbf:=(j_0, j_1, \ldots, j_r)$, topological couplings $\qsf_\upsilon\in\C^*$ and $\vec\xi_\upsilon \in (\C^*)^{k-1}$ for
$\upsilon=0,1, \ldots, r$, and masses $\mubf:=(\mu_0,\mu_1, \ldots, \mu_r)$. We prove the following AGT relation.
\begin{theorem*}[AGT relation for
$\Ncal=2$ $U(1)^{r+1}$ quiver gauge theory of type $\widehat{A}_{r}$]
The partition function of the $\widehat{A}_r$-theory on $X_k$ is given
by
\begin{equation}
\Zcal_{X_k}^{\widehat{A}_r}\big(\varepsilon_1,\varepsilon_2,\mubf;\qbf, \vec\xibf
\ \big)_{\jbf} = \Trr_{\W_{j_0}} \, \qsf^{L_0}\ \vec\xi{}^{\ C^{-1}\vec h} \ \prod_{\upsilon=0}^r\,
\V_{\mu_\upsilon}^{j_\upsilon,j_{\upsilon+1}}(\vec x_\upsilon,z_\upsilon)\ \delta^{\rm
  conf}_{\upsilon,{\upsilon+1}} \ ,
\end{equation}
where
$\qsf:= \qsf_0\, \qsf_1\cdots \qsf_r$, $( \, \vec\xi\ )_i:=
(\vec\xi_0)_i\, (\vec\xi_1)_i\, \cdots (\vec\xi_r)_i$, $z_\upsilon:=
z_0\, \qsf_1\cdots \qsf_\upsilon$, and $(\vec x_\upsilon)_i:= (\vec
x_0)_i\, (\vec\xi_1)_i\cdots$ $ (\vec\xi_\upsilon)_i$ for
$\upsilon=1,\dots,r$ and $i=1,\dots,k-1$. Here $L_0$ is the Virasoro energy operator associated to $\glfrakhat_k$, $\vec h=(h_1,\dots,h_{k-1})$ are the generators of the Cartan subalgebra of $\slfrak_k$, and $\delta^{\rm
  conf}_{\upsilon,{\upsilon+1}}$ is the conformal restriction operator defined in Equation \eqref{eq:confrestrop}.
\end{theorem*}
We also get a characterization of $\Zcal_{X_k}^{\widehat{A}_r}\big(\varepsilon_1,\varepsilon_2,\mubf;\qbf, \vec\xibf
\ \big)_{\jbf}$ in terms of the corresponding partition function on $\C^2$ and a part depending only on $\slfrakhat_k$.
\begin{corollary*} 
Let $\V_\mu(\vec{v}_{21},\vec x,
z):=  \, z^{\Delta_{\vec u_2}-\Delta_{\vec u_1}} {\bar \V}_\mu(\vec{v}_{21}, \vec x, z)\ \exp\big(\log z\ \cfrak-\gamma_{21} \big)\, \exp\big(\gamma_{21}\,\big)$. Then we have
\begin{multline}
\Zcal_{X_k}^{\widehat{A}_r}\big(\varepsilon_1,\varepsilon_2,\mubf;\qbf, \vec\xibf
\ \big)_{\jbf}
= \Zcal^{\widehat{A}_{r}}_{\C^2}(\varepsilon_1,\varepsilon_2,\mubf;\qsf)^{\frac{1}{k}}\ \qsf^{\frac{1}{24}\,(1-\frac{1}{k})}\,\eta(\qsf)^{\frac{1}{k}-1} \\
\times\ \Trr_{\Vcal(\,\widehat{\omega}_{j_0}\,)} \,
\qsf^{L_0^{\slfrakhat_k}}\ \vec\xi{}^{\ C^{-1}\vec h}\ \prod_{\upsilon=0}^r\ 
\sum_{(\vec{u}_\upsilon\in\Ufrak^{\rm conf}_{j_\upsilon})} \,
\V_{\mu_\upsilon}(\vec{v}_{\upsilon,\upsilon+1},\vec x_\upsilon, z_\upsilon) \, \big\vert
_{\W_{\vec{u}_\upsilon,j_\upsilon}} \ ,
\end{multline}
where $\eta(\qsf)$ is the Dedekind function, $\Vcal(\,\widehat{\omega}_{j_0}\,)$ is the $j_0$-th dominant representation of $\slfrakhat_k$ and $\Ufrak^{\rm conf}_{j_\upsilon}$ is the subset of $\Ufrak_{j_\upsilon}$ defined in Equation \eqref{eq:confchargeset}.
\end{corollary*}

Let $\Zcal_{X_k}^{{A}_r}\big(\varepsilon_1,\varepsilon_2,\mubf;\qbf, \vec\xibf
\ \big)_{\jbf}$ be the instanton partition function for the
$\Ncal=2$ superconformal $U(1)^{r+1}$ quiver gauge theory of type ${A}_{r}$ with holonomy at infinity associated with $\jbf:=(j_0,j_1, \ldots, j_r)$. We also prove the following AGT relation.
\begin{theorem*}[AGT relation for $\Ncal=2$ $U(1)^{r+1}$ quiver gauge theory of type ${A}_{r}$]
The partition function of the $A_r$-theory on $X_k$ is given
by
\begin{multline}
\Zcal_{X_k}^{{A}_r}\big(\varepsilon_1,\varepsilon_2,\mubf;\qbf, \vec\xibf
\ \big)_{\jbf} \\
= \Big\langle |0\rangle_{\rm conf} \, , \, V_{\mu_0}(\vec x_0,z_0) \, \Big(\, \prod_{\upsilon=1}^{r}\,
\V_{\mu_{\upsilon}}^{j_{\upsilon}-1, j_{\upsilon}}(\vec x_\upsilon,z_\upsilon)\ \delta^{\rm
  conf}_{{\upsilon-1},{\upsilon}} \, \Big)\, \V_{\mu_{r+1}}(\vec x_{r+1},z_{r+1}) |0\rangle_{\rm conf} \Big\rangle_{\bigoplus_{j=0}^{k-1}\,\W_j} \ ,
\end{multline}
where
$z_\upsilon:=
z_0\, \qsf_0\, \qsf_1\cdots \qsf_\upsilon$ and $(\vec x_\upsilon)_i:= (\vec
x_0)_i\, (\vec\xi_0)_i\, (\vec\xi_1)_i$ $\cdots (\vec\xi_{\upsilon-1})_i$ for
$\upsilon=1,\dots,r+1$, $i=1,\dots,k-1$, and $|0\rangle_{\rm conf}:= \prod_{\upsilon=0}^r\,\delta^{\rm
  conf}_{0,\upsilon}\triangleright [\emptyset, \vec{0}\,]$ with
$[\emptyset, \vec{0}\,]$ the vacuum vector of the fixed point basis of
$\bigoplus_{j=0}^{k-1}\,\W_j$.
\end{theorem*}
Denote by $\Vcal$ the direct sum of the $k$ level one dominant representations of $\slfrakhat_k$. Similarly to before, we have the following characterization.
\begin{corollary*} We have
\begin{multline}
\Zcal_{X_k}^{{A}_r}\big(\varepsilon_1,\varepsilon_2,\mubf;\qbf, \vec\xibf
\ \big)_{\jbf} 
= \Zcal_{\C^2}^{{A}_r}(\varepsilon_1,\varepsilon_2,\mubf;\qbf)^{\frac1k} \\ \shoveleft{ \times \
\Big\langle |0\rangle_{\rm conf} \, , \, \Big(\, \sum_{j_0,j_0'=0}^{k-1} \
\sum_{\vec{u}_0\in\Ufrak_{j_0},\vec{u}_0'\in\Ufrak_{j_0'}} \
\V_{\mu_0}(\vec{v}_{0',0},\vec x_0, z_0)\, \big\vert
_{\W_{\vec{u}_0,j_0}} \, \Big) } \\ \times \
\prod_{\upsilon=1}^{r} \ \sum_{(\vec u_\upsilon\Ufrak^{\rm conf}_{j_\upsilon})}\, 
\V_{\mu_{\upsilon}}(\vec v_{\upsilon-1,\upsilon},\vec x_\upsilon,z_\upsilon)\, \big|_{\W_{\vec u_\upsilon,j_\upsilon}} \\ \times \ \Big(\, 
\sum_{j_{r+1},j_{r+1}'=0}^{k-1} \
\sum_{\vec{u}_1\in\Ufrak_{j_1},\vec{u}_1'\in\Ufrak_{j_1'}} \
\V_{\mu_{r+1}}(\vec{v}_{1',1},\vec x_{r+1}, z_{r+1})\, \big\vert
_{\W_{\vec{u}_1,j_1}}\, \Big)
|0\rangle_{\rm conf} \Big\rangle_{\Vcal} \ .
\end{multline}
\end{corollary*}

Another important aspect of the AGT correspondence that we address in
this paper is the relation of our construction with quantum integrable
systems. In particular, for any $j=0,1, \ldots, k-1$ we define an
infinite system of commuting operators which are diagonalized in the
fixed point basis of $\W_j$; geometrically these operators correspond
to multiplication by equivariant cohomology classes (see Section
\ref{sec:integrals}). The eigenvalues of these operators with respect
to this basis can be decomposed into a part associated with $k$
non-interacting Calogero-Sutherland models and a part which can be
interpreted as particular matrix elements of the vertex operators
$\V_\mu(\vec x,z)$ in highest weight vectors of $\glfrakhat_k$. The
significance of this property is that this special orthogonal
basis manifests itself in the special integrable structure of the
two-dimensional conformal field theory and yields completely factorized matrix elements of composite vertex
operators explicitly in terms of simple rational functions of the
basic parameters, which from the gauge theory perspective represent the contributions of bifundamental matter fields.

The study of the AGT relation for pure $\Ncal=2$ $U(1)$ gauge theories and the problem of constructing commuting operators associated with $\glfrakhat_k$ is also addressed in \cite{art:belavinbershteintarnopolsky2013} from another point of view: there they consider the ``conformal" limit of the Ding-Iohara algebra, depending on parameters $q,t$, for $q,t $ approaching a primitive $k$-th root of unity and relate the representation theory of this limit to the AGT correspondence. However, their point of view is completely algebraic, so unfortunately it is not clear to us how to geometrically construct the action of the conformal limit on the equivariant cohomology groups.

\subsection{Outline}

This paper is structured as follows. In Section \ref{sec:combprelim} we briefly recall the relevant combinatorial notions that we use in this paper. In Section \ref{sec:infiniteliealg} we collect preliminary material on Heisenberg algebras and affine Lie algebras of type $\widehat{A}_{k-1}$, giving particular attention to the Frenkel-Kac construction of level one dominant representations of $\slfrakhat_k$ and $\glfrakhat_k$. In Section \ref{sec:AGTonR4} we review the AGT relations for $\Ncal=2$ superconformal abelian quiver gauge theories on $\R^4$. In Section \ref{sec:sheaves} we briefly recall the construction of the orbifold compactification $\Xscr_k$ and of moduli spaces of framed sheaves on $\Xscr_k$ from \cite{art:bruzzopedrinisalaszabo2013}. Section \ref{sec:glkreps} addresses the construction of the action of $\glfrakhat_k$ on $\W_j$ for $j=0,1, \ldots, k-1$: we perform a vertex algebra construction of the representation by using the Frenkel-Kac theorem. In Section \ref{sec:chiralvertex} we define the virtual bundle $\Ebf_\mu$ and the vertex operator $\V_\mu(\vec x,z)$, and we characterize it in terms of vertex operators of an infinite-dimensional Heisenberg algebra $\hfrak$ and primary fields of $\slfrakhat_k$ under the decomposition $\glfrakhat_k=\hfrak\oplus \slfrakhat_k$; moreover, we geometrically define an infinite system of commuting operators. In Section \ref{sec:quivergaugeXk} we prove our AGT relations, and furthermore provide expressions for our partition functions in terms of the corresponding partition functions on $\C^2$ and a part depending only on $\slfrakhat_k$. The paper concludes with two Appendices containing some technical details of the constructions from the main text: in Appendix \ref{app:virasoroprimary} we give the proof that the vertex operator $\bar \V_\mu(\vec{v}_{21},\vec x, z)$ is a primary field, while in Appendix \ref{app:edgecontributions} we recall the expressions from \cite{art:bruzzopedrinisalaszabo2013} for the edge factors which appear in the definition of $\bar \V_\mu(\vec{v}_{21},\vec x, z)$ as well as in the eigenvalues of the integrals of motion.

\subsection{Acknowledgements} 

We are grateful to M.\ Bershtein, A.\ Konechny, O.\ Schiffmann and E.\ Vasserot for helpful discussions. Also, we are indebted to the anonymous referee, whose remarks helped to improve the paper. This work was supported in part by PRIN ``Geometria delle variet\`a algebriche", by GNSAGA-INDAM, by the Grant RPG-404 from the Leverhulme Trust, and by the Consolidated Grant ST/J000310/1 from the UK Science and Technology Facilities Council. The bulk of this paper was written while the authors were staying at Heriot-Watt University in Edinburgh and at SISSA in Trieste. The last draft of the paper was written while the first and second authors were staying at IHP in Paris under the auspices of the RIP program. We thank these institutions for their hospitality and support.

\bigskip \section{Combinatorial preliminaries\label{sec:combprelim}}

\subsection{Partitions and Young tableaux}\label{sec:Young}

A \emph{partition} of a positive integer $n$ is a nonincreasing
sequence of positive numbers $\lambda=(\lambda_1\geq \lambda_2\geq
\cdots\geq \lambda_\ell> 0)$ such that
$\vert\lambda\vert:=\sum_{a=1}^\ell \, \lambda_a=n.$ We call
$\ell=\ell(\lambda)$ the \emph{length} of the partition $\lambda.$
Another description of a partition $\lambda$ of $n$ uses the notation
$\lambda=(1^{m_1}\, 2^{m_2}\, \cdots)$, where
$m_i=\#\{a\in\N\,\vert \, \lambda_a =i\}$ with $\sum_i \, i\,
m_i=n$ and $\sum_i\, m_i=\ell(\lambda)$. On the set of all partitions there is a natural
partial ordering called
\emph{dominance ordering}:
For two partitions $\mu$ and $\lambda$, we write $\mu\leq\lambda$ if
and only if $|\mu|=|\lambda|$ and
$\mu_1+\cdots+\mu_a\leq\lambda_1+\cdots+\lambda_a$ for all
$a\geq1$. We write $\mu<\lambda$ if and only if $\mu\leq\lambda$ and $\mu\neq\lambda.$ 

One can associate with a partition $\lambda$ its \emph{Young tableau},
which is the set $Y_\lambda=\{(a,b)\in\N^2\,\vert\, 1\leq a \leq
\ell(\lambda)\, ,\, 1\leq b \leq \lambda_a\}$. Then $\lambda_a$ is the
length of the $a$-th column of $Y_\lambda$; we write $\vert
Y_\lambda\vert=\vert\lambda\vert$ for the \emph{weight} of the Young
tableau $Y_\lambda$. We shall identify a partition $\lambda$ with its
Young tableau $Y_\lambda$. For a partition $\lambda$, the
\emph{transpose partition} $\lambda'$ is the partition whose Young
tableau is $Y_{\lambda'}:=\{(b,a)\in \N^2\,\vert\,(a,b)\in Y_\lambda\}$. 

The elements of a Young tableau $Y$ are called the \emph{nodes} of
$Y$. For a node $s=(a,b)\in Y$, the \emph{arm length} of $s$ is
the quantity $A(s):=A_Y(s)=\lambda_a-b$ and the \emph{leg length} of $s$ the quantity $L(s):= L_Y(s)=\lambda_b'-a$. The \emph{arm colength} and \emph{leg colength} are respectively given by $A'(s):= A_Y'(s)=b-1$ and $L'(s):= L_Y'(s)=a-1$. 

\subsection{Symmetric functions}\label{sec:symmetric}

Here we recall some preliminaries about the theory of symmetric functions in infinitely many variables which we shall use later on. Our main reference is \cite{book:macdonald1995}.

Let $\F$ be a field of characteristic zero. The \emph{algebra
  of symmetric polynomials in} $N$ \emph{variables} is the subspace
$\Lambda_{\F,N}$ of $\F[x_1,\ldots,x_N]$ which is invariant under the
action of the group of permutations $\sigma_N$ on $N$ letters. Then $\Lambda_{\F,N}$ is a graded ring: $\Lambda_{\F,N}=\bigoplus_{n\geq0}\, \Lambda_{\F,N}^n$, where $\Lambda_{\F,N}^n$ is the ring of homogeneous symmetric polynomials in $N$ variables of degree $n$ (together with the zero polynomial).

For any $M>N$ there are morphisms
$\rho_{MN}\colon\Lambda_{\F,M}\to \Lambda_{\F,N}$ that map the
variables $x_{N+1},\ldots,x_M$ to zero. They preserve the grading, and
hence we can define $\rho_{MN}^n\colon\Lambda^n_{\F,M}\to
\Lambda^n_{\F,N}$; this allows us to define the inverse limits
\begin{equation}
\Lambda^n_\F:= \lim_{\stackrel{\scriptstyle\longleftarrow}{\scriptstyle
    N}} \, \Lambda^n_{\F,N}\ ,
\end{equation}
and the \emph{algebra of symmetric functions in infinitely many variables} as $\Lambda_\F:=\bigoplus_{n\geq0}\, \Lambda^n_\F.$ In the following when no confusion is possible we will denote $\Lambda_\F$ (resp.\ $\Lambda_\F^n$) simply by $\Lambda$ (resp. $\Lambda^n$). 

Now we introduce a basis for $\Lambda.$ For this, we start by defining a basis in $\Lambda_N.$ Let $\mu=(\mu_1,\ldots,\mu_t)$ be a partition with $t\leq N$, and define the polynomial
\begin{equation}
m_\mu(x_1,\ldots,x_N) =\sum_{\tau\in \sigma_N}\, x_1^{\mu_{\tau(1)}}\cdots x_N^{\mu_{\tau(N)}}\ ,
\end{equation}
where we set $\mu_j=0$ for $j=t+1,\ldots, N$. The polynomial $m_\mu$ is symmetric, and the set of $m_\mu$ for all partitions $\mu$ with $\vert \mu \vert\leq N$ is a basis of $\Lambda_N.$ Then the set of $m_\mu$, for all partitions $\mu$ with $\vert \mu \vert\leq N$ and $\sum_i\, \mu_i=n$, is a basis of $\Lambda^n_N$. Since for $M>N\geq t$ we have $\rho^n_{MN}(m_\mu(x_1,\ldots,x_M))=m_\mu(x_1,\ldots,x_N)$, by using the definition of inverse limit we can define the \emph{monomial symmetric functions} $m_\mu.$ By varying over the partitions $\mu$ of $n$, these functions form a basis for $\Lambda^n$.

Next we define the $n$\emph{-th power sum symmetric function} $p_n$ as
\begin{equation}
p_n := m_{(n)}= \sum_{i} \, x_i^n\ .
\end{equation}
The set consisting of symmetric functions $p_\mu:=p_{\mu_1}\dots p_{\mu_t}$, for all partitions $\mu=(\mu_1,\ldots,\mu_t)$, is another basis of $\Lambda$. 

We now set $\F=\C$ throughout and we fix a parameter $\beta\in \C$ (though
everything works for any field extension $\C\subseteq \F$ and $\beta\in
\F$). Define an inner product $\langle-,-\rangle_{\beta}$ on the vector space
$\Lambda\otimes\Q(\beta)$ with respect to which the basis of power sum symmetric functions $p_\lambda(x)$ are orthogonal with the normalization
\begin{equation}\label{eq:jackinnerprod}
\langle p_\lambda,p_\mu\rangle_{\beta} = \delta_{\lambda,\mu}\ z_\lambda\, \beta^{-\ell(\lambda)} \ ,
\end{equation}
where $\delta_{\lambda,\mu}:=\prod_{a}\, \delta_{\lambda_a,\mu_a}$ and
\begin{equation}
z_\lambda:= \prod_{j\geq1}\, j^{m_j}\, m_j! \ .
\end{equation}
This is called the \emph{Jack inner product}.

\begin{definition}\label{def:macdonald}
The monic forms of the {\em Jack functions} $J_\lambda( x;\beta^{-1} )\in\Lambda\otimes\Q(\beta)$ for $x=(x_1, x_2, \ldots)$ are uniquely defined by the following two conditions \cite{book:macdonald1995}:
\begin{itemize}
\item[(i)] Triangular expansion in the basis $m_\mu(x)$ of monomial symmetric functions:
\begin{equation}
J_\lambda(x;\beta^{-1}) = m_\lambda (x) + \sum_{\mu < \lambda} \,
\psi_{\lambda, \mu} (\beta) \, m_\mu (x) \qquad\mbox{with} \quad \psi_{\lambda, \mu} (\beta) \in \Q(\beta)\ .
\end{equation}
\item[(ii)] Orthogonality:
\begin{equation}\label{eq:orthogonalityjack}
\langle J_\lambda, J_\mu\rangle_{\beta} = \delta_{\lambda,\mu} \ \prod_{s\in Y_\lambda} \ \frac{\beta\,L(s) + A(s)+1}{\beta\,\big(L(s)+1 \big) + A(s)} \ .
\end{equation}
\end{itemize}
\end{definition}

\begin{lemma}\label{lem:p1}
For any integer $n\geq 1$ we have
\begin{equation}\label{eq:p1}
(p_1)^n=n!\, \sum_{\vert\lambda\vert=n} \ \prod_{s\in Y_\lambda} \ \frac{1}{\beta\,L(s) + A(s)+1}\ J_\lambda\ .
\end{equation}
\end{lemma}
\proof
The assertion follows straightforwardly from \cite[Proposition 2.3 and
Theorem 5.8]{art:stanley1989} after normalizing our Jack functions:
the Jack functions considered in \cite{art:stanley1989} are given by
\begin{equation}
\tilde J_\lambda=\beta^{-\vert\lambda\vert}\ \prod_{s\in Y_\lambda}\
\big[\beta\,\big(L(s)+1 \big) + A(s)\big] \ J_\lambda\ ,
\end{equation}
where the normalization factor is computed by using \cite[Theorem 5.6]{art:stanley1989}.
\endproof

\bigskip \section{Infinite-dimensional Lie algebras\label{sec:infiniteliealg}}

\subsection{Heisenberg algebras}

In this section we recall the representation theory of Heisenberg
algebras and the affine Lie algebras $\slfrakhat_k$. Since the Lie algebra $\glfrak_k$ coincides with $\F\,\mathrm{id}\oplus\slfrak_k$, as a by-product we get the representation theory of $\glfrakhat_k$.

Let $\C\subseteq \F$ be a field extension of $\C$. Let $\Lfrak$ be a
{lattice}, i.e., a free abelian group of finite rank $d$ equipped with a symmetric nondegenerate bilinear form $\langle -,- \rangle_{\Lfrak} \colon \Lfrak\times \Lfrak \to \Z$. Fix a basis $\gamma_1, \ldots, \gamma_d$ of $\Lfrak$.
\begin{definition}
The \emph{lattice Heisenberg algebra} $\hfrak_{\F,\Lfrak}$ associated with $\Lfrak$ is the infinite-dimensional Lie algebra over $\F$ generated by $\mathfrak{q}^i_m $, for $m\in\mathbb{Z}\setminus\{0\}$ and $i\in\{1, \ldots, d\}$, and the {central element} $\mathfrak{c}$ satisfying the relations
\begin{equation}\label{eq:commutationrelations}
\left\{
\begin{array}{ll}
\left[\mathfrak{q}^i_m,\cfrak\right]=0 & \mbox{for} \quad m\in
\Z\setminus\{0\} \ ,\ i\in\{1, \ldots, d\}\ , \\[6pt]
\big[\mathfrak{q}^i_m,\mathfrak{q}^j_n\big]=m\, \delta_{m,-n} \,
\langle \gamma_i, \gamma_j \rangle_{\Lfrak}\ \cfrak & \mbox{for}
\quad m,n\in \Z\setminus\{0\} \ ,\ i,j\in\{1, \ldots, d\}\ .
\end{array}
\right.
\end{equation}
\end{definition}
For any element $v\in \Lfrak$ we define the element
$\mathfrak{q}^v_m\in\hfrak_{\F,\Lfrak}$ by linearity, with $\mathfrak{q}^i_m :=\mathfrak{q}_m^{\gamma_i}$. Set
\begin{equation}\label{eq:triangularheis}
\hfrak_{\F,\Lfrak}^+:=\bigoplus_{m>0} \ \bigoplus_{i=1}^d \,  \F\mathfrak{q}^i_m\qquad
\mbox{and} \qquad
\hfrak_{\F,\Lfrak}^-:=\bigoplus_{m<0} \ \bigoplus_{i=1}^d \, \F\mathfrak{q}^i_m \ .
\end{equation}

Let us denote by $\Ucal(\hfrak_{\F, \Lfrak})$ (resp.\
$\Ucal(\hfrak^\pm_{\F, \Lfrak})$) the {universal enveloping
  algebra} of $\hfrak_{\F, \Lfrak}$ (resp.\ $\hfrak_{\F,\Lfrak}^\pm$),
i.e., the unital associative algebra over $\F$ generated by
$\hfrak_{\F, \Lfrak}$ (resp.\ $\hfrak_{\F,\Lfrak}^\pm$).

We introduce some terminology similar to that used in \cite[Section 5.2.5]{book:frenkelben-zvi2004}.
\begin{definition}
For $v\in\Lfrak$, define \emph{free bosonic fields} as the elements
\begin{equation}
\varphi_-^v(z):=\sum_{m=1}^\infty\,\frac{z^m}{m}\,
\qfrak^v_{-m}\qquad\mbox{and}\qquad\varphi_+^v(z):=
\sum_{m=1}^\infty\, \frac{z^{-m}}{m}\, \qfrak^v_{m}
\end{equation}
in $\hfrak_{\F,\Lfrak}^-[[z]]$ and $\hfrak_{\F,\Lfrak}^+[[z^{-1}]]$,
respectively. For $\alpha,\beta\in\F$, define the \emph{generalized
  bosonic exponential}
\begin{equation}
\V_{\alpha,\beta}^v(z):=\exp\big(\alpha\,\varphi_-^v(z)\big)\, \exp\big(\beta\,\varphi_+^v(z)\big) =: \
:\, \exp\big(\alpha\,\varphi_-(z)+\beta\, \varphi_+(z)\big)\, :
\end{equation}
in $\hfrak_{\F,\Lfrak}[[z,z^{-1}]]$, where the symbol $:-:$ denotes
normal ordering with respect to the decomposition $\hfrak_{\F,\Lfrak}=
\hfrak_{\F,\Lfrak}^-\oplus \hfrak_{\F,\Lfrak}^+$, i.e., all negative
generators $\qfrak_{-m}^v$ are moved to the left of all positive
generators $\qfrak_m^v$ for $m>0$. When $\beta=-\alpha$, we call
$\V_{\alpha,-\alpha}^v(z)$ a \emph{normal-ordered bosonic exponential}.
\label{def:bosonops}\end{definition}

\begin{remark}
The bosonic exponentials are \emph{vertex operators}, i.e., they are
uniquely characterized by their commutation relations in the
Heisenberg algebra $\hfrak_{\F,\Lfrak}$: For $v,v'\in\Lfrak$ one has
\begin{equation}
\big[\qfrak_m^v,\V_{\alpha,\beta}^{v'}(z) \big] = \left\{
\begin{array}{ll}
\alpha\, \langle v,v'\,\rangle_\Lfrak\, z^m\,\V_{\alpha,\beta}^{v'}(z)
& \mbox{for} \quad m>0 \ , \\[6pt]
-\beta\, \langle v,v'\,\rangle_\Lfrak\, z^m\,\V_{\alpha,\beta}^{v'}(z)
 & \mbox{for} \quad m<0 \ .
\end{array} \right.
\end{equation}
The compositions of vertex operators
$\V_{\alpha_1,\beta_1}^{v_1}(z_1)\cdots\V_{\alpha_n,\beta_n}^{v_n}(z_n)$
in $\hfrak_{\F,\Lfrak}[[z_1^{\pm\,1},\dots,z_n^{\pm\,1}]]$ can be
easily calculated as
\begin{equation}\label{eq:VOnprod}
\prod_{i=1}^n\, \V_{\alpha_i,\beta_i}^{v_i}(z_i) = \bigg(\,
\prod_{1\leq j<i\leq n}\, \Big(\,1-\frac{z_i}{z_j}\,
\Big)^{-\alpha_i\, \beta_j\, \langle v_i,v_j\rangle_\Lfrak}\, \bigg) \
:\, \prod_{i=1}^n\, \V_{\alpha_i,\beta_i}^{v_i}(z_i)\, : \ ,
\end{equation}
where the factors $(1-\frac{z_i}{z_j}\,
)^{-\alpha_i\, \beta_j\, \langle v_i,v_j\rangle_\Lfrak}$ are
understood as formal power series in $\frac{z_i}{z_j}$.
\end{remark}

\begin{remark}
When $v=\gamma_i$ for $i=1, \ldots, d$, we simply denote
$\varphi_\pm^i(z):= \varphi_\pm^{\gamma_i}(z)$; if $d=1$, we further
simply write $\varphi_\pm(z)$. We use analogous notation for the generalized free boson exponentials.
\end{remark}

\begin{example}\label{ex:rankkheis}
Consider the lattice $\Lfrak:=\Z^k$ with the symmetric nondegenerate
bilinear form $\langle v, w \rangle_\Lfrak = \sum_{i=1}^k \,
v_i\,w_i$. In this case $\hfrak_{\F,\Lfrak}$ is called the {\em
  Heisenberg algebra of rank $k$} over $\F$, and we denote it by $\hfrak_\F^k$. It is generated by elements $\pfrak^i_m$, $m\in\Z\setminus \{0\}$, $i=1,\ldots,k$, and the central element $\cfrak$ satisfying the relations \eqref{eq:commutationrelations} with $\langle \gamma_i, \gamma_j \rangle_{\Lfrak}=\delta_{ij}$. When $k=1$, $\hfrak_{\F,\Lfrak}$ is simply the infinite-dimensional \emph{Heisenberg algebra} $\hfrak_\F$ over the field $\F$.
\end{example}

\begin{example}\label{ex:latticeheisenebrg}
Fix an integer $k\geq 2$ and let $\Qfrak$ be the {root lattice} of
type $A_{k-1}$ endowed with the standard bilinear form $\langle
-,-\rangle_{\Qfrak}$ (see Remark \ref{rmk:rootlattice} below). Let
$\hfrak_{\F,\Qfrak}$ be the lattice Heisenberg algebra over $\F$
associated to $\Qfrak$; we call $\hfrak_{\F,\Qfrak}$ the \emph{Heisenberg algebra of type} $A_{k-1}$ over $\F$. It can be realized  as the Lie algebra over $\F$ generated by $\qfrak_m^i$ for $m\in\Z\setminus \{ 0 \}$, $i=1,\ldots, k-1$, and the central element $\cfrak$ satisfying the relations 
\begin{equation}
\left\{
\begin{array}{ll}
\left[ \qfrak_m^i, \cfrak \right] = 0 \quad & \mbox{for}\quad
m\in\Z\setminus \{ 0 \} \ , \ i=1,\ldots, k-1\ ,\\[6pt]
\big[ \qfrak_m^i, \qfrak_n^j \big] = m\ \delta_{m,-n}\ C_{ij}\ \cfrak
\quad & \mbox{for}\quad m\in\Z\setminus \{ 0 \} \ , \ i,j=1,\ldots, k-1\ ,
\end{array}
\right.
\end{equation}
where $C=(C_{ij})$ is the Cartan matrix of type $A_{k-1}$.
\end{example}

\subsubsection{Virasoro generators}\label{sec:virasoro-heisenberg}

We construct the Viraroso algebra associated with the Heisenberg algebra $\hfrak_{\F}$. Define elements
\begin{align}
L^\hfrak_0=\sum_{m=1}^\infty\, \qfrak_{-m}\ \qfrak_{m} \qquad \mbox{and} \qquad
L^\hfrak_n = \frac{1}{2}\ \sum_{m\in\Z}\, \qfrak_{-m}\ \qfrak_{m+n}\quad
\mbox{for}\quad n\in\Z\setminus\{0\} 
\end{align}
in the completion of the enveloping algebra $\Ucal(\hfrak_{\F})$, where we set $\qfrak_0:=0$. They satisfy the relations
\begin{equation}
\big[L^\hfrak_n, L^\hfrak_m \big]=(n-m)\, L^\hfrak_{n+m}+\frac{n}{12}\,\big(n^2-1 \big)\, \delta_{m+n,0}\,\cfrak\ ,
\end{equation}
hence $\cfrak$ and $L^\hfrak_n$ with $n\in\Z$ generate a Virasoro
algebra $\virfrak_\F$ over $\F$.

\begin{remark}\label{rem:genbosprimary}
It is well-known (see Appendix~\ref{app:virasoroprimary}) that the
generalized bosonic exponential $\V_{\alpha,\beta}(z)$ is a primary
field of the Virasoro algebra $\virfrak_\F$ generated by $L_n^\hfrak$
with conformal dimension $\Delta(\alpha,\beta)=-\frac12\,
\alpha\,\beta$, i.e., it satisfies the commutation relations
\begin{equation}
\big[L_n^\hfrak,\V_{\alpha,\beta}(z)\big] = z^n\,
\big(z\,\partial_z+\Delta(\alpha,\beta)\, (n+1) \big)\, \V_{\alpha,\beta}(z) \ .
\end{equation}
\end{remark}

\subsubsection{Fock space}

We are interested in a special type of representation of a given
lattice Heisenberg algebra $\hfrak_{\F,\Lfrak}$ over $\F$.
\begin{definition}\label{def:fock}
Let $W$ be the trivial representation of $\hfrak_{\F,\Lfrak}^{+}$
(i.e., the one-dimensional $\F$-vector space with trivial
$\hfrak_{\F,\Lfrak}^{+}$-action). The \emph{Fock space} representation
of the Heisenberg algebra $\hfrak_{\F,\Lfrak}$ is the induced representation $\mathcal{F}_{\F,\Lfrak}:=\hfrak_{\F,\Lfrak}\otimes_{\hfrak_{\F,\Lfrak}^+} W$.
\end{definition}
The Fock space is an irreducible \emph{highest weight representation} whereby any element $w_0\in W$ is a \emph{highest weight vector}, i.e., $\hfrak_{\F,\Lfrak}^{+}$ annihilates $w_0$ and the elements in $W$ of the form $\qfrak^v_{-m_1}\cdots \qfrak^v_{-m_l}\triangleright w_0$ generate $\mathcal{F}_{\F,\Lfrak}$ for $v\in \Lfrak$, $l\geq 1$ and $m_i\geq 1$ for $i=1, \ldots, l$.

\begin{example}
For the Heisenberg algebra $\hfrak_\F$, the Fock space
$\mathcal{F}_\F$ is isomorphic to the polynomial algebra
$\Lambda_\F=\F\left[p_1, p_2, \ldots \right]$ in the power sum
symmetric functions introduced in Section \ref{sec:symmetric}. In this realization, the actions of the generators are given for $m>0$ by
\begin{equation}\label{eq:fockheisenberg}
\pfrak_{-m}\triangleright f:=p_mf \ , \qquad \pfrak_m\triangleright
f:=m\, \frac{\partial f}{\partial p_m} \qquad \mbox{and} \qquad \cfrak\triangleright f:=f 
\end{equation}
for any $f\in \Lambda_{\F}$.
\end{example}

\begin{example}
The Fock space $\Fcal^k_\F$ of the rank $k$ Heisenberg algebra $\hfrak_\F^k$ can be realized as the tensor product of $k$ copies of the polynomial algebra $\Lambda_\F$:
\begin{equation}
\Fcal^k_\F \simeq \Lambda_\F^{\otimes k}\ .
\end{equation}
In this realization, the action of the generators $\pfrak_m^i$ is obvious: each copy of the Heisenberg algebra generated by $\pfrak_m^i$ for $m\in\Z\setminus\{0\}$ acts on the $i$-th factor $\Lambda_\F$ as in Equation \eqref{eq:fockheisenberg}.
\end{example}

\subsubsection{Whittaker vectors}

We give the definition of Whittaker vector for Heisenberg algebras
following \cite[Section 3]{art:christodoulopoulou2008}; in conformal
field theory it has the meaning of a \emph{coherent state}.
\begin{definition}\label{def:whittakerheisenberg}
Let $\chi\colon \Ucal(\hfrak_{\F, \Lfrak}^+)\to \F$ be an algebra homomorphism such that $\chi\vert_{\hfrak_{\F, \Lfrak}^+}\neq 0$, and let $V$ be a $\Ucal(\hfrak_{\F, \Lfrak})$-module. A nonzero vector $w\in V$ is called a \emph{Whittaker vector of type} $\chi$ if $\eta\triangleright w=\chi(\eta)\, w$ for all $\eta\in \Ucal(\hfrak_{\F, \Lfrak}^+)$.
\end{definition}
\begin{remark}\label{rmk:unicitawhittakerheis}
By \cite[Proposition~10]{art:christodoulopoulou2008}, if $w,w'$
are Whittaker vectors of the same type $\chi$, then $w'=\lambda \, w$ for some nonzero $\lambda\in \F$.
\end{remark}

\subsection{Affine algebra of type $\widehat{A}_{k-1}$}

Let $k\geq 2$ be an integer and let $\slfrak_k:=\slfrak(k,\F)$ denote the finite-dimensional Lie algebra of rank
$k-1$ over $\F$
generated in the Chevalley basis by $E_i, F_i, H_i$ for $i=1, \ldots, k-1$ satisfying the relations
\begin{equation}
\begin{array}{ll}
\left[E_i, F_j\right] = \delta_{ij} \, H_j\ , & [H_i, H_j]=0\ ,\\[4pt]
\left[H_i, E_j\right] = C_{ij} \, E_j\ , & [H_i, F_j]=-C_{ij}\, F_j\ ,
\end{array}
\end{equation}
where $C=(C_{ij})$ is the Cartan matrix type $A_{k-1}$ (see Remark \ref{rmk:rootlattice} below).

An explicit realization of the generators of $\slfrak_k$ in the algebra
$M(k,\F)$ of $k\times k$ matrices over $\F$ is given in the following way. Let
$\boldsymbol{E}_{i,j}$ denote the $k\times k$ matrix unit with 1 in
the $(i,j)$ entry and 0 everywhere else for $i,j=1, \ldots, k$. Define 
\begin{equation}
E_i:=\boldsymbol{E}_{i,i+1}\ ,\qquad F_i:= \boldsymbol{E}_{i+1,i}
\qquad \mbox{and} \qquad H_i:= \boldsymbol{E}_{i,i}- \boldsymbol{E}_{i+1,i+1}
\end{equation}
for $i=1, \ldots, k-1$. One sees immediately that $E_i, F_i, H_i$
satisfy the defining relations for $\slfrak_k$.

Let us denote by $\mathfrak{t}$ the Lie subalgebra of $\slfrak_k$
generated by $H_i$ for $i=1, \ldots, k-1$ and by $\mathfrak{n}_+$
(resp.\ $\mathfrak{n}_{-}$) the Lie subalgebra of $\slfrak_k$
generated by $E_i$ (resp.\ $F_i$) for $i=1, \ldots, k-1$. Then there
is a {triangular decomposition}
\begin{equation}
\slfrak_k=\mathfrak{n}_{-}\oplus\mathfrak{t}\oplus\mathfrak{n}_+
\end{equation}
as a direct sum of vector spaces.

\begin{remark}\label{rmk:rootlattice}

For $i=1, \ldots, k$, define $\boldsymbol{e}_i\in \mathfrak{t}^\ast$ by
\begin{equation}
\boldsymbol{e}_i\big(\mathrm{diag}(a_1, \ldots, a_k) \big)=a_i \ .
\end{equation}
The elements $\gamma_i:=\boldsymbol{e}_i-\boldsymbol{e}_{i+1}$ for
$i=1, \ldots, k-1$ form a basis of $\mathfrak{t}^\ast$. The
\emph{root lattice} $\Qfrak$ is the lattice
$\Qfrak:=\bigoplus_{i=1}^{k-1}\, \Z\gamma_i$. The elements of
$\Qfrak$ are called \emph{roots}, and in particular $\gamma_i$ are called the
\emph{simple roots}. The lattice of positive roots is
$\Qfrak_+:=\bigoplus_{i=1}^{k-1} \, \N \gamma_i$. Since $\boldsymbol{e}_i$ corresponds to the $i$-th
coordinate vector in $\Z^{k}$, there is a description of $\Qfrak$ and
$\Qfrak_+$ in $\Z^k$ given by
\begin{equation}
\Qfrak=\big\{\boldsymbol{e}_{i}-\boldsymbol{e}_{j}\, \big\vert\, i,j=1, \ldots,
k \big\}\qquad \mbox{and} \qquad
\Qfrak_+=\big\{\boldsymbol{e}_{i}-\boldsymbol{e}_{j}\, \big\vert\,
1\leq i<j\leq k \big\}\ .
\end{equation}
By setting $\langle \gamma_i, \gamma_j \rangle_\Qfrak:=\gamma_i(H_j)=C_{ij}$, we define a nondegenerate symmetric bilinear form $\langle-,- \rangle_\Qfrak$ on $\Qfrak$. 

The \emph{fundamental weights $\omega_i$ of type $A_{k-1}$} are the elements of $\mathfrak{t}^\ast$ defined by $\omega_i(H_j)=\delta_{ij}$ for $i,j=1, \ldots, k-1$. In the standard basis of $\Z^k$, they are given explicitly by
\begin{equation}
\omega_i:=\sum_{l=1}^{i} \, \boldsymbol{e}_l-\frac{i}{k}\,
\sum_{l=1}^{k}\, \boldsymbol{e}_l
\end{equation}
for $i=1,\dots,k-1$. Let $\Pfrak:=\bigoplus_{i=1}^{k-1}\, \Z \omega_i$
be the \emph{weight lattice}. Then $\Qfrak\subset \Pfrak$, as
$\gamma_i=\sum_{j=1}^{k-1} \, C_{ij}\, \omega_j$. The set of
\emph{dominant weights} is $\Pfrak_+:=\bigoplus_{i=1}^{k-1}\, \N
\omega_i$. There is a coset decomposition of $\Pfrak$ given by
\begin{equation}
\Pfrak=\bigcup_{j=0}^{k-1}\, (\Qfrak+\omega_j)\ ,
\label{eq:lateralweightlattice}\end{equation}
where we set $\omega_0:=\boldsymbol{0}$.

The \emph{coroot lattice} is the lattice $\Qfrak^\vee:=\bigoplus_{i=1}^{k-1}\, \Z H_i$.
\end{remark}

We now introduce the Kac-Moody algebra $\slfrakhat_k$ of type $\widehat{A}_{k-1}$, first via its {canonical generators} and then as a central extension of the {loop algebra} of $\slfrak_k$.
\begin{definition}
The \emph{Kac-Moody algebra} $\slfrakhat_k$ \emph{of type} $\widehat{A}_{k-1}$ over $\F$ is the Lie algebra over $\F$ generated by $e_i, f_i, h_i$ for $i=0,1, \ldots, k-1$ satisfying the relations
\begin{equation}
\begin{array}{ll}
\left[e_i, f_j\right]= \delta_{ij} \, h_j\ ,& [h_i, h_j]=0\ ,\\[4pt]
\left[h_i, e_j\right] = \widehat{C}_{ij} \, e_j\ ,& [h_i,
f_j]=-\widehat{C}_{ij} \, f_j\ ,
\end{array}
\end{equation}
where $\widehat{C}=\big(\widehat{C}_{ij} \big)$ is the Cartan matrix of the extended Dynkin diagram of type $\widehat{A}_{k-1}$.
\end{definition}
The matrix $\widehat{C}$ is given for $k\geq 3$ by
\begin{equation}
\widehat{C}=\big(\widehat{C}_{ij} \big)=\begin{pmatrix}
2 & -1 & 0 & \dots & -1 \\
-1 & 2& -1&   \dots & 0 \\
0 & -1& 2& \dots & 0 \\
\vdots & \vdots & \vdots & \ddots & \vdots \\
-1 &0 &0 & \dots & 2
\end{pmatrix}
\end{equation}
and for $k=2$ by
\begin{equation}
\widehat{C}=\big(\widehat{C}_{ij} \big)=\begin{pmatrix}
2 & -2 \\
-2 & 2
\end{pmatrix}.
\end{equation}

Let us denote by $\tfrakhat$ the Lie subalgebra of $\slfrakhat_k$
generated by $h_i$ for $i=0, 1, \ldots, k-1$ and by $\nfrakhat_+$
(resp.\ $\nfrakhat_{-}$) the Lie subalgebra of $\slfrakhat_k$
generated by $e_i$ (resp.\ $f_i$) for $i=0, 1, \ldots, k-1$. Then
there is a {triangular decomposition}
\begin{equation}\label{eq:triangularslk}
\slfrakhat_k=\nfrakhat_{-}\oplus\tfrakhat\oplus\nfrakhat_+
\end{equation}
as a direct sum of vector spaces.

Now we describe the relation between $\slfrak_k$ and
$\slfrakhat_k$. Define in $\slfrak_k$ the elements
\begin{equation}
E_0:=\boldsymbol{E}_{k,1}\ ,\qquad F_0:= \boldsymbol{E}_{1,k}\qquad \mbox{and} \qquad
H_0:=\boldsymbol{E}_{k,k}- \boldsymbol{E}_{1,1} \ .
\end{equation}
Consider next the \emph{loop algebra} $\widetilde{\slfrak}_k:=\slfrak_k\otimes \F[t,t^{-1}]$. Set 
\begin{equation}
\begin{array}{ll}
\tilde{e}_0:=E_0\otimes t\ , & \tilde{e}_i:=E_i\otimes 1\ ,\\[4pt]
\tilde{f}_0:=F_0\otimes t^{-1}\ , & \tilde{f}_i:=F_i\otimes 1\ , \\[4pt]
\tilde{h}_0:= H_0\otimes 1\ , & \tilde{h}_i:=H_i\otimes 1\ ,
\end{array}
\end{equation}
for $i=1, \ldots, k-1$. Let us denote by $\cfrak$ the {central element} of $\slfrakhat_k$
given by $\cfrak=\sum_{i=0}^{k-1} \, h_i$. Then we can realize $\slfrakhat_k$ as a one-dimensional central extension
\begin{equation}
0\ \longrightarrow \ \F \cfrak\ \longrightarrow \ \slfrakhat_k \
\xrightarrow{ \ \pi \ } \ \widetilde{\slfrak}_k \ \longrightarrow \ 0\ ,
\end{equation}
where the homomorphism $\pi$ is defined by
\begin{equation}
\pi\, :\, e_i\ \longmapsto \tilde{e}_i\ ,\qquad f_i \ \longmapsto\
\tilde{f}_i\ ,\qquad h_i\ \longmapsto \ \tilde{h}_i\ ,
\end{equation}
for $i=0, 1, \ldots, k-1$, and the Lie algebra structure of $\slfrakhat_k$ is obtained through
\begin{equation}\label{eq:algebrastructslk}
[M\otimes t^m, N\otimes t^n]=[M, N]\otimes t^{m+n} +m \, \delta_{m,
  -n}\,\mathrm{tr}(M\, N) \ \cfrak
\end{equation}
for every $M, N\in \slfrak_k$ and $m,n\in\mathbb{Z}$. Thus the canonical generators of $\slfrakhat_k$ are
\begin{equation}
\begin{array}{ll}
e_0:= E_0\otimes t\ ,& e_i:=E_i\otimes 1 \ ,\\[4pt]
f_0:= F_0\otimes t^{-1}\ , & f_i:= F_i\otimes 1\ ,\\[4pt]
h_0:= H_0\otimes 1 +\cfrak\ ,& h_i:=H_i\otimes 1 \ ,
\end{array}
\end{equation}
and we can realize $\tfrakhat$ as the one-dimensional extension
\begin{equation}
0\ \longrightarrow \ \F \cfrak\ \longrightarrow \ \tfrakhat \
\xrightarrow{ \ \pi \ } \ \mathfrak{t} \ \longrightarrow \ 0\ . 
\end{equation}

\begin{remark}
Let $\gamma_0:=-\sum_{i=1}^{k-1} \gamma_i$. For $i=1,
\ldots, k-1$, let $\boldsymbol{e}_i$ be as in Remark
\ref{rmk:rootlattice}; then $\gamma_0=\boldsymbol{e}_k-\boldsymbol{e}_1$. We extend $\boldsymbol{e}_i$ from
$\tfrak^\ast$ to $\tfrakhat^\ast$ by setting
$\boldsymbol{e}_i(\cfrak)=0$. Then $\gamma_i(\cfrak)=0$ for $i=0,1, \ldots, k-1$. Thus the root lattice $\widehat{\Qfrak}$ of
$\slfrakhat_k$ is the lattice $\widehat\Qfrak= \bigoplus_{i=0}^{k-1} \, \Z\gamma_i = \Z\gamma_0 \oplus \Qfrak$. In a similar way, one can define the lattice of positive roots and a nondegenerate symmetric bilinear form on $\widehat{\Qfrak}$.

Let $\widehat{\omega}_0$ be the element in $\tfrakhat^\ast$ defined by $\widehat{\omega}_0(\tfrak^\ast)=0$ and $\widehat{\omega}_0(\cfrak)=1$. Define 
\begin{equation}
\widehat{\omega}_i:=\omega_i+\widehat{\omega}_0\qquad\mbox{for} \quad i=1, \ldots, k-1\ .
\end{equation}
We call $\widehat{\omega}_0, \widehat{\omega}_1, \ldots,
\widehat{\omega}_{k-1}$ the \emph{fundamental weights of type
  $\widehat{A}_{k-1}$}. Set $\Pfrakhat:=\bigoplus_{i=0}^{k-1}\, \Z
\widehat{\omega}_i$. Any weight $\widehat{\lambda}=\sum_{i=0}^{k-1}\,
\lambda_i \, \widehat{\omega}_i\in \Pfrakhat$ can be written as
$\widehat{\lambda}=\lambda+k_{\,\widehat{\lambda}}\,
\widehat{\omega}_0$, where $\lambda\in \Pfrak$ and
$k_{\,\widehat{\lambda}}=\widehat{\lambda}(\cfrak)=\sum_{i=0}^{k-1} \,
\lambda_i$ is the \emph{level} of $\widehat{\lambda}$.
\end{remark}

\subsubsection{Highest weight representations}

By declaring the degrees of generators $\deg e_i= -\deg f_i =1$ and $\deg h_i=0$ for $i=0,1, \ldots, k-1$, we endow $\slfrakhat_k$ with the \emph{principal grading}
\begin{equation}
\slfrakhat_k= \bigoplus_{n\in \Z} \, \big(\, \slfrakhat_k \big)_n \ .
\end{equation}
The principal grading of
$\slfrakhat_k$ induces a $\Z$-grading of its universal enveloping
algebra $\Ucal\big(\, \slfrakhat_k \big)$ over $\F$, which is written
as
\begin{equation}
\Ucal\big(\, \slfrakhat_k \big)=\bigoplus_{n\in \Z} \, \Ucal_n\ .
\end{equation}
Set $\bfrakhat:=\tfrakhat\oplus\nfrakhat_+$. Let $\widehat{\lambda}$ be a linear
form on $\tfrakhat$. We define a one-dimensional $\bfrakhat$-module
$\F v_{\, \widehat{\lambda}}$ by 
\begin{equation}
\nfrakhat_+\triangleright v_{\, \widehat{\lambda}}=0 \qquad\mbox{and}\qquad
h_i\triangleright v_{\, \widehat{\lambda}}=\widehat{\lambda}(h_i)\, v_{\, \widehat{\lambda}} \qquad\mbox{for}\quad i=0,1, \ldots, k-1\ .
\end{equation}
Consider the induced $\slfrakhat_k$-module
\begin{equation}
\tilde{\Vcal}(\,\widehat{\lambda}\, ):=\Ucal\big(\,\slfrakhat_k \big)\otimes_{\Ucal(\,\bfrakhat\,
  )} \F v_{\,\widehat{\lambda}}\ .
\end{equation}
Setting $\tilde{\Vcal}_n:=\Ucal_n\triangleright v_{\,\widehat{\lambda}}$, we define
the principal grading $\tilde{\Vcal}(\,\widehat{\lambda}\,
)=\bigoplus_{n\in \Z}\, \tilde{\Vcal}_n$. The $\slfrakhat_k$-module
$\tilde{\Vcal}(\,\widehat{\lambda}\, )$ contains a unique maximal
proper graded $\slfrakhat_k$-submodule $I(\,\widehat{\lambda}\, )$.

\begin{definition}\label{def:basicrepslk}
The quotient module
\begin{equation}
\Vcal(\,\widehat{\lambda}\,):=\tilde{\Vcal}(\,\widehat{\lambda}\,
)\,\big/\, I(\,\widehat{\lambda}\, )
\end{equation}
is called the \emph{highest weight representation of $\slfrakhat_k$ at level $k_{\,\widehat{\lambda}}$}. The nonzero multiples of the image of $v_{\,\widehat{\lambda}}$ in
$\Vcal(\,\widehat{\lambda}\, )$ are called the \emph{highest weight vectors} of
$\Vcal(\,\widehat{\lambda}\, )$.
\end{definition}

The principal grading on $\tilde{\Vcal}(\,\widehat{\lambda}\, )$ induces an $\N$-grading
\begin{equation}
\Vcal(\, \widehat{\lambda}\, )=\bigoplus_{n\geq 0}\, \Vcal_{-n}
\end{equation}
called the principal grading of $\Vcal(\,\widehat{\lambda}\, )$.

\begin{definition}
The \emph{$i$-th dominant representation of $\slfrakhat_k$ at level
  one} is the highest weight representation
$\Vcal(\,\widehat{\omega}_i)$ of $\slfrakhat_k$ for $i=0,1, \ldots,
k-1$. The module $\Vcal(\,\widehat{\omega}_0)$ is also called the \emph{basic representation of $\slfrakhat_k$}.\\
\end{definition}

\begin{remark}
One can define the Lie algebra $\glfrakhat_k$ as the one-dimensional extension
\begin{equation}
0\ \longrightarrow \ \F \cfrak\ \longrightarrow \ \glfrakhat_k \
\xrightarrow{ \ \pi \ } \ \mathfrak{gl}_k\otimes\F[t,t^{-1}] \ \longrightarrow \ 0\ . 
\end{equation}
Since $\mathfrak{gl}_k=\F\,\mathrm{id}\oplus\mathfrak{sl}_k$, the
representation theory of  $\glfrakhat_k$ is obtained by combining the
representation theory of the Heisenberg algebra $\hfrak_\F$ with that of
$\slfrakhat_k$. For example, all highest weight representations of
$\glfrakhat_k$ are of the form $\Fcal_{\F}\otimes \Vcal(\, \widehat
\lambda\, )$ for some weight $\widehat \lambda\in \widehat \Pfrak$.
\end{remark}

\subsubsection{Whittaker vectors}

Let us denote by $\qfrak_m^i$ the element $H_i\otimes t^m$ for $i\in\{1, \ldots, k-1\}$ and $m\in\Z$. By Equation \eqref{eq:algebrastructslk}, these elements satisfy
\begin{equation}
\big[\qfrak_m^i, \qfrak_n^j\big]=m\ \delta_{m+n, 0}\ C_{ij}\ \cfrak \qquad\mbox{and}\qquad \big[\qfrak_m,\cfrak\big]=0\ ,
\end{equation}
for $i,j\in\{1, \ldots, k-1\}$ and $m,n\in\Z$. For a root $\gamma$, we
denote by $\qfrak_m^\gamma$ the element $H_\gamma\otimes t^m$ where
$H_\gamma\in\tfrak$ is defined by $\langle H,
H_\gamma\rangle_{\Qfrak^\vee\otimes_\Z \R} = \gamma(H)$ for any $H\in \tfrak$.

The subalgebra of $\slfrakhat_k$ generated by $\qfrak_m^i$, for $i\in\{1, \ldots, k-1\}$ and $m\in\Z\setminus\{0\}$, and $\cfrak$ is isomorphic to the
Heisenberg algebra $\hfrak_{\F,\Qfrak}$. This motivates the following
definition of Whittaker vector for $\slfrakhat_k$ (cf.\ \cite[Section 6]{art:christodoulopoulou2008}).
\begin{definition}\label{def:whittakersl}
Let $\chi\colon \Ucal(\hfrak_{\F, \Qfrak}^+)\to \F$ be an algebra
homomorphism such that $\chi \vert_{\hfrak_{\F, \Qfrak}^+}\neq 0$, and let $V$ be a $\Ucal\big(\, \slfrakhat_k \big)$-module. A nonzero vector $w\in V$ is called a \emph{Whittaker vector of type} $\chi$ if $\eta\triangleright w=\chi(\eta)\, w$ for all $\eta\in \Ucal(\hfrak_{\F, \Qfrak}^+)$.
\end{definition}

\subsection{Frenkel-Kac construction}\label{sec:frenkelkac}

Let $\Vcal$ be a representation of $\hfrak_{\F,\Qfrak}$. We say that
it is a \emph{level one} representation if the central element $\cfrak$ acts by the
identity map. Henceforth we let $\Vcal$ be a level one representation of $\hfrak_{\F,\Qfrak}$ such that for any $v\in \Vcal$ there exists an integer $m(v)$ for which 
\begin{equation}\label{eq:condition-fock}
\big(\mathfrak{q}^{l_1}_{m_1}\cdots \mathfrak{q}^{l_a}_{m_a}\big) \triangleright v=0
\end{equation}
if $m_i>0$ and $\sum_i\, m_i>m(v)$.

Fix an index $j\in\{0,1, \ldots, k-1\}$ and consider the coset $\Qfrak+\omega_j$. Denote by $\F[\Qfrak+\omega_j]$ the group algebra of $\Qfrak+\omega_j$
over $\F$. For a root $\gamma\in \Qfrak$ we define the generating function
$V(\gamma,z)\in \mathrm{End}(\Vcal\otimes \F[\Qfrak+\omega_j])[[z,z^{-1}]]$ of
operators on $\Vcal\otimes \F[\Qfrak+\omega_j]$ by the
{bosonic vertex operator}
\begin{align}
V(\gamma,z)
&=\V_{1,-1}^\gamma(z) \, \exp(\log z\ \cfrak+\gamma)\\[4pt] & =\exp\Big(\sum_{m=1}^\infty\, \frac{z^m}{m}\,
\qfrak^\gamma_{-m}\, \Big)\, \exp\Big(-\sum_{m=1}^\infty\,
\frac{z^{-m}}{m}\, \qfrak^\gamma_{m}\Big)\, \exp(\log z\ \cfrak+\gamma)\ ,
\end{align}
where $\exp(\log z\ \cfrak+\gamma)$ is the operator defined by
\begin{equation}
\exp(\log z\ \cfrak+\gamma)\triangleright(v \otimes
[\beta+\omega_j]):=z^{\frac12\langle \gamma,\gamma\rangle_{\Qfrak} + \langle
  \gamma, \beta+\omega_j\rangle_{\Qfrak}}\, (v \otimes [\beta+\gamma+\omega_j])
\end{equation}
for $v \otimes [\beta+\omega_j] \in \Vcal\otimes
\F[\Qfrak+\omega_j]$. 
\begin{remark}
Here for the operator $\exp(\log z\ \cfrak+\gamma)$ we follow the notation in \cite[Section 3.2.1]{art:nagao2009}. In the existing literature, this operator is denoted in various different ways.
\end{remark}

Let $V_m(\gamma)\in {\rm End}(\Vcal\otimes
\F[\Qfrak+\omega_j])$
denote the operator defined by the formal Laurent series expansion $V(\gamma,z)=\sum_{m\in\Z} \,
V_m(\gamma)\, z^m$. 

We define a map $\epsilon\colon \Qfrak\times \Qfrak\to \{\pm 1\}$ by
\begin{equation}
\epsilon(\gamma_i,\gamma_j)=\left\{
\begin{array}{ll}
-1 \ , & j=i, i+1\ ,\\
1 \ , & \mbox{otherwise} \ ,
\end{array}\right.
\end{equation}
with the properties
$\epsilon(\gamma+\gamma',\beta)=\epsilon(\gamma,\beta)\, \epsilon(\gamma',\beta)$
and $\epsilon(\gamma,\beta+\beta'\,
)=\epsilon(\gamma,\beta)\, \epsilon(\gamma,\beta'\, )$.

\begin{theorem}[{\cite[Theorem 1]{art:frenkelkac1980}}]\label{thm:frenkelkac}
Let $j\in\{0,1, \ldots, k-1\}$ and let $\Vcal$ be a level one representation of $\hfrak_{\F,\Qfrak}$ satisfying the condition \eqref{eq:condition-fock}. Then the vector space $\Vcal \otimes \F[\Qfrak+\omega_j]$ carries a level one $\slfrakhat_k$-module structure given by
\begin{align}
(H_i\otimes 1)\triangleright(v \otimes [\beta+\omega_j])&=\big(\langle
\gamma_i,\beta\rangle_\Qfrak+\delta_{ij}\big)\, (v \otimes [\beta+\omega_j])\ ,\\[4pt]
(H_i\otimes t^m)\triangleright(v \otimes [\beta+\omega_j])&= \big(\mathfrak{q}^i_m\triangleright v\big) \otimes [\beta+\omega_j]\ ,\\[4pt]
(E_i\otimes t^m)\triangleright(v \otimes
[\beta+\omega_j])&=\epsilon(\gamma_i,\beta)\, V_{m+\delta_{ij}}(\gamma_i)\triangleright(v \otimes [\beta+\omega_j])\ ,\\[4pt]
(F_i\otimes t^m)\triangleright(v \otimes [\beta+\omega_j])&=
\epsilon(\beta,\gamma_i)\, V_{-m-\delta_{ij}}(-\gamma_i)\triangleright(v \otimes [\beta+\omega_j])\ ,
\end{align}
for $i\in\{1, \ldots, k-1\}$ and $m\in\Z\setminus\{0\}$. If $\Vcal$ is
the Fock space of $\hfrak_{\F,\Qfrak}$, then $\Vcal \otimes \F[\Qfrak+\omega_j]$ is the $j$-th dominant representation of $\slfrakhat_k$.
\end{theorem}

\subsubsection{Virasoro operators}\label{sec:virasoro-affine}

Let $\{\eta_i\}_{i=1}^{ k-1}$ be an orthonormal basis of the vector
space $\Qfrak\otimes_\Z \R$. The Virasoro algebra associated with $\hfrak_{\F, \Qfrak}\subset\slfrakhat_k$ has generators $\cfrak$ and $L^{\slfrakhat_k}_n$ for $n\in\Z$ defined by \cite[Section 2.8]{art:frenkelkac1980}
\begin{align}
L^{\slfrakhat_k}_0&=\sum_{i=1}^{k-1}\ \sum_{m=1}^\infty\, \qfrak_{-m}^{\eta_i}\
\qfrak_{m}^{\eta_i}+\frac{1}{2}\ \sum_{i=1}^{k-1}\,
\big(\qfrak_0^{\eta_i}\big)^2 \ ,\\[4pt] 
L^{\slfrakhat_k}_n &= \frac{1}{2}\ \sum_{i=1}^{k-1}\ \sum_{m\in\Z}\,
\qfrak_{-m}^{\eta_i}\ \qfrak_{m+n}^{\eta_i}\qquad \mbox{for}\quad n\in\Z\setminus\{0\}\ .
\end{align}
Note that distinct orthonormal bases of $\Qfrak\otimes_\Z \R$ give
rise to the
same Virasoro algebra $\virfrak_\F$.

\bigskip \section{AGT relations on $\R^4$\label{sec:AGTonR4}}

\subsection{Equivariant cohomology of $\hilb{n}{\C^2}$}\label{sec:C2}

In the following we shall give a brief survey of results concerning the
equivariant cohomology of the Hilbert schemes $\hilb{n}{\C^2}$ and representations of Heisenberg algebras thereon \cite{art:nakajima1996, art:grojnowski1996, book:nakajima1999, art:vasserot2001, art:liqinwang2004, art:qinwang2007, art:nakajima2014}.

Let us consider the action of the torus $T:=(\C^*)^2$ on the
complex affine plane $\C^2$ given by $(t_1,t_2)\triangleright
(x,y)=(t_1\,x,t_2\,y)$, and the induced $T$-action on the
{Hilbert scheme of $n$ points} $\hilb{n}{\C^2}$ which is the
{fine} moduli space parameterizing zero-dimensional subschemes of
$\C^2$ of length $n$; it is a smooth quasi-projective variety of
dimension $2n$. Following \cite{art:ellingsrudstromme1987}, the $T$-fixed points of $\hilb{n}{\C^2}$ are zero-dimensional subschemes of $\C^2$ of length $n$ supported at the origin $0\in\C^2$ which correspond to partitions $\lambda$ of $n$. We shall denote by $Z_\lambda$ the fixed point in $\hilb{n}{\C^2}^T$ corresponding to the partition $\lambda$ of $n$. 

For $i=1,2$ denote by $t_i$ the $T$-modules corresponding to the
characters $\chi_i\colon (t_1,t_2)\in T \mapsto t_i \in\C^*$, and by
$\varepsilon_i$ the equivariant first Chern class of $t_i$. Then $H_T^*({\rm
  pt}; \C)=\C[\varepsilon_1,\varepsilon_2]$ is the coefficient ring
for the $T$-equivariant cohomology. The equivariant Chern
character of the tangent space to $\hilb{n}{\C^2}$ at a fixed point
$Z_\lambda$ is given by
\begin{equation}
\mathrm{ch}_T\big(T_{Z_\lambda}\hilb{n}{\C^2} \big)=\sum_{s\in
  Y_\lambda}\, \big(\e^{(L(s)+1)\, \varepsilon_1-A(s)\,
  \varepsilon_2}+\e^{-L(s)\, \varepsilon_1+(A(s)+1)\, \varepsilon_2}\big)\ .
\end{equation}
The equivariant Euler class is therefore given by
\begin{equation}
\eu_T\big(T_{Z_\lambda}\hilb{n}{\C^2} \big)=(-1)^n\,
\eu_+(\lambda)\, \eu_-(\lambda)\ ,
\end{equation}
where 
\begin{equation}
\eu_+(\lambda)=\prod_{s\in Y_\lambda}\, \big((L(s)+1)\,
\varepsilon_1-A(s)\, \varepsilon_2\big)\qquad \mbox{and} \qquad
\eu_-(\lambda) =\prod_{s\in Y_\lambda}\, \big(L(s)\,
\varepsilon_1-(A(s)+1)\, \varepsilon_2\big)\ .
\end{equation}
\begin{remark}
By \cite[Corollary~3.20]{art:nakajima2014}, $\eu_+(\lambda)$ is the equivariant Euler class of the nonpositive part $T^{\leq 0}_{Z_\lambda}$ of the tangent space to $\hilb{n}{\C^2}$ at the fixed point $Z_\lambda$.
\end{remark}

Let $\imath_{\lambda}\colon \{Z_\lambda\}\hookrightarrow
\hilb{n}{\C^2}$ be the inclusion morphism and define the class
\begin{equation}\label{eq:fixed-class}
[\lambda]:={\imath_\lambda}_*(1) \ \in \ H^{4n}_T\big(\hilb{n}{\C^2}\big)\ .
\end{equation}
By the projection formula we get
\begin{equation}
[\lambda]\cup[\mu]=\delta_{\lambda,\mu}\,
\eu_T\big(T_{Z_\lambda}\hilb{n}{\C^2}
\big)\, [\lambda]=(-1)^n\, \delta_{\lambda,\mu}\,
\eu_+(\lambda)\, \eu_-(\lambda)\, [\lambda]\ .
\end{equation}
Denote 
\begin{equation}
\imath_n:=\bigoplus_{Z_\lambda \in\hilb{n}{\C^2}^T}\, \imath_{\lambda}
\, \colon \, \hilb{n}{\C^2}^T\ \longrightarrow \ \hilb{n}{\C^2}\ .
\end{equation}
Let $\imath_n^!\colon H_T^\ast\left(\hilb{n}{\C^2}^T\right)_{\mathrm{loc}}\to
H_T^\ast\big(\hilb{n}{\C^2} \big)_{\mathrm{loc}}$ be the induced Gysin map, where
\begin{equation}
H_T^\ast(-)_{\mathrm{loc}}:=H_T^\ast(-)\otimes_{\C[\varepsilon_1,\varepsilon_2]}\C(\varepsilon_1,\varepsilon_2)
\end{equation}
is the localized equivariant cohomology. By the localization theorem, $\imath_n^!$ is an isomorphism and its inverse is given by
\begin{equation}
\big(\imath_n^! \big)^{-1}\, \colon \, A \ \longmapsto \ \Big(\,
\frac{\imath_{\lambda}^\ast(A)}{\eu_T\big(T_{Z_\lambda}\hilb{n}{\C^2}
  \big)}\, \Big)_{Z_\lambda\in\hilb{n}{\C^2}^T}\ .
\end{equation}
Henceforth we denote $\mathbb{H}_{\C^2,n}:=H_T^{\ast}(\hilb{n}{\C^2})_{\mathrm{loc}}$. Define the bilinear form 
\begin{align}\label{eq:bilinearform-C2}
\langle -,- \rangle_{\mathbb{H}_{\C^2,n}} \, \colon \,
\mathbb{H}_{\C^2,n}\times \mathbb{H}_{\C^2,n} & \ \longrightarrow \
\C(\varepsilon_1, \varepsilon_2)\ ,\\ \nonumber 
(A, B)& \ \longmapsto \ (-1)^n \, p_n^! \, \big(\imath_n^!
\big)^{-1}(A \cup B) \ ,
\end{align}
where $p_n$ is the projection of $\hilb{n}{\C^2}^T$ to a point.
\begin{remark}
Our sign convention in defining the bilinear form is different from
the one used e.g. in \cite{book:nakajima1999,art:carlssonokounkov2012}. We choose this convention because,
under the isomorphism to be introduced later on in
\eqref{eq:isomorphism}, the form \eqref{eq:bilinearform-C2} becomes
exactly the Jack inner product \eqref{eq:jackinnerprod}. This convention produces various sign changes compared to
previous literature. Hence every time we state that a given result coincides with what is known in the literature, the reader should keep in mind ``up to the sign convention we choose".
\end{remark}
Following \cite[Section 2.2]{art:liqinwang2004} we define the distinguished classes
\begin{equation}
[\alpha_\lambda]:=\frac{1}{\eu_+(\lambda)}\, [\lambda] \
\in \ H^{2n}_T\big(\hilb{n}{\C^2} \big)_{\mathrm{loc}}\ .
\end{equation}
For $\lambda,\mu$ partitions of $n$ one has
\begin{align}\label{eq:product-C2}
\big\langle[\alpha_\lambda]\,,\,[\alpha_\mu]
\big\rangle_{\mathbb{H}_{\C^2,n}}&=\delta_{\lambda,\mu}\, \frac{\eu_-(\lambda)}{\eu_+(\lambda)}
\\[4pt] &=\delta_{\lambda,\mu}\, \prod_{s\in Y_\lambda} \,
\frac{L(s)\, \varepsilon_1-\big(A(s)+1 \big)\,
  \varepsilon_2}{\big(L(s)+1 \big)\, \varepsilon_1-A(s)\varepsilon_2}
=\delta_{\lambda,\mu}\, \prod_{s\in Y_\lambda} \, \frac{L(s)\,
  \beta+A(s)+1}{\big(L(s)+1 \big)\, \beta+A(s)}\ ,\nonumber
\end{align}
where
\begin{equation}\label{eq:beta}
\beta=-\frac{\varepsilon_1}{\varepsilon_2} \ .
\end{equation}
\begin{remark}
In \cite[Section 3(v)]{art:nakajima2014}, Nakajima gives a geometric interpretation of the class $[\alpha_\lambda]$.
\end{remark}

By the localization theorem and Equation \eqref{eq:product-C2}, the
classes $[\alpha_\lambda]$ form a
$\C(\varepsilon_1,\varepsilon_2)$-basis for the infinite-dimensional vector space
$\mathbb{H}_{\C^2}:=\bigoplus_{n\geq 0} \, \mathbb{H}_{\C^2,n}$. Hence the symmetric bilinear form \eqref{eq:bilinearform-C2} is nondegenerate. The forms $\langle -,- \rangle_{\mathbb{H}_{\C^2,n}}$ define a symmetric bilinear form 
\begin{equation}
\langle -,- \rangle_{\mathbb{H}_{\C^2}}\, \colon \,
\mathbb{H}_{\C^2}\times \mathbb{H}_{\C^2} \ \longrightarrow \ \C(\varepsilon_1, \varepsilon_2)
\end{equation}
by imposing that $\mathbb{H}_{\C^2,n_1}$ and $\mathbb{H}_{\C^2,n_2}$
are orthogonal for $n_1\neq n_2$. Then $\langle -,-
\rangle_{\mathbb{H}_{\C^2}}$ is also nondegenerate.

The unique partition of $n=1$ is $\lambda=(1)$. Let us denote by $[\alpha]:=[\alpha_{(1)}]$ the corresponding class. Then
\begin{equation}
\big\langle[\alpha]\,,\,[\alpha] \big\rangle_{\HH_{\C^2}}=\beta^{-1}\ .
\end{equation}
Let us denote by $D_x$ and $D_y$ respectively the $x$ and $y$ axes of
$\C^2$. By localization, the corresponding equivariant cohomology
classes in $H_T^\ast(\C^2)_{\mathrm{loc}}$ are given by
\begin{align}
\left[D_x\right]_T=\frac{[0]}{\varepsilon_1}=
\frac{[0]}{\eu_{+}(1)}=[\alpha]\qquad \mbox{and} \qquad 
\left[D_y\right]_T=\frac{[0]}{\varepsilon_2}=-\beta\,
[\alpha]\ .
\end{align}

\subsection{Heisenberg algebra}

Following \cite{art:nakajima1996, book:nakajima1999}, for an integer $m>0$ define the
Hecke correspondences
\begin{equation}
D_x(n,m):= \big\{(Z,Z'\, )\in\hilb{n+m}{\C^2}\times\hilb{n}{\C^2}\
\big| \ Z'\subset Z \ , \
\mathrm{supp}(\mathcal{I}_{Z'}/\mathcal{I}_Z)=\{y\}\subset D_x \big\}\ ,
\end{equation}
where $\mathcal{I}_Z,\mathcal{I}_{Z'}$ are the ideal sheaves corresponding to $Z,Z'$ respectively. Let $q_1,q_2$ denote the projections of $\hilb{n+m}{\C^2}\times\hilb{n}{\C^2}$ to the two factors, respectively. Define linear operators $\pfrak_{-m}([D_x]_T)\in\mathrm{End}(\mathbb{H}_{\C^2})$ by
\begin{equation}
\pfrak_{-m}([D_x]_T)A :=q_1^!\big(q_2^\ast A\cup[D_x(n,m )]_T \big)
\end{equation}
for $A\in H_T^\ast(\hilb{n}{\C^2})_{\mathrm{loc}}$. We also define
$\pfrak_{m}([D_x]_T)\in\mathrm{End}(\mathbb{H}_{\C^2})$ to be the
adjoint operator of $\pfrak_{-m}([D_x]_T)$ with respect to the
inner product $\langle-,-\rangle_{\mathbb{H}_{\C^2}}$ on
$\mathbb{H}_{\C^2}$. As the
class $[D_x]_T$ spans $H_T^\ast(\C^2)_{\mathrm{loc}}$ over the field
$\C(\varepsilon_1,\varepsilon_2)$, we can define operators
$\pfrak_{m}(\eta)\in\mathrm{End}(\mathbb{H}_{\C^2})$ for every class
$\eta\in H_T^\ast(\C^2)_{\mathrm{loc}}$.
\begin{theorem}[{see \cite{art:nakajima1996, art:nakajima2014}}]\label{thm:nakajimaoperators}
The linear operators $\pfrak_{m}(\eta)$, for $m\in\Z\setminus\{0\}$ and $\eta\in H_T^\ast(\C^2)_{\mathrm{loc}}$, satisfy the Heisenberg commutation relations
\begin{equation}
\big[\pfrak_{m}(\eta_1)\,,\,\pfrak_{n}(\eta_2) \big]=m\,
\delta_{m,-n}\,
\langle\eta_1,\eta_2\rangle_{\mathbb{H}_{\C^2,1}}\,\mathrm{id}\qquad\mbox{and}\qquad
\big[\pfrak_{m}(\eta)\,,\,\mathrm{id} \big]=0\ .
\end{equation}
The vector space $\HH_{\C^2}$ becomes the Fock space of the Heisenberg algebra $\hfrak_{\HH_{\C^2,1}}$ modelled on $\HH_{\C^2,1}=H_T^{\ast}(\C^2)_{\mathrm{loc}}$ with the unit $\vert0\rangle$ in $H_T^0(\hilb{0}{\C^2})_{\mathrm{loc}}$ as highest weight vector.
\end{theorem}
\begin{remark}
Since $[D_x]_T=[\alpha]$, we have $\pfrak_{m}([\alpha])=\pfrak_{m}([D_x]_T)$.
\end{remark}
Henceforth we denote by $\hfrak_{\C^2}$ the Heisenberg algebra
$\hfrak_{\HH_{\C^2,1}}$, and we define
\begin{equation}\label{eq:heisenberg-C2}
\pfrak_{m}:=\pfrak_{m}([D_x]_T) \qquad \mbox{for} \quad
m\in\Z\setminus \{0\} \ ,
\end{equation}
so that one has the nonzero commutation relations
\begin{equation}\label{eq:commutheisenberg-C2}
[\pfrak_{-m},\pfrak_m]=m \, \beta^{-1}\,\mathrm{id}\ .
\end{equation}
Since $[D_x]_T$ generates $H_T^\ast(\C^2)_{\mathrm{loc}}$ over $\C(\varepsilon_1, \varepsilon_2)$, the operators $\pfrak_m$ generate $\hfrak_{\C^2}$.

Let $\lambda=(1^{m_1}\,2^{m_2}\,\cdots)$ be a partition. Define $\pfrak_\lambda:=\prod_i\, \pfrak_{-i}^{m_i}$. Then
\begin{equation}
\big\langle\mathfrak{p}_\lambda\vert0\rangle\,,\,\pfrak_\mu\vert0\rangle
\big\rangle_{\mathbb{H}_{\C^2}}=\delta_{\lambda,\mu}\, z_\lambda \, \beta^{-\ell(\lambda)}\ .
\end{equation}
Let us denote by $\Lambda_\beta$ the ring of symmetric functions in
infinitely many variables $\Lambda_{\C(\varepsilon_1,\varepsilon_2)}$
over the field $\C(\varepsilon_1,\varepsilon_2)$, equiped with the
Jack inner product \eqref{eq:jackinnerprod}. 
%The following result is established in \cite[Theorem
%3.2]{art:liqinwang2004}, in the antidiagonal action case (i.e., $t=t_1=t_2^{-1}$),
%and in \cite[Section 1.5]{art:carlssonokounkov2012} (see also \cite{art:nakajima1996}).
\begin{theorem}[{see \cite{art:nakajima1996, art:liqinwang2004, art:carlssonokounkov2012}}]\label{thm:actionC2}
There exists a $\C(\varepsilon_1,\varepsilon_2)$-linear isomorphism 
\begin{equation}\label{eq:isomorphism}
\phi\, \colon\, \mathbb{H}_{\C^2} \ \longrightarrow \ \Lambda_\beta
\end{equation}
preserving bilinear forms such that
\begin{equation}
\phi(\pfrak_\lambda\vert0\rangle)=p_\lambda(x) \qquad\mbox{and}\qquad \phi([\alpha_\lambda])=
J_\lambda(x;\beta^{-1})\ . 
\end{equation}
Via the isomorphism $\phi$, the operators $\pfrak_m$ act on $\Lambda_\beta$
as multiplication by $p_{-m}$ for $m<0$ and as $m\, \beta^{-1}\,
\frac{\partial}{\partial p_m}$ for $m >0$.
\end{theorem}

\subsubsection{Whittaker vectors}

We characterize a particular class of Whittaker vectors
(cf. Definition \ref{def:whittakerheisenberg}) which will be useful in
our studies of gauge theories.
\begin{proposition}\label{prop:whittakerC2}
Let $\eta \in \C (\varepsilon_1,\varepsilon_2)$. In the completed Fock space $\prod_{n\geq 0}\,\mathbb{H}_{\C^2,n}$, every vector of the form
\begin{equation}
G(\eta):=\exp\left( \eta \ \pfrak_{-1}\right)\vert0\rangle
\end{equation}
is a Whittaker vector of type $\chi_\eta$, where the algebra
homomorphism $\chi_\eta\colon
\Ucal(\hfrak^+_{\C^2} )\to\C(\varepsilon_1,\varepsilon_2)$ is defined by
\begin{equation}
\chi_\eta(\pfrak_1)=\eta \, \beta^{-1} \qquad
\mbox{and} \qquad \chi_\eta(\pfrak_n)=0 \qquad \mbox{for} \quad n>1 \ .
\end{equation}
\end{proposition}
\proof
The statement follows from the formal expansion
\begin{equation}\label{eq:formalexpansion}
G(\eta)= \sum_{n=0}^\infty\, \frac{\eta^n}{n!}\, (\pfrak_{-1})^n\vert 0\rangle
\end{equation}
with respect to the vector $\vacuum$, together with the relation $\pfrak_m\vert 0\rangle=0$ for $m>0$ and the identity
\begin{equation}
\pfrak_m\,(\pfrak_{-1})^n= n\, \beta^{-1}\, \delta_{m,1}\,(\pfrak_{-1})^{n-1}+ (\pfrak_{-1})^n\, \pfrak_m
\end{equation}
in $\Ucal(\hfrak_{\C^2} )$ for $m\geq 1$. 
\endproof
 
\subsection{Vertex operators}\label{sec:carlsson-okounkov}

Let $T_\mu=\C^\ast$ and $H^\ast_{T_\mu}(\mathrm{pt};\C)=\C[\mu]$. Let us denote by $\Ocal_{\C^2}(\mu)$ the trivial line bundle on $\C^2$
on which $T_\mu$ acts by scaling the fibers. In~\cite{art:carlssonokounkov2012}, Carlsson
and Okounkov define a vertex operator $\V(\Lcal,z)$ for any
smooth quasi-projective surface $X$ and any line bundle $\Lcal$ on
$X$. Here we shall describe only $\V(\Ocal_{\C^2}(\mu),z)$;
see~\cite{art:carlssonokounkov2012} for a complete description of such
types of vertex operators.

Let $\Zbf_n\subset \hilb{n}{\C^2}\times \C^2$ be the universal
subscheme, whose fiber over a point $Z\in\hilb{n}{\C^2}$ is the
subscheme $Z\subset\C^2$ itself. Consider
\begin{equation}
{\boldsymbol\Zcal}_i:=p_{i3}^\ast(\Ocal_{\Zbf_{n_i}}) \ \in \ K\big(\hilb{n_1}{\C^2}\times
\hilb{n_2}{\C^2}\times \C^2\big) \qquad \mbox{for} \quad i=1, 2\ ,
\end{equation}
where $p_{ij}$ denotes the projection to the $i$-th and $j$-th factors. Define the virtual vector bundle
\begin{equation}
\Ebf_\mu^{n_1, n_2}={p_{12}}_\ast\left(({\boldsymbol\Zcal}_1^\vee+{\boldsymbol\Zcal}_2-{\boldsymbol\Zcal}_1^\vee\cdot
  {\boldsymbol\Zcal}_2)\cdot p_3^\ast(\Ocal_{\C^2}(\mu))\right) \ \in \
K\big(\hilb{n_1}{\C^2}\times \hilb{n_2}{\C^2} \big)\ ,
\end{equation}
where $p_3$ is the projection to $\C^2$. The fibre of $\Ebf_\mu^{n_1, n_2}$ over $(Z_1, Z_2)\in
\hilb{n_1}{\C^2}^T\times \hilb{n_2}{\C^2}^T$ is given by
\begin{equation}
\Ebf_\mu^{n_1, n_2}\big\vert_{(Z_1, Z_2)}=\chi\big(\Ocal_{\C^2}\,,\,
\Ocal_{\C^2}(\mu)
\big)-\chi\big(\Ical_{Z_1}\,,\,\Ical_{Z_2}\otimes\Ocal_{\C^2}(\mu) \big)\ ,
\end{equation}
where $\chi(E,F):=\sum_{i=0}^2 \, (-1)^i\, \mathrm{Ext}^i(E,F)$ for any pair of
coherent sheaves $E,F$ on $\C^2$, while its rank is
\begin{equation}
\rk\big(\Ebf_\mu^{n_1, n_2} \big)= n_1+n_2 \ .
\end{equation}

Define the operator
$\V(\Ocal_{\C^2}(\mu),z)\in\mathrm{End}(\mathbb{H}_{\C^2})[[z,
z^{-1}]]$ by its matrix elements
\begin{multline}\label{eq:carlssonokounkov}
(-1)^{n_2} \, \big\langle \V(\Ocal_{\C^2}(\mu),z)A_1 \,,\, A_2
\big\rangle_{\mathbb{H}_{\C^2}} \\ :=z^{n_2-n_1} \,
\int_{\hilb{n_1}{\C^2}\times \hilb{n_2}{\C^2}} \,
\eu_T\big(\Ebf_\mu^{n_1, n_2}\big)\cup p_1^\ast(A_1) \cup p_2^\ast(A_2) \ ,
\end{multline} 
where $A_i\in H_T^\ast(\hilb{n_i}{\C^2})_{\mathrm{loc}}$ and $p_i$ is the projection from
$\hilb{n_1}{\C^2}\times \hilb{n_2}{\C^2}$ to the $i$-th factor for
$i=1,2$. By \cite[Lemma
6]{art:carlssonokounkov2012}, the matrix elements \eqref{eq:carlssonokounkov} in the fixed point basis are given by
\begin{align}\label{eq:carlokfixedbasis} 
\big\langle \V(\Ocal_{\C^2}(\mu),z)[\lambda_1] \,,\, [{\lambda_2} ]
\big\rangle_{\mathbb{H}_{\C^2}}& =(-1)^{|\lambda_2|}\,
z^{|\lambda_2|-|\lambda_1|} \ 
\eu_T\big(\Ebf_\mu^{n_1, n_2}\big\vert_{(Z_{\lambda_1},Z_{\lambda_2})}
\big) \\[4pt] \nonumber
& = (-1)^{|\lambda_2|}\, z^{|\lambda_2|-|\lambda_1|} \ m_{Y_{\lambda_1}, Y_{\lambda_2}}(\varepsilon_1, \varepsilon_2,\mu)\ ,
\end{align}
where
\begin{align}\label{eq:m}
m_{Y_1, Y_2}(\varepsilon_1, \varepsilon_2,a):= &\prod_{s_1\in Y_1}\,
\big(a-L_{Y_2}(s_1)\,
\varepsilon_1+(A_{Y_1}(s_1)+1)\, \varepsilon_2\big) \\ \nonumber & \times \ \prod_{s_2 \in
  Y_2}\, \big(a+(L_{Y_1}(s_2)+1)\,
\varepsilon_1-A_{Y_2}(s_2)\, \varepsilon_2 \big)
\end{align}
for a pair of Young tableaux $Y_1,Y_2$ and $a\in\C(\varepsilon_1,\varepsilon_2)$. In gauge theory this factorized expression for the matrix elements represents the contribution of the
\emph{bifundamental hypermultiplet}.

We shall now describe the operator $\V(\Ocal_{\C^2}(\mu),z)$ in terms
of the operators $\pfrak_{m}$ defined in Equation \eqref{eq:heisenberg-C2} for $m\in\Z\setminus\{0\}$. In our setting, \cite[Theorem 1]{art:carlssonokounkov2012} assumes the following form.
\begin{theorem}\label{thm:CO}
The operator $\V(\Ocal_{\C^2}(\mu),z)$ is a vertex operator in
Heisenberg operators given by the generalized bosonic exponential
associated with the Heisenberg algebra $\hfrak_{\C^2}$ as
\begin{equation}\label{eq:carlssonokounkov-heisenberg}
\V(\Ocal_{\C^2}(\mu),z)=\V_{-\frac{\mu}{\varepsilon_2}, \frac{\mu+ \varepsilon_1+\varepsilon_2}{\varepsilon_2}}(z)\ .    
\end{equation}
\end{theorem}
%\proof
%By \cite[Theorem 1]{art:carlssonokounkov2012}, the operator $\V(\Ocal_{\C^2}(\mu),z)$ can be expressed as
%\begin{multline}\label{eq:carlsson-okounkov}
%\V(\Ocal_{\C^2}(\mu),z) = \exp\Big(- \sum_{m=1}^\infty \, \frac{z^m}{m}\, \pfrak_{-m}\big(\eu_{T\times
%    T_\mu}(\Ocal_{\C^2}(\mu)) \big)\Big) \\ 
%\times \ \exp\Big(-\sum_{m=1}^\infty\, \frac{z^{-m}}{m}\, \pfrak_{m}\big(\eu_{T\times
%    T_\mu}(\mathcal{K}_{\C^2}\otimes\Ocal_{\C^2}(\mu)^\vee) \big)\Big)
%\end{multline}
%where $\mathcal{K}_{\C^2}$ is the canonical line bundle of $\C^2$. Since
%\begin{equation}\label{eq:eulerclasses}
%\imath_0^\ast\, \eu_{T\times
%  T_\mu}\big(\Ocal_{\C^2}(\mu)\big)=\mu \qquad\mbox{and}\qquad
%\imath_0^\ast \, \eu_{T\times
%  T_\mu}\big(\mathcal{K}_{\C^2}\otimes\Ocal_{\C^2}(\mu)^\vee \big)=-\mu-\varepsilon_1-\varepsilon_2\ ,
%\end{equation}
%and $1=\frac{[\alpha]}{\varepsilon_2}=\frac{[D_x]_T}{\varepsilon_2}$
%in $H_T^\ast(\C^2)_{\mathrm{loc}}$, the assertion follows.
%\endproof

\subsection{Integrals of motion}\label{sec:integralC2}

Let $\Vbf^n$ be the pushforward of $\Ebf_0^{n,n}$ with respect to the
projection of the product $\hilb{n}{\C^2}\times \hilb{n}{\C^2}$ to the
second factor. It is a $T$-equivariant vector bundle on
$\hilb{n}{\C^2}$ of rank $n$, which we shall call the natural bundle over $\hilb{n}{\C^2}$. The $T$-equivariant Chern character of $\Vbf^n$ at the fixed point
$Z_\lambda$ is given by
\begin{equation}
\ch_T\big(\Vbf^n\big\vert_{Z_\lambda}\big) = \sum_{s\in
  Y_\lambda}\, \e^{-L'(s)\, \varepsilon_1- A'(s)\, \varepsilon_2} \ .
\end{equation}

\begin{remark}
The vector bundle $\Vbf^n$ can equivalently be defined as the pushforward with
respect to the projection $\hilb{n}{\C^2}\times \C^2\to\hilb{n}{\C^2}$
of the structure sheaf of the universal subscheme $\Zbf_n$. In the
literature $\Vbf^n$ is also called the tautological sheaf and denoted $\Ocal^{[n]}$ (cf.\ \cite{art:lehn1999, book:lehn2004,art:okounkovpandharipande2010, art:nakajima2014}). 
\end{remark}

Let us denote by $\Vbf$ the natural bundle over $\coprod_{n\geq 0}\, \hilb{n}{\C^2}$. The operators
of multiplication by $\boldsymbol{I}_1:=\rk\big(\Vbf \big)$ and $\boldsymbol{I}_p:=
(\cc_{p-1} )_T\big(\Vbf \big)$ for $p\geq 2$ on $\prod_{n\geq
  0}\,\mathbb{H}_{\C^2,n}$ have even
degrees, are self-adjoint with respect to the inner product on $\prod_{n\geq
  0}\,\mathbb{H}_{\C^2,n}$, and commute with each other; they can thus be simultaneously
diagonalized in the
fixed point basis $[\lambda]$ of $\prod_{n\geq
  0}\,\mathbb{H}_{\C^2,n}$ (see \cite{art:okounkovpandharipande2010},
and also \cite[Section 4]{art:nakajima2014} where our operator
$\pfrak_m$ is denoted $P_m(\varepsilon_2)$ for  $m\in\Z\setminus \{0\}$). For example, one has
\begin{align}\label{eq:integralofmotions-fixedpoints}
\boldsymbol{I}_1\triangleright [\lambda]= |\lambda|\, [\lambda] \qquad \mbox{and} \qquad
\boldsymbol{I}_2\triangleright [\lambda]= -\sum_{s\in Y_\lambda}\,
\big(L'(s)\, \varepsilon_1+ A'(s)\, \varepsilon_2\big) \
[\lambda] \ .
\end{align}
As a consequence, the operators of multiplication by $\boldsymbol{I}_p$ for $ p\geq 1$ can be written in terms of the Heisenberg
operators \eqref{eq:heisenberg-C2} as elements of a commutative subalgebra of $\Ucal(\hfrak_{\C^2})$. For example, we have
\begin{align}\label{eq:rank-C2}
\boldsymbol{I}_1 =& \ \beta\,
\sum_{m=1}^\infty \, \pfrak_{-m}\, \pfrak_{m}\ ,\\[4pt] 
\boldsymbol{I}_2 =
& \ \varepsilon_1\, \Big(\, \frac{\beta}{2}\, \sum_{m,n=1}^\infty \, \big(
\pfrak_{-m}\, \pfrak_{-n}\, \pfrak_{m+n}+\pfrak_{-m-n}\, \pfrak_{n}\,
\pfrak_{m} \big) - \frac{\beta-1}{2}\, \sum_{m=1}^\infty\, (m-1)\,
\pfrak_{-m}\, \pfrak_{m}\,\Big)\ . \label{eq:c1-C2}
\end{align}
Note that the \emph{energy operator} $\boldsymbol{I}_1$ coincides with the Virasoro generator
$L_0^{\hfrak}$ from Section \ref{sec:virasoro-heisenberg}, while the
operator $\boldsymbol{I}_2$ is equal to $\varepsilon_1\, \Box^{\beta^{-1}}$, where $\Box^{\beta^{-1}}$ is the bosonized Hamiltonian of the
quantum trigonometric Calogero-Sutherland model with infinitely many particles and coupling constant $\beta^{-1}$.

\subsection{$\Ncal=2$ gauge theory}\label{sec:R4pure}

The Nekrasov partition function for pure $\Ncal=2$ $U(1)$ gauge
theory on $\mathbb{R}^4$ is given by the generating function \cite{art:nekrasov2003, art:bruzzofucitomoralestanzini2003}
\begin{align}
\Zcal_{\C^2}(\varepsilon_1,\varepsilon_2; \qsf) :=& \
\sum_{n=0}^\infty\, \qsf^n \ \int_{\hilb{n}{\C^2}} \, \big[\hilb{n}{\C^2}\big]_T\\[4pt]
 = & \ \sum_{n=0}^\infty \, (-\qsf)^n\, \big\langle
 [\hilb{n}{\C^2}]_T\,,\,[\hilb{n}{\C^2}]_T \big\rangle_{\mathbb{H}_{\C^2}}
\end{align}
where $\qsf\in\C^*$ with $|\qsf|<1$.
By the localization theorem we obtain
\begin{equation}
\Zcal_{\C^2}(\varepsilon_1,\varepsilon_2; \qsf) = \sum_{\lambda} \,
\Big(-\frac{\qsf}{\varepsilon_2^2}\, \Big)^{\vert\lambda\vert} \
\prod_{s\in Y_\lambda}\, \frac{1}{\big((L(s)+1 )\,\beta+A(s)\big)\,
  \big(L(s)\, \beta+(A(s)+1)\big)} 
\end{equation}
as in \cite[Equation (3.16)]{art:bruzzofucitomoralestanzini2003},
where the sum runs over all partitions $\lambda$.
\begin{remark}\label{rmk:exponentialnekrasov}
By \cite[Equation (4.7)]{art:nakajimayoshioka2005-I}, the partition
function can be summed explicitly and written in the closed form
\begin{equation}\label{eq:nakajimayoshioka}
\Zcal_{\C^2}(\varepsilon_1,\varepsilon_2; \qsf) = \exp\Big(\, \frac{\qsf}{\varepsilon_1\, \varepsilon_2}\, \Big) \ .
\end{equation}
\end{remark}

\subsubsection{Gaiotto state}\label{sec:gaiottoC2}

In \cite{art:gaiotto2009}, Gaiotto considers the inducing state of the
(completed) Verma module of the Virasoro algebra. It has the property
that it is a Whittaker vector for the Verma module, and the norm of
its $\qsf$-deformation coincides with the Nekrasov partition
function of pure $\Ncal=2$ $SU(2)$ gauge theory on $\R^4$. Here we
consider the analogous vector for $U(1)$ gauge theory on $\R^4$.

Following \cite{art:schiffmannvasserot2013}, we define the \emph{Gaiotto state} to be the sum of all fundamental classes
\begin{equation}
G := \sum_{n\geq 0} \, \left[ \hilb{n}{\C^2} \right]_T
\end{equation}
in the completed Fock space $\prod_{n \geq 0} \, \HH_{\C^2,n}$. We also introduce the \emph{weighted Gaiotto state} as the formal power series
\begin{equation}
G_\qsf := \sum_{n\geq 0} \, \qsf^n \, \left[ \hilb{n}{\C^2} \right]_T
\ \in \ \mbox{$\prod\limits_{n \geq 0}$}\, \qsf^n \, \HH_{\C^2,n}\ .
\end{equation}
Consider the bilinear form
\begin{equation}
\langle-,-\rangle_{\HH_{\C^2},\qsf} \,:\, 
\mbox{$\prod\limits_{n \geq 0}$}\, \qsf^n \, \HH_{\C^2,n} \times
\mbox{$\prod\limits_{n \geq 0}$}\, \qsf^n \, \HH_{\C^2,n} \
\longrightarrow \ \C(\varepsilon_1,\varepsilon_2)[[\qsf]]
\end{equation}
defined by
\begin{equation}
\Big\langle \, \mbox{$\sum\limits_{n\geq 0}$} \, \qsf^n\, \eta_n \,,\,
\mbox{$\sum\limits_{n\geq 0}$} \, \qsf^n \, \nu_n \,
\Big\rangle_{\HH_{\C^2},\qsf} :=  \sum_{n= 0}^\infty\, \qsf^n \
\int_{\hilb{n}{\C^2}}\, \eta_n \cup \nu_n =  \sum_{n= 0}^\infty\,
(-\qsf)^n \, \langle \eta_n, \nu_n \rangle_{\HH_{\C^2}}\ .
\end{equation}
It follows immediately that the norm of the weighted Gaiotto state is the Nekrasov partition function for $\Ncal=2$ $U(1)$ gauge theory on $\R^4$:
\begin{equation}\label{eq:gaiottonekrasovC2}
\Zcal_{\C^2}(\varepsilon_1,\varepsilon_2;\qsf) = \langle G_\qsf,
G_\qsf \rangle_{\HH_{\C^2},\qsf} \ .
\end{equation}
By Proposition \ref{prop:whittakerC2} we have the following result.
\begin{proposition}\label{prop:gaiottowhittakerC2}
The Gaiotto state $G$ is a Whittaker vector of type $\chi$, where the
algebra homomorphism $\chi\colon
\Ucal(\hfrak^+_{\C^2}) \to\C(\varepsilon_1,\varepsilon_2)$ is defined by
\begin{equation}
\chi(\pfrak_1)=-\frac{1}{\varepsilon_1} \qquad
\mbox{and} \qquad \chi(\pfrak_n)=0 \qquad \mbox{for} \quad n>1 \ .
\end{equation}
\end{proposition}
\proof
Let $\eta \in \C (\varepsilon_1,\varepsilon_2)$. 
By using the formal expansion \eqref{eq:formalexpansion} and the
isomorphism $\phi$, we can write $\phi(G(\eta))$ in terms of powers
$p_1^n$. By Lemma \ref{lem:p1} and simple algebraic manipulations we
can then rewrite the vector $G(\eta)$ as
\begin{equation}\label{eq:whittaker}
G(\eta)=\sum_{n\geq 0} \, (\eta\,\varepsilon_2)^n\,
\big[\hilb{n}{\C^2} \big]_T
\end{equation}
and the result follows.
\endproof

\subsection{Quiver gauge theories\label{sec:R4quivers}}

We now add matter fields to the $\Ncal=2$ gauge theory on $\R^4$. We consider the most general $\Ncal=2$ superconformal quiver gauge
theory with gauge group $U(1)^{r+1}$ for $r\geq0$, following the
general ADE classification of \cite[Chapter
3]{art:nekrasovpestun2012}. 

Let $\quiv=(\quiv_0,\quiv_1)$ be a quiver,
i.e., an oriented graph with a finite set of vertices $\quiv_0$,
a finite set of edges
$\quiv_1\subset\quiv_0\times\quiv_0$, and two projection maps
$\source,\tail:\quiv_1\rightrightarrows \quiv_0$ which assign to each oriented edge
its source and target vertex respectively. Representations of the
quiver encode the matter field content of the gauge theory. Fix a vector
$(n_\upsilon)_{\upsilon\in \quiv_0}\in\N^{\quiv_0}$ of integers labelled by the nodes of
the quiver $\quiv$, and consider the product of Hilbert schemes
$\prod_{\upsilon\in\quiv_0}\, \hilb{n_\upsilon}{\C^2}$. The vertices $\upsilon\in \quiv_0$
label $U(1)$ gauge groups and $m_\upsilon\geq0$ (resp. $\bar m_\upsilon\geq 0$)
fundamental (resp. antifundamental) hypermultiplets of masses
$\mu^s_{\upsilon}$, $s=1,\dots,m_\upsilon$ (resp. $\bar\mu^{\bar s}_{\upsilon}$,
$\bar s=1,\dots,\bar m_\upsilon$) which correspond to the $T$-equivariant vector
bundles $\Vbf_{\mu^s_{\upsilon}}^{n_\upsilon}$
(resp. $\bar\Vbf_{\bar\mu^{\bar s}_{\upsilon}}^{n_\upsilon}$) of rank $n_\upsilon$ on $\hilb{n_\upsilon}{\C^2}$ obtained by
pushforward of $\Ebf_{\mu^s_\upsilon}^{n_\upsilon,n_\upsilon}$
(resp. $\Ebf_{\bar\mu^{\bar s}_\upsilon}^{n_\upsilon,n_\upsilon}$) with respect to the
projection of $\hilb{n_\upsilon}{\C^2} \times\hilb{n_\upsilon}{\C^2}$ to the
second (resp. first) factor. The edges
$e\in \quiv_1$ label $U(1)\times U(1)$ bifundamental hypermultiplets of masses $\mu_e$ which correspond to
the vector bundles $\Ebf_{\mu_e}^{n_{\source(e)},n_{\tail(e)}}$ of
rank $n_{\source(e)}+n_{\tail(e)}$ on $\hilb{n_{\source(e)}}{\C^2}
\times\hilb{n_{\tail(e)}}{\C^2}$; if the edge $e$ is a vertex loop, i.e., $\source(e)=\tail(e)$, then the restriction of
$\Ebf_{\mu_e}^{n_{\source(e)},n_{\source(e)}}$ to the diagonal of $\hilb{n_{\source(e)}}{\C^2}
\times\hilb{n_{\source(e)}}{\C^2}$ describes an adjoint hypermultiplet
of mass $\mu_e$. The total matter field content of
the $\Ncal=2$ quiver gauge theory associated to $\quiv$
in the sector labelled by $(n_\upsilon)_{\upsilon\in \quiv_0}\in\N^{\quiv_0}$
is thus described by the bundle on $\prod_{\upsilon\in\quiv_0}\,
\hilb{n_\upsilon}{\C^2}$ given by
\begin{equation}
\boldsymbol{M}^{(n_\upsilon)}_{(\mu^s_{\upsilon}),(\bar\mu^{\bar s}_\upsilon),
  (\mu_e)}:= \bigoplus_{\upsilon\in\quiv_0}\, p_\upsilon^* \Big(\, \bigoplus_{s=1}^{m_\upsilon}\,
\Vbf_{\mu^s_{\upsilon}}^{n_\upsilon} \ \oplus \ \bigoplus_{\bar
  s=1}^{\bar m_\upsilon}\, \bar\Vbf_{\bar\mu^{\bar s}_{\upsilon}}^{n_\upsilon}\, \Big) \ \oplus \ \bigoplus_{e\in\quiv_1}\,
p_e^*\Ebf_{\mu_e}^{n_{\source(e)},n_{\tail(e)}} \ ,
\end{equation}
where $p_\upsilon$ is the projection of $\prod_{\upsilon\in\quiv_0}\,
\hilb{n_\upsilon}{\C^2}$ to the $\upsilon$-th factor and $p_e$ the projection
  to $\hilb{n_{\source(e)}}{\C^2}
\times\hilb{n_{\tail(e)}}{\C^2}$. 

For each vertex $\upsilon\in\quiv_0$, the degree of the Euler class of the
pushforward of this
bundle to the $\upsilon$-th factor is the integer
\begin{align}
d_\upsilon := & \ \dim_\C\hilb{n_\upsilon}{\C^2}- \rk\Big( 
\boldsymbol{M}^{(n_\upsilon)}_{(\mu^s_{\upsilon}),(\bar\mu^{\bar s}_\upsilon),
  (\mu_e)}\Big\vert_{\hilb{n_\upsilon}{\C^2}}\,\Big) \\[4pt] = & \ n_\upsilon\,
\big(2-m_\upsilon-\bar m_\upsilon-\#\{e\in\quiv_1\, | \, \source(e)=\upsilon\}-
\#\{e\in\quiv_1\, | \, \tail(e)=\upsilon\}\big) \ .
\end{align}
The $\Ncal=2$ quiver gauge theory is said to be \emph{conformal} if
$d_\upsilon=0$ for all $\upsilon\in\quiv_0$; it is
\emph{asymptotically free} if $d_\upsilon>0$. Note that with this
definition the pure $\Ncal=2$ gauge theory of Section \ref{sec:R4pure}
is asymptotically free.
As explained in \cite[Chapter 3]{art:nekrasovpestun2012}, $\Ncal=2$ asymptotically free quiver gauge theories can be recovered from conformal theories, so in the following we restrict our
attention to superconformal quiver gauge theories. 

Introduce coupling constants $\qsf_\upsilon\in\C^*$ with
$|\qsf_\upsilon|<1$ at each vertex
$\upsilon\in\quiv_0$, and let $T_{\mubf}$ be the maximal torus of the
total flavour symmetry group
\begin{equation}
G_{\rm f}= 
\prod_{\upsilon\in\quiv_0}\, GL(m_\upsilon,\C)\times GL(\bar m_\upsilon,\C)  \
\times \ \prod_{e\in\quiv_1}\, \C^*
\end{equation}
with $H_{T_{\mubf}}^*({\rm
  pt};\C)=\C[(\mu_e),(\mu_\upsilon^s),(\bar\mu_\upsilon^{\bar s})]$. Then the quiver
gauge theory partition function is defined by the generating function
\begin{multline}
\Zcal_{\C^2}^{\quiv}(\varepsilon_1,\varepsilon_2,\mubf;\qbf):=
\sum_{(n_\upsilon)\in\N^{\quiv_0}}\, \qbf^{\boldsymbol n} \
  \int_{\prod\limits_{\upsilon\in\quiv_0}\, \hilb{n_\upsilon}{\C^2}}\, \eu_{T\times
  T_{\mubf}}\big(\boldsymbol{M}^{(n_\upsilon)}_{(\mu^s_{\upsilon}),(\bar\mu^{\bar
    s}_\upsilon),
  (\mu_e)} \big) \\[4pt] = \sum_{(n_\upsilon)\in\N^{\quiv_0}} \, \qbf^{\boldsymbol n} \
  \int_{\prod\limits_{\upsilon\in\quiv_0}\, \hilb{n_\upsilon}{\C^2}}\
  \prod_{\upsilon\in\quiv_0}\, p_\upsilon^*\Big(\, \prod_{s=1}^{m_\upsilon}\,
  \eu_T\big(\Vbf_{\mu^s_{\upsilon}}^{n_\upsilon}\big) \
  \prod_{\bar s=1}^{\bar m_\upsilon}\,
  \eu_T\big(\bar\Vbf_{\bar\mu^{\bar s}_{\upsilon}}^{n_\upsilon} \big)\, \Big) \\
  \times \ \prod_{e\in\quiv_1}\,
  p_e^*\, \eu_T\big(\Ebf_{\mu_e}^{n_{\source(e)},n_{\tail(e)}}
  \big) \ ,
\end{multline}
where $\qbf^{\nbf}:= \prod_{\upsilon\in\quiv_0}\, \qsf_\upsilon^{n_\upsilon}$. By
the localization theorem, we obtain
\begin{multline}
\Zcal_{\C^2}^{\quiv}(\varepsilon_1,\varepsilon_2,\mubf;\qbf) =
\sum_{(\lambda^\upsilon)}\, (-\qbf)^{\lambdabf} \ \prod_{\upsilon\in\quiv_0}\,
\frac{\prod\limits_{s=1}^{m_\upsilon}\,
  m_{Y_{\lambda^\upsilon}}\big(\varepsilon_1,\varepsilon_2,\mu_\upsilon^s \big) \
  \prod\limits_{\bar s=1}^{\bar m_\upsilon}\,
  m_{Y_{\lambda^\upsilon}}\big(\varepsilon_1,\varepsilon_2,\bar
  \mu_\upsilon^{\bar s} + \varepsilon_1 +
  \varepsilon_2\big)}{m_{Y_{\lambda^\upsilon},
    Y_{\lambda^\upsilon}}(\varepsilon_1,\varepsilon_2
  ,0)} \label{eq:ZcalC2quiver} \\ \times \
\prod_{e\in\quiv_1}\,
m_{Y_{\lambda^{\source(e)}},Y_{\lambda^{\tail(e)}}}(
\varepsilon_1,\varepsilon_2 ,\mu_e)
\end{multline}
where $\qbf^{\lambdabf}:= \prod_{\upsilon\in\quiv_0}\,
\qsf_\upsilon^{|\lambda^\upsilon|}$ for a collection of partitions
$\lambda^\upsilon$ associated to the vertices of the quiver, and
\begin{equation}
m_Y(\varepsilon_1,\varepsilon_2,a):= \prod_{s\in Y}\, \big(a-L'(s)\,
\varepsilon_1-A'(s)\, \varepsilon_2\big)
\end{equation}
for a Young tableau $Y$ and $a\in\C(\varepsilon_1,\varepsilon_2)$.

The conformal constraint
\begin{equation}
m_\upsilon+\bar m_\upsilon+ \#\big\{e\in\quiv_1\, \big| \, \source(e)=\upsilon \big\} +
\#\big\{e\in\quiv_1\, \big| \, \tail(e)=\upsilon \big\} = 2
\label{eq:confconstr}\end{equation}
for each $\upsilon\in\quiv_0$ severely restricts the possible quivers in the abelian gauge
theory. It is easy to check that the only admissible quivers in the ADE
classification of \cite[Chapter
3]{art:nekrasovpestun2012} are the linear (or chain)
quivers of the finite-dimensional $A_r$-type Dynkin diagram and the cyclic (or necklace)
quivers of the affine $\widehat{A}_{r}$-type extended Dynkin
diagram for some $r\geq0$\footnote{Here the $A_0$-type Dynkin diagram is the trivial quiver consisting of a single vertex with no arrows, and the $\widehat{A}_0$-type Dynkin diagram is the quiver consisting of a single vertex with a vertex edge loop.}. We consider in detail each case in turn.

\subsection{$\widehat{A}_{r}$ theories\label{sec:hatArtheories}}

For the cyclic quivers of type $\widehat{A}_{r}$
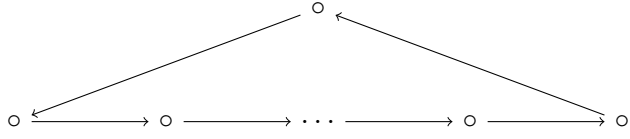
\begin{equation}
  \begin{tikzpicture}[xscale=2,yscale=-1]
\node (A0_3) at (3,0) {$\circ$};
\node (A1_1) at (1,1.5) {$\circ$}; 
\node (A1_2) at (2,1.5) {$\circ$}; 
\node (A1_3) at (3,1.5) {$\ldots$}; 
\node (A1_4) at (4,1.5) {$\circ$}; 
\node (A1_5) at (5,1.5) {$\circ$}; 
    \path (A1_1) edge [->]node [auto] {$\scriptstyle{}$} (A1_2);
    \path (A1_2) edge [->]node [auto] {$\scriptstyle{}$} (A1_3);
    \path (A1_3) edge [->]node [auto] {$\scriptstyle{}$} (A1_4);
    \path (A1_4) edge [->]node [auto] {$\scriptstyle{}$} (A1_5);
    \path (A0_3) edge [->]node [auto] {$\scriptstyle{}$} (A1_1);
    \path (A1_5) edge [->]node [auto] {$\scriptstyle{}$} (A0_3);
  \end{tikzpicture}
\end{equation}
with $r+1$ vertices and arrows, one has $m_\upsilon=\bar m_\upsilon=0$ by
Equation \eqref{eq:confconstr}. We
label the vertices $\quiv_0$ by $\upsilon=0,1,\dots, r$ with
counterclockwise orientation and read modulo $r+1$, and similarly for
the edges $e=(\upsilon,\upsilon+1) \in\quiv_1$. The partition function for the $\Ncal=2$
quiver gauge theory of type $\widehat{A}_{r}$ reads as
\begin{equation}\label{eq:ZC2hatAr}
\Zcal^{\widehat{A}_{r}}_{\C^2}(\varepsilon_1,\varepsilon_2,\mubf;\qbf)
= \sum_{\lambdabf}\, (-\qbf)^{\lambdabf} \ \prod_{\upsilon=0}^r\
\frac{m_{Y_{\lambda^\upsilon},Y_{\lambda^{\upsilon+1}}}(\varepsilon_1,\varepsilon_2,\mu_\upsilon)}{m_{Y_{\lambda^\upsilon},Y_{\lambda^{\upsilon}}}(\varepsilon_1,\varepsilon_2,0)}
\ ,
\end{equation}
where the sum is over all $r+1$-vectors of partitions
$\lambdabf=(\lambda^0,\lambda^1,\dots,\lambda^{r})$ with
$\lambda^{r+1}:= \lambda^0$ and
$\qbf^{\lambdabf}:=\prod_{\upsilon=0}^r \,
\qsf_\upsilon^{|\lambda^\upsilon|}$.

\subsubsection{Conformal blocks}

We will relate the partition function \eqref{eq:ZC2hatAr} to the trace of 
vertex operators
$\V(\Ocal_{\C^2}(\mu_\upsilon),z_\upsilon)$. We shall also denote by $\hfrak$ the Heisenberg algebra $\hfrak_{\C^2}$ to simplify the presentation.
\begin{proposition}\label{prop:hatArtrace}
The partition function of the $\widehat{A}_r$-theory on $\R^4$ is
given by
\begin{equation}
\Zcal^{\widehat{A}_{r}}_{\C^2}(\varepsilon_1,\varepsilon_2,\mubf;\qbf)
= \mathrm{Tr}_{\mathbb{H}_{\C^2}}\, \qsf^{L_0^\hfrak} \ \prod_{\upsilon=0}^r\,
\, \V(\Ocal_{\C^2}(\mu_\upsilon),z_\upsilon)
\end{equation}
independently of $z_0\in\C^*$, where $\qsf:=\qsf_0\,
\qsf_1\cdots\qsf_r$ and $z_\upsilon:= z_0\, \qsf_1\cdots \qsf_\upsilon$
for $\upsilon=1,\dots,r$.
\end{proposition}
\proof
By
Equation \eqref{eq:integralofmotions-fixedpoints}, the Virasoro
operator $L_0^{\hfrak}$ acts in $\mathbb{H}_{\C^2}$ as
\begin{equation}
L_0^{\hfrak} \, \big\vert_{\mathbb{H}_{\C^2,n}}=n\
\mathrm{id}_{\mathbb{H}_{\C^2,n}} \ ,
\end{equation}
and so the trace of products of Ext vertex operators $\V(\Ocal_{\C^2}(\mu_\upsilon),z_\upsilon)$ is given by the sum of their matrix elements over the fixed point basis as
\begin{equation}
\mathrm{Tr}_{\mathbb{H}_{\C^2}}\, \qsf^{L_0^{\hfrak}} \ \prod_{\upsilon=0}^r
\, \V(\Ocal_{\C^2}(\mu_\upsilon),z_\upsilon) =\sum_{\nbf\in\N^{r+1}}
\, \qsf^{n_0} \
\sum_{\lambdabf\, :\, \vert\lambda^\upsilon\vert=n_\upsilon}\
\prod_{\upsilon=0}^r\, \frac{\big\langle
  \V(\Ocal_{\C^2}(\mu_\upsilon),z_\upsilon)[\lambda^\upsilon]\,,\, [\lambda^{\upsilon+1}] \big\rangle_{\mathbb{H}_{\C^2}}}{\big\langle [\lambda^\upsilon]\,,\, [\lambda^\upsilon] \big\rangle_{\mathbb{H}_{\C^2}}}\ .
\end{equation}
By Equation \eqref{eq:carlokfixedbasis} we obtain
\begin{align}
\mathrm{Tr}_{\mathbb{H}_{\C^2}} \, \qsf^{L_0^{\hfrak}} \ \prod_{\upsilon=0}^r\,
\V(\Ocal_{\C^2}(\mu_\upsilon),1)& =\sum_{\nbf\in\N^{r+1}}
\, \qsf^{n_0} \ \prod_{\upsilon=0}^r\, (-1)^{n_\upsilon}\,
z_\upsilon^{n_\upsilon-n_{\upsilon+1}} \\ & \qquad \qquad \times \
\sum_{\lambdabf\, :\, \vert\lambda^\upsilon\vert=n_\upsilon} \
\prod_{\upsilon=0}^r\, 
\frac{\eu_T\big(
  \boldsymbol{E}_{\mu_\upsilon}^{n_\upsilon,n_{\upsilon+1}} \big\vert_{
    (Z_{\lambda^\upsilon} ,Z_{\lambda^{\upsilon+1}})}
  \big)}{\eu_T\big(T_{Z_{\lambda^\upsilon}}\hilb{n_\upsilon}{\C^2} \big)}\\[4pt]
&=\sum_{\nbf\in\N^{r+1}}
\, (-\qbf)^{\nbf} \
\sum_{\lambdabf\, :\, \vert\lambda^\upsilon\vert=n_\upsilon}\
\prod_{\upsilon=0}^r\
\frac{m_{Y_{\lambda^\upsilon},Y_{\lambda^{\upsilon+1}}}(\varepsilon_1,\varepsilon_2,\mu_\upsilon)}{m_{Y_{\lambda^\upsilon},Y_{\lambda^{\upsilon}}}(\varepsilon_1,\varepsilon_2,0)}
\ ,
\end{align}
and the result follows.
\endproof

\begin{remark}
Proposition \ref{prop:hatArtrace} shows that the partition function of
the $\widehat{A}_r$-theory coincides with the conformal block of the
Heisenberg algebra $\hfrak_{\C^2}$ on the
elliptic curve with nome $\qsf$ and $r+1$ punctures at
$z_0,z_1,\dots,z_r$; we can set $z_0=1$ without loss of
generality. The conformal dimension of the primary field inserted at
the $\upsilon$-th puncture is
\begin{equation}
\Delta(\mu_\upsilon;\varepsilon_1,\varepsilon_2)= \frac{\mu_\upsilon\,
  (\mu_\upsilon+\varepsilon_1+\varepsilon_2)}{2\varepsilon_1\,
  \varepsilon_2} \ .
\end{equation}
This elliptic curve is the Seiberg-Witten curve of the $\Ncal=2$
$U(1)^{r+1}$ quiver gauge theory on $\R^4$.
\end{remark}

By using the same arguments as in the proof of \cite[Corollary
1]{art:carlssonokounkov2012}, an explicit formula for the trace in
this case can be obtained using Equation
\eqref{eq:carlssonokounkov-heisenberg} and we arrive at the explicit evaluation of
the partition function.  In the following $\eta(\qsf):=\qsf^{\frac{1}{24}}\,
\prod_{n=1}^\infty\, (1-\qsf^n)$ denotes the Dedekind function.

\begin{proposition}\label{prop:tracededekind}
\begin{equation}
\Zcal^{\widehat{A}_{r}}_{\C^2}(\varepsilon_1,\varepsilon_2,\mubf;\qbf)
= \prod_{\upsilon=0}^r\,\big(\qsf_\upsilon^{-\frac{1}{24}}\, \eta(\qsf_\upsilon)\big)^{-\frac{\mu_\upsilon\,
    (\mu_\upsilon+\varepsilon_1+\varepsilon_2)}{\varepsilon_1\,
    \varepsilon_2}}\ \qsf^{\frac{1}{24}}\,\eta(\qsf)^{-1}\ .
\end{equation}
\end{proposition}
A similar formula for the $U(1)^{r+1}$ quiver gauge theory partition
function is conjectured in \cite[Appendix
C.2]{art:aldaygaiottotachikawa2010}.

\subsubsection{$\widehat{A}_{0}$ theory\label{sec:hatA0theory}}

The degenerate case $r=0$ of the $\widehat{A}_r$ quiver gauge theory
corresponds to the quiver consisting of a single
node with a vertex edge loop
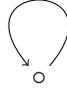
\begin{equation}
  \begin{tikzpicture}[auto]
\node (A0_0) at (0,0) {$\circ$};

    \path[->] (A0_0) edge[loop above, in=130,out=50,looseness=10, shorten >=-2pt, shorten <=-2pt] node {$\scriptstyle{}$} (A0_0);
  \end{tikzpicture}
\end{equation}
and is known as the $\Ncal=2^*$ gauge theory; it describes a single adjoint matter hypermultiplet of mass $\mu$ in
the $U(1)$ $\Ncal=2$ gauge theory on $\R^4$. Then the quiver gauge
theory partition function is the Nekrasov partition function
for $\Ncal=2^\ast$ gauge theory \cite{art:nekrasov2003,art:bruzzofucitomoralestanzini2003} and is given
by
\begin{equation}
\Zcal_{\C^2}^{\widehat{A}_{0}}(\varepsilon_1,\varepsilon_2, \mu; \qsf) =
\sum_{\lambda} \, \qsf^{\vert\lambda\vert} \ \prod_{s\in Y_\lambda}
\frac{\big( (L(s)+1)\, \varepsilon_1-A(s)\, \varepsilon_2+\mu \big)\,
  \big( L(s)\, \varepsilon_1-(A(s)+1)\, \varepsilon_2-\mu \big)}{\big(
  (L(s)+1)\, \varepsilon_1-A(s)\, \varepsilon_2 \big) \, \big( L(s)\,
  \varepsilon_1-(A(s)+1)\, \varepsilon_2 \big)}
\end{equation}
as in \cite[Equation (3.26)]{art:bruzzofucitomoralestanzini2003}.

By Proposition \ref{prop:tracededekind}, we have
\begin{equation}\label{eq:Ahat0dedekind}
\Zcal_{\C^2}^{\widehat{A}_{0}} (\varepsilon_1,\varepsilon_2, \mu;
\qsf)=\big(\qsf^{-\frac{1}{24}}\, \eta(\qsf)\big)^{-\frac{\mu\,
    (\mu+\varepsilon_1+\varepsilon_2)}{\varepsilon_1\, \varepsilon_2}-1}\ .
\end{equation}
A similar formula is written in \cite[Equation (2.28)]{art:wyllard2009}.
In the case of an antidiagonal torus action, i.e.,
$\varepsilon_1=-\varepsilon_2$, this result coincides with the formula
derived in
\cite[Equation (6.12)]{art:nekrasovokounkov2006}.
We can then rewrite Proposition \ref{prop:hatArtrace} in the following way.
\begin{corollary}
\begin{equation}
\Zcal^{\widehat{A}_{r}}_{\C^2}(\varepsilon_1,\varepsilon_2,\mubf;\qbf)
= \prod_{\upsilon=0}^r\, \Zcal_{\C^2}^{\widehat{A}_{0}} (\varepsilon_1,\varepsilon_2, \mu_\upsilon ;
\qsf_\upsilon)\ \frac{\prod\limits_{\upsilon=0}^r\,\eta(\qsf_\upsilon)}{\eta(\qsf)}\ .
\end{equation}
\end{corollary}

\subsection{${A}_{r}$ theories\label{sec:Artheories}}

Consider now the linear quivers of type ${A}_{r}$
\begin{equation}
  \begin{tikzpicture}[xscale=2.5,yscale=-.7]
\node (A1_0) at (0,0) {$\bullet$}; 
\node (A1_1) at (1,0) {$\circ$}; 
\node (A1_2) at (2,0) {$\ldots$}; 
\node (A1_3) at (3,0) {$\circ$}; 
\node (A1_4) at (4,0) {$\bullet$};  
    \path (A1_0) edge [->]node [auto] {$\scriptstyle{}$} (A1_1);
    \path (A1_1) edge [->]node [auto] {$\scriptstyle{}$} (A1_2);
    \path (A1_2) edge [->]node [auto] {$\scriptstyle{}$} (A1_3);
    \path (A1_3) edge [->]node [auto] {$\scriptstyle{}$} (A1_4);
  \end{tikzpicture}
\end{equation}
with $r+1$ vertices and $r$ arrows, where the solid nodes indicate the
insertion of a single fundamental or antifundamental
hypermultiplet. In this case we label vertices $\quiv_0$ from left to
right with $\upsilon=0,1,\dots,r$ and edges $\quiv_1$ with $e=(\upsilon,\upsilon+1)$;
for definiteness we take $\bar m_\upsilon=0$, so that $m_0=m_{r}=1$ by the
conformal constraints \eqref{eq:confconstr}. The partition function for the $\Ncal=2$ quiver
gauge theory of type $A_r$ for $r\geq1$ reads as
\begin{equation}\label{eq:ZC2Ar}
\Zcal^{{A}_{r}}_{\C^2}(\varepsilon_1,\varepsilon_2,\mubf;\qbf)
= \sum_{\lambdabf}\, (-\qbf)^{\lambdabf} \ 
\frac{m_{Y_{\lambda^0}}(\varepsilon_1,\varepsilon_2,\mu_0) \
  \prod\limits_{\upsilon=0}^{r-1}\,
  m_{Y_{\lambda^\upsilon},Y_{\lambda^{\upsilon+1}}}(\varepsilon_1,\varepsilon_2,\mu_{\upsilon+1})
\
m_{Y_{\lambda^r}}(\varepsilon_1,\varepsilon_2,\mu_{r+1})}{\prod\limits_{\upsilon=0}^r\,
m_{Y_{\lambda^\upsilon},Y_{\lambda^{\upsilon}}}(\varepsilon_1,\varepsilon_2,0)}
\ .
\end{equation}

\subsubsection{Conformal blocks}

We will express the partition function \eqref{eq:ZC2Ar} as a
particular matrix element of
Ext vertex operators.
\begin{proposition}\label{prop:fundamental-matter}
The partition function of the $A_r$-theory on $\R^4$ is given by
\begin{equation}
\Zcal^{{A}_{r}}_{\C^2}(\varepsilon_1,\varepsilon_2,\mubf;\qbf) =
\Big\langle \vacuum\,,\,
\prod_{\upsilon=0}^{r+1} \, \V\big(\Ocal_{\C^2}(\mu_{\upsilon}
),z_\upsilon \big)\vacuum\Big\rangle_{\mathbb{H}_{\C^2}}
\end{equation}
independently of $z_0\in\C^*$, where $z_\upsilon:=z_0\, \qsf_0\,
\qsf_1\cdots \qsf_{\upsilon-1}$ for $\upsilon=1,\dots,r+1$.
\end{proposition}
\proof
Arguing as in the proof of Proposition \ref{prop:hatArtrace}, we write
\begin{multline}
\Big\langle \vacuum\,,\,
\prod_{\upsilon=0}^{r+1} \, \V\big(\Ocal_{\C^2}(\mu_{\upsilon}
),z_\upsilon \big)\vacuum\Big\rangle_{\mathbb{H}_{\C^2}} \\ =\sum_{\nbf\in\N^{r+1}}
\ \sum_{\lambdabf\, :\, \vert\lambda^\upsilon\vert=n_\upsilon} \, \big\langle
  \V(\Ocal_{\C^2}(\mu_0),z_0)\,[\lambda^{0}] ,\, |0\rangle 
  \big\rangle_{\mathbb{H}_{\C^2}} \ 
\frac{\prod\limits_{\upsilon=0}^{r-1} \, \big\langle
  \V(\Ocal_{\C^2}(\mu_{\upsilon+1}),z_{\upsilon+1} )[\lambda^{\upsilon+1}]\,,\,
  [\lambda^{\upsilon}]
  \big\rangle_{\mathbb{H}_{\C^2}}}{\prod\limits_{\upsilon=0}^r\,
  \big\langle [\lambda^\upsilon]\,,\, [\lambda^\upsilon]
  \big\rangle_{\mathbb{H}_{\C^2}}} \\ \times \ \big\langle
  \V(\Ocal_{\C^2}(\mu_{r+1}),z_{r+1} )|0\rangle \,,\, [\lambda^{r}]
  \big\rangle_{\mathbb{H}_{\C^2}} \ ,
\end{multline}
and by Equation \eqref{eq:carlokfixedbasis} and the orthogonality relation \eqref{eq:product-C2} the result then follows.
\endproof

\begin{remark}
Proposition \ref{prop:fundamental-matter} expresses the partition
function of the $A_r$-theory as a conformal block of the Heisenberg
algebra $\hfrak_{\C^2}$ on the Riemann
sphere with $r+4$ punctures at $\infty,z_0,z_1,\dots,z_{r+1},0$; again
we can set $z_0=1$ without loss of generality. The conformal dimension
of the primary field at the insertion point $z_\upsilon$ is
$\Delta(\mu_\upsilon;\varepsilon_1,\varepsilon_2)$, while at
$\infty,0$ they are given respectively by
$\Delta(\tilde\mu_{\infty,0}; \varepsilon_1,\varepsilon_2)$, where the
masses $\tilde\mu_{\infty,0}$ obey
\begin{equation}
\tilde\mu_\infty+\tilde\mu_0=\varepsilon_1+\varepsilon_2
+\sum_{\upsilon=0}^{r+1}\, \mu_\upsilon \ .
\end{equation}
The Seiberg-Witten curve of the $\Ncal=2$ $U(1)^{r+1}$ quiver gauge
theory on $\R^4$ is a branched cover of this $(r+2)$-punctured Riemann
sphere, ramified over the points $\infty,0$.
\end{remark}

Using the vertex operator representation, we can again get a closed
formula for the combinatorial expansion \eqref{eq:ZC2Ar}.
\begin{proposition}\label{prop:ZC2Arexpl}
\begin{equation}
\Zcal^{{A}_{r}}_{\C^2}(\varepsilon_1,\varepsilon_2,\mubf;\qbf) =
\prod_{0\leq \upsilon<\upsilon'\leq r+1}\,
\big(1-\qsf_{\upsilon+1}\cdots
\qsf_{\upsilon'}\big)^{-\frac{\mu_{\upsilon'}\,
    (\mu_\upsilon+\varepsilon_1+\varepsilon_2)}{\varepsilon_1\,
    \varepsilon_2}} \ .
\end{equation}
\end{proposition}
\proof
Using Equation \eqref{eq:VOnprod} to express the product of vertex
operators in Proposition \ref{prop:fundamental-matter} in normal
ordered form, we can write
\begin{multline}
\prod_{\upsilon=0}^{r+1}\,
\V\big(\Ocal_{\C^2}(\mu_\upsilon),z_\upsilon\big) \vacuum = 
\prod_{0\leq \upsilon<\upsilon'\leq r+1}\, \Big(\,
1-\frac{z_{\upsilon'}}{z_\upsilon}\, \Big)^{-\frac{\mu_{\upsilon'}\,
    (\mu_\upsilon+\varepsilon_1+\varepsilon_2)}{\varepsilon_1\,
    \varepsilon_2}} \\ \times \ \exp\Big(-\sum_{\upsilon=0}^{r+1}\,
\frac{\mu_\upsilon}{\varepsilon_2} \ \sum_{m=1}^\infty\,
\frac{z_\upsilon^m}{m}\, \pfrak_{-m}\Big)\vacuum
\end{multline}
since $\pfrak_m\vacuum=0$ for all $m>0$. Since $\pfrak_m$ is the
adjoint operator of $\pfrak_{-m}$ with respect to the scalar product
on $\mathbb{H}_{\C^2}$, we have
$\big\langle\vacuum,(\pfrak_{-m})^n\vacuum\big\rangle_{\mathbb{H}_{\C^2}}=0$
for all $m,n\geq1$ and the result follows.
\endproof
A similar formula for the $U(1)^{r+1}$ quiver gauge theory partition
function is conjectured in \cite[Appendix
C.1]{art:aldaygaiottotachikawa2010}.

\subsubsection{$A_0$ theory}

The degenerate limit $r=0$ of the $A_r$ quiver gauge theory is built
on the trivial quiver consisting of a single vertex with no arrows
\begin{equation}
\bullet
\end{equation}
and $m_0=2$ fundamental matter fields by Equation
\eqref{eq:confconstr}. Then the quiver gauge theory partition function
is the Nekrasov partition function for $\Ncal=2$ gauge theory with two
fundamental matter hypermultiplets of masses $\mu_0,\mu_1$ 
\cite{art:nekrasov2003, art:bruzzofucitomoralestanzini2003} which is
given by
\begin{equation}
\Zcal_{\C^2}^{A_0}(\varepsilon_1,\varepsilon_2, \mu_0, \mu_1;
\qsf) = \sum_{\lambda} \, (-\qsf)^{\vert\lambda\vert} \ \prod_{s\in Y_\lambda}\, \frac{\big(L'(s)\, \beta -A'(s)+\tilde \mu_0\big)\, \big(L'(s)\, \beta- A'(s)+\tilde \mu_1\big)}{\big((L(s)+1)\, \beta+A(s) \big)\, \big(L(s)\, \beta + A(s)+1\big)}
\end{equation}
as in \cite[Equation (3.22)]{art:bruzzofucitomoralestanzini2003},
where $\tilde \mu_0 =\mu_0/\varepsilon_2$ and $\tilde
\mu_1=\mu_1/\varepsilon_2$. By Proposition
\ref{prop:fundamental-matter} this partition function computes the
four-point conformal block for the Heisenberg algebra
$\hfrak_{\C^2}$ on the Riemann sphere with
primary field insertions at $\infty,1,\qsf,0$, and by Proposition
\ref{prop:ZC2Arexpl} the combinatorial sum can be evaluated explicitly
with the result
\begin{equation}
\Zcal_{\C^2}^{A_0}(\varepsilon_1,\varepsilon_2, \mu_0, \mu_1;
\qsf) = (1-\qsf)^{- \frac{\mu_1\, (\mu_0+\varepsilon_1+\varepsilon_2)}{\varepsilon_1\, \varepsilon_2}}\ .
\label{eq:fundamental}\end{equation}
A similar expression is written in \cite[Equation (2.27)]{art:wyllard2009}.
In the antidiagonal limit $\beta=1$, this formula coincides with the
partition function expression derived in \cite[Equation (49)]{art:marshakovmironovmorozov2010}.

\bigskip \section{Moduli spaces of framed sheaves\label{sec:sheaves}}

\subsection{Orbifold compactification of $X_k$}

In this subsection we recall the construction of the orbifold
compactification of the minimal resolution of $\C^2/\Z_k$ from
\cite[Section 3]{art:bruzzopedrinisalaszabo2013} and describe the main
results that we will use in this paper. For background to the theory of root and toric stacks used in the construction, see \cite[Section 2]{art:bruzzopedrinisalaszabo2013}, and to the theory of framed sheaves on (projective) Deligne-Mumford stacks, see \cite{art:bruzzosala2015}.

Fix an integer $k\geq 2$ and denote by $\mu_k$ the group of $k$-th
roots of unity in $\C$. A choice of a primitive $k$-th root of unity
$\omega$ defines an isomorphism of groups $\mu_k\simeq \Z_k$. We
define an action of $\mu_k\simeq \Z_k$ on $\C^2$ as $\omega\triangleright (x,y):= (\omega\, x, \omega^{-1}\, y)$.
The quotient $\C^2/\Z_k$ is a normal affine toric surface. The origin is the only singular point of $\C^2/\Z_k$, and is a particular case of a rational double point or du~Val singularity \cite[Definition~10.4.10]{book:coxlittleschenck2011}.

Let $\varphi_k\colon X_k\to\C^2/\Z_k$ be the minimal resolution of the
singularity of $\C^2/\Z_k$; it is a smooth toric surface with $k$
torus-fixed points $p_1, \ldots, p_k$ and $k+1$ torus-invariant
divisors $D_0, D_1, \ldots, D_k$ which are smooth projective curves of genus zero. For any $i=1, \ldots, k$ the divisors $D_{i-1}$ and $D_i$ intersect at the point $p_i$. Moreover, $D_1, \ldots, D_{k-1}$ are the irreducible components of the exceptional divisor $\varphi_k^{-1}(0)$. 
By the McKay correspondence, there is a one-to-one
correspondence between the irreducible representations of $\mu_k$ and
the divisors $D_1, \ldots, D_{k-1}$ \cite[Corollary~10.3.11]{book:coxlittleschenck2011}. By \cite[Equation~(10.4.3)]{book:coxlittleschenck2011}, the intersection matrix $(D_i\cdot D_j)_{1\leq i,j\leq k-1}$ is given by minus the Cartan matrix $C$ of type $A_{k-1}$, i.e., one has
\begin{equation}
\left( D_i \cdot D_j \right)_{1\leq i,j\leq k-1} = -C= 
\begin{pmatrix}
-2 & 1 & \cdots & 0 \\
1 & -2 & \cdots & 0 \\
\vdots & \vdots & \ddots & \vdots\\
0 & 0 & \cdots & -2
\end{pmatrix} \ .
\end{equation}
The surface $X_k$ is an ALE space of type $A_{k-1}$.

Let $U_i$ be the torus-invariant affine open subset of $X_k$ which is
a neighbourhood of the torus-fixed point $p_i$ for $i=1, \ldots,
k$. Its coordinate ring is given by $\C[U_i]:=\C[T_1^{2-i}\, T_2^{1-i},T_1^{i-1}\, T_2^i]$ for $i=1, \ldots, k$. By imposing the change of variables $T_1=t_1^k$ and $T_2=t_2 \, t_1^{1-k}$, we have
\begin{equation}\label{eq:coordinatering-i}
\C[U_i]=\C[t_1^{k-i+1}\,
t_2^{1-i},t_1^{i-k}\, t_2^{i}]\ .
\end{equation}
Define
\begin{equation}
\chi^i_1(t_1,t_2)=t_1^{k-i+1}\, t_2^{1-i}\qquad
\mbox{and}\qquad\chi_2^i(t_1,t_2)=t_1^{i-k}\, t_2^{i}\ .
\end{equation}
After identifying characters of $T$ with one-dimensional $T$-modules, let $\varepsilon_j^{(i)}$ denote the equivariant first Chern class of
$\chi^i_j$ for $i=1,\dots,k$ and $j=1,2$. Then
\begin{equation}
\varepsilon_1^{(i)}(\varepsilon_1,\varepsilon_2)=(k-i+1)\,
\varepsilon_1-(i-1)\,
\varepsilon_2\qquad\mbox{and}\qquad\varepsilon_2^{(i)}(\varepsilon_1,\varepsilon_2)=
-(k-i)\, \varepsilon_1+i\, \varepsilon_2\ .
\end{equation}

One can compactify the ALE space $X_k$ to a normal projective toric surface $\bar
X_k$ by adding a torus-invariant divisor $D_\infty\simeq \PP^1$ such
that for $k=2$ the surface $\bar X_2$ coincides with the second
Hirzebruch surface $\F_2$. For $k\geq3$ the surface $\bar X_k$ is
singular, but one can associate with $\bar X_k$ its canonical toric
stack $\Xscr_k^{\mathrm{can}}$ which is a two-dimensional projective
toric orbifold with Deligne-Mumford torus $T$ and coarse moduli space
$\pi_k^{\mathrm{can}}\colon \Xscr_k^{\mathrm{can}}\to \bar{X}_k$. By
\emph{canonical} we mean that the locus over which
$\pi_k^{\mathrm{can}}$ is not an isomorphism has non-positive
dimension; for $k=2$ one has $\Xscr_2^{\mathrm{can}}\simeq
\F_2$. Consider the one-dimensional, torus-invariant, integral closed
substack
$\tilde{\Dscr}_\infty:=(\pi_k^{\mathrm{can}})^{-1}(D_\infty)_{\mathrm{red}}\subset\Xscr_k^{\mathrm{can}}$.
By performing the $k$-th root construction on $\Xscr_k^{\mathrm{can}}$
along $\tilde{\Dscr}_\infty$ to extend the automorphism group of a
generic point of $\tilde{\Dscr}_\infty$ by $\mu_k$, we obtain a two-dimensional projective toric orbifold $\Xscr_k$ with Deligne-Mumford torus $T$ and coarse moduli space $\pi_k\colon\Xscr_k\to\bar X_k$. The surface $X_k$ is isomorphic to the open subset $\Xscr_k \setminus\Dscr_\infty$ of $\Xscr_k$, where $\Dscr_\infty :=\pi_k^{-1}(D_\infty)_{\mathrm{red}}$. Let $\Dscr_i:=\pi_k^{-1}(D_i)_{\mathrm{red}}$ be the divisors in $\Xscr_k$ corresponding to $D_i$ for $i=1, \ldots, k-1$. The classes
\begin{equation}
- \sum_{j=1}^{k-1}\, \big(C^{-1}\big)^{ij} \, \Dscr_j
\end{equation}
are integral for $i=1, \ldots, k-1$, where the inverse of the Cartan
matrix $C$ is given by
\begin{equation}
\big(C^{-1}\big)^{ij}=\frac{i\, (k-j)}k \qquad \mbox{for} \quad i\leq
j \ .
\end{equation}
Denote by $\Rcal_i$ the associated line bundles on $\Xscr_k$; the restrictions of $\Rcal_i$ to $X_k$ are precisely the tautological line bundles of Kronheimer and Nakajima \cite{art:kronheimernakajima1990}. 
\begin{proposition}[{\cite[Proposition 3.25]{art:bruzzopedrinisalaszabo2013}}]
The Picard group $\Pic(\Xscr_k)$ of $\Xscr_k$ is freely generated over $\Z$ by $\Ocal_{\Xscr_k}(\Dscr_\infty)$ and $\Rcal_i$ with $i=1, \ldots, k-1$.
\end{proposition}

The divisor $\Dscr_\infty$ can be characterized as a toric Deligne-Mumford stack with Deligne-Mumford torus $\C^\ast\times \Bscr\mu_k$ and coarse moduli space $r_k\colon \Dscr_\infty\to D_\infty$.
\begin{proposition}[{\cite[Proposition 3.27]{art:bruzzopedrinisalaszabo2013}}]
The divisor $\Dscr_\infty$ is isomorphic as a toric Deligne-Mumford stack to the global toric quotient stack
\begin{equation}
\left[\frac{\C^2\setminus\{0\}}{\C^\ast\times\mu_k}\right]\ ,
\end{equation}
where the group action is given in \cite[Equation (3.28)]{art:bruzzopedrinisalaszabo2013}.
\end{proposition}
\begin{corollary}[{\cite[Corollary 3.29]{art:bruzzopedrinisalaszabo2013}}]
The Picard group $\Pic(\Dscr_\infty)$ is isomorphic to $\Z\oplus
\Z_k$. It is generated by the line bundles $\Lcal_1$ and $\Lcal_2$
corresponding respectively to the characters
\begin{equation}
\chi_1\,\colon\, (t,\omega)\in\C^\ast\times\mu_k \ \longmapsto\
t\in\C^\ast\quad \mbox{and}\quad \chi_2\, \colon\,
(t,\omega)\in\C^\ast\times\mu_k \ \longmapsto\ \omega\in\C^\ast\ .
\end{equation}
\end{corollary}
For $j=0,1, \ldots, k-1$ define the line bundles
\begin{equation}
\Ocal_{\Dscr_\infty}(j)=\left\{
\begin{array}{ll}
\Lcal_2^{\otimes j} & \mbox{for even $k$}\ ,\\[8pt]
\Lcal_2^{\otimes j\, \frac{k+1}{2}} & \mbox{for odd $k$}\ .
\end{array}
\right.
\end{equation}
\begin{proposition}[{\cite[Corollary 3.34]{art:bruzzopedrinisalaszabo2013}}]
The restrictions of the tautological line bundles $\Rcal_j$ to
$\Dscr_\infty$ are given by
\begin{equation}
{\Rcal_j}\big\vert_{ \Dscr_\infty}\simeq \Ocal_{\Dscr_\infty}(j)\ .
\end{equation}
\end{proposition}
\begin{remark}
In \cite{art:bruzzopedrinisalaszabo2013} the line bundles $\Ocal_{\Dscr_\infty}(j)$ are the line bundles $\Ocal_{\Dscr_\infty}(s,j)$ for $s=0$. Indeed, one can prove that the degree of $\Ocal_{\Dscr_\infty}(j)$ is zero. Moreover, $\Ocal_{\Dscr_\infty}(0,j)$ can be endowed with a unitary flat connection associated with the $j$-th 
irreducible unitary representation $\rho_j$ of $\Z_k$ for $j=0,1, \ldots, k-1$ (cf.\ \cite[Remark 6.5]{art:eyssidieuxsala2013}).
\end{remark}

\subsection{Rank one framed sheaves}

\begin{definition}
Fix $j\in\{0,1,\dots,k-1\}$. A \emph{rank one $(\Dscr_\infty,\Ocal_{\Dscr_\infty}(j))$-framed sheaf} on $\Xscr_k$ is a pair $(\Ecal, \phi_{\Ecal})$, where $\Ecal$ is a
torsion-free sheaf on $\Xscr_k$ of rank one which is locally free in a
neighbourhood of $\Dscr_\infty$, and $\phi_{\Ecal}\colon
\Ecal\big\vert_{\Dscr_\infty}\xrightarrow{\sim} \Ocal_{\Dscr_\infty}(j)$
is an isomorphism. We call $\phi_{\Ecal}$ a \emph{framing} of
$\Ecal$. A \emph{morphism} between
$(\Dscr_\infty,\Ocal_{\Dscr_\infty}(j))$-framed sheaves
$(\Ecal,\phi_\Ecal)$ and $(\Gcal,\phi_\Gcal)$ of rank one is a morphism  
$f\colon \Ecal\to\Gcal$ such that $\phi_\Gcal\circ
f\big\vert_{\Dscr_\infty} = \phi_\Ecal$.
\end{definition}
\begin{remark}
By \cite[Remark 4.3]{art:bruzzopedrinisalaszabo2013}, the Picard group
of $\Xscr_k$ is isomorphic to the second singular cohomology group of
$\Xscr_k$ with integral coefficients via the first Chern class map $\crm_1$. Thus fixing the determinant line bundle of a coherent sheaf $\Ecal$ on $\Xscr_k$ is equivalent to fixing its first Chern class. 
\end{remark}
Given a vector $\vec{u}= (u_1,\ldots,u_{k-1}) \in\Z^{k-1}$, we denote by $\Rcal^{\vec{u}}$ the line bundle $\bigotimes_{i=1}^{k-1}\, \Rcal_i^{\otimes u_i}$ and by $\Rcal_0$ the trivial line bundle $\Ocal_{\Xscr_k}$.
\begin{lemma}[{\cite[Lemma 4.4]{art:bruzzopedrinisalaszabo2013}}]
Let $(\Ecal, \phi_{\Ecal})$ be a rank one $(\Dscr_\infty,\Ocal_{\Dscr_\infty}(j))$-framed sheaf on $\Xscr_k$. Then the determinant $\det(\Ecal)$ of $\Ecal$ is of the form $\Rcal^{\vec{u}}$, where the vector $\vec{u}\in \Z^{k-1}$ satisfies the condition
\begin{equation}\label{eq:condition-determinant}
\sum_{i=1}^{k-1}\, i\, u_i= j \ \bmod{k}  \ .
\end{equation}
\end{lemma}
\begin{remark}\label{rem:firstchern}
Set $\vec{v}:=C^{-1}\vec{u}$. Then Equation \eqref{eq:condition-determinant} implies the relations
\begin{equation}\label{eq:v-condition}
k\, v_l=-l\, j \ \bmod{k}
\end{equation}
for $l=1, \ldots, k-1$. Note that a component $v_l$ is integral if and
only if every component of $\vec v$ is integral. We subdivide the vectors $\vec{
u}\in\Z^{k-1}$ according to Equation \eqref{eq:condition-determinant} as
\begin{equation}
\Ufrak_j:=\Big\{ \vec{u}\in\Z^{k-1} \ \Big\vert \
\mbox{$\sum\limits_{i=1}^{k-1}$}\, i\, u_i = j \ \bmod{k} 
\Big\}\ .
\end{equation}
Define now a bijective map by identifying a vector
$\vec{u}\in\Z^{k-1}$ with $\sum_{i=1}^{k-1} \, u_i \crm_1(\Rcal_i) =
\sum_{i=1}^{k-1}\, u_i \, \omega_i$ as
\begin{equation}
\psi\, \colon\, \vec{u}\in \Z^{k-1} \ \longmapsto \ \sum_{i=1}^{k-1}
\, u_i \, \omega_i\in \Pfrak\ .
\end{equation}
It is natural to split this map according to the coset decomposition
\eqref{eq:lateralweightlattice} as
\begin{equation}
\psi^{-1}(\Qfrak + \omega_j) = \Ufrak_j\ ,
\end{equation}
which means that $\psi(\vec{u}\, )$ for $\vec{u}\in \Ufrak_j$ is naturally written as
 a sum of the fundamental weight $\omega_j$ and an element $\gamma_{\vec{u}}$ of the root
 lattice $\Qfrak$, which is given by
\begin{equation}
\gamma_{\vec{u}}:=\sum_{i=1}^{k-1}\, \Big(\, \sum_{l=1}^{k-1} \, \big(C^{-1}\big)^{
il}\, u_l-\big(C^{-1}\big)^{ij}\, \Big)\, \gamma_i = 
\sum_{i=1}^{k-1}\, \Big( v_i -\big(C^{-1}\big)^{ij}\Big)\, \gamma_i \
\in \ \Qfrak\ .
\end{equation}
We write
\begin{equation}
\psi_j:=\psi\big\vert_{\Ufrak_j} \, \colon \, \Ufrak_j \ \longrightarrow \
\Qfrak + \omega_j \ .
\label{eq:psij}\end{equation}
\end{remark}

Following \cite[Section 4]{art:bruzzopedrinisalaszabo2013}, let
$\Mcal(\vec{u},n,j)$ be the fine moduli space parameterizing
$(\Dscr_\infty,\Ocal_{\Dscr_\infty}(j))$-framed sheaves of rank one on
$\Xscr_k$ with determinant line bundle $\Rcal^{\vec{u}}$ and second
Chern class $n\in \Z$; the vector $\vec{u}$ belongs to $\Ufrak_j$. Let $p_{\Xscr_k}:\Xscr_k\times \Mcal(\vec{u},n,j) \to \Xscr_k$ be the projection.
As explained in \cite[Remark 4.7]{art:bruzzopedrinisalaszabo2013}, by ``fine'' one means that there exists a \emph{universal framed sheaf} $(\boldsymbol{\Ecal}, \boldsymbol{\phi}_{\boldsymbol{\Ecal}})$, where
$\boldsymbol{\Ecal}$ is a coherent sheaf on $\Mcal(\vec{u},n,j)\times \Xscr_k$ which is flat over $\Mcal(\vec{u},n,j)$, and $\boldsymbol{\phi}_{\boldsymbol{\Ecal}}\colon \boldsymbol{\Ecal}\to
p_{\Xscr_k}^\ast(\Ocal_{\Dscr_\infty}(j))$ is a morphism such that
its restriction to  $\Mcal(\vec{u},n,j)\times \Dscr_\infty$ is an isomorphism; the fibre over $[(\Ecal,\phi_\Ecal)]\in \Mcal(\vec{u},n,j)$ is itself the $(\Dscr_\infty,\Ocal_{\Dscr_\infty}(j))$-framed sheaf $(\Ecal,\phi_\Ecal)$ on $\Xscr_k$. In the following we shall call $\boldsymbol{\Ecal}$ the \emph{universal sheaf}.
\begin{theorem}[{\cite[Theorem 4.13]{art:bruzzopedrinisalaszabo2013}}]
The moduli space $\Mcal(\vec{u},n,j)$ is a smooth quasi-projective variety of dimension $2n$. The Zariski tangent
space of
$\Mcal(\vec{u},n,j)$ at a point $[(\Ecal,\phi_{\Ecal})]$ is $\mathrm{Ext}^1(\Ecal,\Ecal\otimes \Ocal_{\Xscr_k}(-\Dscr_\infty))$.
\end{theorem}

As explained in \cite[Section 4.3]{art:bruzzopedrinisalaszabo2013}, the Hilbert
scheme of $n$ points $\hilb{n}{X_k}$ of $X_k$ is isomorphic to
$\Mcal(\vec u,n,j)$ for any $\vec{u}\in \Ufrak_j$. For this, let $\imath:X_k\hookrightarrow \Xscr_k$ be the
inclusion morphism. If $Z$ is a point of $\hilb{n}{X_k}$ and
$\vec{y}\in\Z^{k-1}$, then the coherent sheaf
$\Ecal:=\imath_\ast(\Ical_Z)\otimes \Rcal^{\boldsymbol e_j-C\vec{y}}$ is a
rank one torsion-free sheaf on $\mathscr{X}_k$ with a framing
$\phi_{\Ecal}$ induced by the canonical isomorphism
$\Rcal^{\boldsymbol e_j-C\vec{y}}\, \vert_{\Dscr_\infty} \xrightarrow{\sim}
\Ocal_{\Dscr_\infty}(j)$ such that $(\Ecal, \phi_{\Ecal})$ is a rank
one $(\Dscr_\infty,\Ocal_{\Dscr_\infty}(j))$-framed sheaf with
determinant line bundle $\Rcal^{\vec{u}}$, where $\vec{u}:=\boldsymbol
e_j-C\vec{y}$, and second Chern class $n$. Thus $Z$ induces a point $[(\mathcal{E},\phi_{\mathcal{E}})]$ in $\Mcal(\vec{u},n,j)$. This defines an inclusion morphism
\begin{equation}
\tilde\imath_{\vec{u},n,j}\, \colon\, \hilb{n}{X_k} \ \hookrightarrow \ \Mcal(\vec{u},n,j)
\end{equation}
which is an isomorphism of fine moduli spaces by \cite[Proposition 4.16]{art:bruzzopedrinisalaszabo2013}.
\begin{remark}\label{rem:quivervariety}
In \cite{art:kuznetsov2007} it is shown that the Hilbert scheme of points $\hilb{n}{X_k}$ is isomorphic to a Nakajima quiver variety of type $\widehat{A}_{k-1}$ with suitable dimension vectors. Thus $\Mcal(\vec{u},n,j)$ is a quiver variety.
\end{remark}

\subsection{Equivariant cohomology}

We define a $T$-action on $\Mcal(\vec{u},n,j)$ in the following way (cf.\ \cite[Section 4.6]{art:bruzzopedrinisalaszabo2013}). For $(t_1,t_2)\in T$ let $F_{(t_1,t_2)}$ be the automorphism of $\Xscr_k$ induced by the torus action; then the $T$-action is given by
\begin{equation}
(t_1,t_2)\triangleright  \big[(\Ecal,\phi_{\Ecal}) \big]:=\big[\big( (F_{(t_1,t_2)}^{-1})^\ast(\Ecal)\,,\, \phi_\Ecal'\big) \big] \ ,
\end{equation}
where $\phi_\Ecal'$ is the composition of isomorphisms
\begin{equation}
\phi_\Ecal' \, \colon\,  \big(F_{(t_1,t_2)}^{-1} \big)^\ast\Ecal \big\vert_{\Dscr_\infty} \ \xrightarrow{(F_{(t_1,t_2)}^{-1})^\ast(\phi_{\Ecal})} \
\big(F_{(t_1,t_2)}^{-1} \big)^\ast \Ocal_{\Dscr_\infty}(j)\ \longrightarrow \ \Ocal_{\Dscr_\infty}(j) \ ;
\end{equation}
here the last arrow is given by the $T$-equivariant structure induced
on $\Ocal_{\Dscr_\infty}(j)$ by restriction of the torus action of
$\Xscr_k$ to $\mathscr D_\infty$. 
Note that the $T$-action on $X_k$ naturally lifts to $\hilb{n}{X_k}$ and the isomorphism $\tilde\imath_{\vec{u},n,j}$ is equivariant with respect to these torus actions.
\begin{proposition}[{\cite[Proposition 4.22]{art:bruzzopedrinisalaszabo2013}}]\label{prop:fixedpoint} 
For a $T$-fixed point $[(\Ecal,\phi_\Ecal)]\in \Mcal(\vec{u},n,j)^T$
the underlying sheaf is of the form $\Ecal=\imath_\ast (\Ical_Z)\otimes\Rcal^{\vec{u}}$, where $\Ical_Z$ is the ideal sheaf of a $T$-fixed zero-dimensional subscheme $Z$ of $X_k$.
\end{proposition}
\begin{remark}\label{rem:fixedpoints}
Let $[(\Ecal,\phi_\Ecal)]$ be a $T$-fixed point of
$\Mcal(\vec{u},n,j)$, with
$\Ecal=\imath_\ast(\Ical_Z)\otimes\Rcal^{\vec{u}}$. The $T$-fixed
subscheme $Z$ of $X_k$ of length $n$ is a disjoint union of $T$-fixed
subschemes $Z_i$ for $i=1, \ldots, k$ supported at the $T$-fixed
points $p_i$ with $\sum_{i=1}^k\,
\mathrm{length}_{p_i}(Z_i)=n$. Put
$n_i=\mathrm{length}_{p_i}(Z_i)$. Since $p_i$ is the $T$-fixed point
of the smooth affine toric surface $U_i\simeq \C^2$, as explained in
Section \ref{sec:C2} the $T$-fixed subscheme $Z_i\in \hilb{n_i}{U_i}$
corresponds to a Young tableau $Y^i$ of weight $\vert Y^i\vert= n_i$
for $i=1, \ldots, k$. Thus the $T$-fixed point $Z$ corresponds to a
$k$-tuple of Young tableaux $\vec{Y}=(Y^1, \ldots, Y^k)$ with
$\vert\vec{Y}\vert:=\sum_{i=1}^k\, \vert Y^i\vert=n$. Hence we can
parametrize the point $[(\Ecal, \phi_{\Ecal})]$ by the pair $(\vec{Y},
\vec{u}\, )$ which we call the \emph{combinatorial datum} of $[(\Ecal, \phi_{\Ecal})]$.
\end{remark}

Consider the $T$-equivariant cohomology of the moduli spaces $\Mcal(\vec{u},n,j)$ and set
\begin{equation}
\W_{\vec{u},n,j}:=H^\ast_T \big( \Mcal(\vec{u},n,j) \big)_{\mathrm{loc}}\ .
\end{equation}
We endow $\W_{\vec{u},n,j}$ with the nondegenerate $\C(\varepsilon_1, \varepsilon_2)$-valued bilinear form
\begin{equation}
\langle A,B\rangle_{\W_{\vec{u},n,j}}:=(-1)^n \, p_{\vec{u}, n,j}^!
\, \big(\imath_{\vec{u}, n,j}^!\big)^{-1} (A \cup B)\ ,
\end{equation}
where $p_{\vec{u}, n,j}$ is the projection from $\Mcal(\vec{u},n,j)$ to a point and $\imath_{\vec{u}, n,j}\colon \Mcal(\vec{u},n,j)^T \hookrightarrow \Mcal(\vec{u},n,j)$ is the inclusion of the fixed-point locus. Thus for $\vec{u}\in\Ufrak_j$ we define
\begin{equation}\label{eq:L0eigenspace}
\W_{\vec{u},j}:= \bigoplus_{n \geq0 } \, \W_{\vec{u},n,j} \ ,
\end{equation}
and the total equivariant cohomology
\begin{equation}
\W_j:=\bigoplus_{\vec{u}\in \Ufrak_j} \ \bigoplus_{n \geq0 }\,
\W_{\vec{u},n,j}
\end{equation}
which is an infinite-dimensional vector space over the field $\C(\varepsilon_1,
\varepsilon_2)$
endowed with the nondegenerate $\C(\varepsilon_1,
\varepsilon_2)$-valued bilinear form $\langle -,- \rangle_{\W_j}$
induced by the symmetric bilinear forms $\langle -,-\rangle_{\W_{\vec{u},n,j}}$.

Let us denote by $[\vec{Y}, \vec{u}\, ]$ the equivariant cohomology
class, defined similarly to \eqref{eq:fixed-class}, associated with the
$T$-fixed point $[(\Ecal, \phi_{\Ecal})]$ with combinatorial datum
$(\vec{Y}, \vec{u}\, )$. By the localization theorem, the classes
$[\vec{Y}, \vec{u}\, ]$ with $\vec{u}\in \Ufrak_j$ form a $\C(\varepsilon_1,\varepsilon_2)$-basis of $\W_j$.

\bigskip \section{Representations of $\glfrakhat_k$\label{sec:glkreps}}

\subsection{Overview}

The results collected so far imply the following result.
\begin{proposition}\label{cor:cohomologymodulihilbert}
There is an isomorphism
\begin{multline}\label{eq:isomWH}
\Psi_j \, \colon\,  \W_j \ \xrightarrow{ \ \sim \ } \
\bigoplus_{\vec{u}\in\Ufrak_j} \ \bigoplus_{ n\geq0 } \,
H^\ast_T\big(\hilb{n}{X_k} \big)_{\mathrm{loc}}\\[4pt]
\simeq \Big(\, \bigoplus_{n\geq0} \, H^\ast_T\big(\hilb{n}{X_k}
\big)_{\mathrm{loc}} \, \Big) \otimes \C(\varepsilon_1, \varepsilon_2)[\Ufrak_j]\\[4pt]
\xrightarrow{ \ \sim \ } \ \Big(\, \bigoplus_{n\geq0 } \,
H^\ast_T\big(\hilb{n}{X_k} \big)_{\mathrm{loc}} \, \Big) \otimes \C(\varepsilon_1, \varepsilon_2)[\Qfrak+\omega_j]\ ,
\end{multline}
where the first arrow is induced by the morphisms
$\tilde\imath_{\vec{u},n,j}^{\,\ast}$ while the last arrow is induced
by the map $\psi_j$ introduced in \eqref{eq:psij}. There is also an isomorphism
\begin{equation}
\Psi:=\bigoplus_{j=0}^{k-1}\, \Psi_j\, \colon \, \W \ \xrightarrow{\
  \sim \ } \ \bigoplus_{j=0}^{k-1}\ \bigoplus_{n\geq0 } \,
H^\ast_T\big(\hilb{n}{X_k} \big)_{\mathrm{loc}} \otimes \C(\varepsilon_1, \varepsilon_2)[\Qfrak+\omega_j]\ ,
\end{equation}
where $\W:=\bigoplus_{j=0}^{k-1}\, \W_j$.
\end{proposition}
In this section we first study the equivariant cohomology of
$\hilb{n}{X_k}$ and construct over it an action of the sum
(identifying central elements) $\hfrak_{\C(\varepsilon_1,\varepsilon_2)}\oplus
\hfrak_{\C(\varepsilon_1,\varepsilon_2),\Qfrak}$. Then we use the
Frenkel-Kac construction (Theorem \ref{thm:frenkelkac}) to obtain an
action of $\glfrakhat_k=\hfrak_{\C(\varepsilon_1,\varepsilon_2)}
\oplus\slfrakhat_k$ on $\W_j$ for $j=0, 1, \ldots, k-1$. 

\subsection{Equivariant cohomology of $\hilb{n}{X_k}$}

In this subsection we derive some results concerning the equivariant
cohomology of the Hilbert schemes $\hilb{n}{X_k}$ by generalizing similar results of \cite[Section 2]{art:qinwang2007} (see also \cite[Section 2]{art:maulikoblomov2009}).

As discussed in Remark \ref{rem:fixedpoints}, a $T$-fixed point
$Z\in\hilb{n}{X_k}$ corresponds to a $k$-tuple $(Z_1, \ldots, Z_k)$
where $Z_i$ is a $T$-fixed point of $\hilb{n_i}{U_i}$ for $i=1,\ldots,
k$ with $\sum_{i=1}^k\, n_i=n$, or equivalently to a $k$-tuple $\vec
Y= (Y^1,
\ldots, Y^k)$ of Young tableaux with $|\vec Y|:= \sum_{i=1}^k\, \vert Y^i\vert = n$. The following result is straightforward to prove.
\begin{lemma}\label{lem:decomposition}
Let $Z$ be a $T$-fixed point of $\hilb{n}{X_k}$. Then there is a $T$-equivariant isomorphism
\begin{equation}
T_Z\hilb{n}{X_k}\simeq \bigoplus_{i=1}^k \, T_{Z_i}\hilb{n_i}{U_i}\ ,
\end{equation}
where $Z=\bigsqcup_{i=1}^k\, Z_i$ and $n_i$ is the length of $Z_i$ at $p_i$ for $i=1, \ldots, k$.
\end{lemma}

By Lemma \ref{lem:decomposition} we get
\begin{equation}
\mathrm{ch}_T\big(T_Z\hilb{n}{X_k} \big)=\sum_{i=1}^k\, \mathrm{ch}_T
\big(T_{Z_i}\hilb{n_i}{U_i} \big)\ .
\end{equation}
By using the description \eqref{eq:coordinatering-i} of the coordinate
ring $\C[U_i]$ of $U_i$, one computes the equivariant Chern characters
\begin{equation}
\mathrm{ch}_T\big(T_{Z_i}\hilb{n_i}{U_i} \big)=\sum_{s\in Y^i}\,
\big(\e^{(L(s)+1)\, \varepsilon_1^{(i)}-A(s)\,
  \varepsilon_2^{(i)}}+\e^{-L(s)\, \varepsilon_1^{(i)}+(A(s)+1)\, \varepsilon_2^{(i)}}\big)\ .
\end{equation}
From now on we identify a torus-fixed point $Z$ of $\hilb{n}{X_k}$ with its $k$-tuple $\vec{Y}$ of Young tableaux. 

Let $\vec{Y}=(Y^1, \ldots, Y^k)$ be a torus-fixed point. Define
\begin{align}
\eu_{+}(\vec{Y})&:=\prod_{i=1}^k \ \prod_{s\in Y^i} \,
\Big( \big(L(s)+1\big)\, \varepsilon_1^{(i)}-A(s)\, \varepsilon_2^{(i)}\Big)\ ,\\[4pt]
\eu_{-}(\vec{Y})&:=\prod_{i=1}^k \ \prod_{s\in Y^i}\,
\Big(L(s)\, \varepsilon_1^{(i)}-\big(A(s)+1\big)\, \varepsilon_2^{(i)}\Big)\ .
\end{align}
Then the equivariant Euler class of the tangent bundle at the fixed
point $\vec{Y}$ is given by
\begin{equation}
\eu_T\big(T_{\vec{Y}}\hilb{n}{X_k}\big)=(-1)^n \,
\eu_{+}(\vec{Y})\, \eu_{-}(\vec{Y})\ .
\end{equation}

\subsubsection{Equivariant basis I: \ Torus-fixed points}\label{sec:basisfixedpoints}

Let $\vec{Y}$ be a $k$-tuple of Young tableaux corresponding to a
fixed point in $\hilb{n}{X_k}$ and $[\vec{Y}]$ the equivariant
cohomology class defined similarly to \eqref{eq:fixed-class}. By the projection formula we get
\begin{equation}
[\vec{Y}]\cup[\vec{Y}'\,]=\delta_{\vec{Y},\vec{Y}'}\,
\eu_T\big(T_{\vec{Y}}\hilb{n}{X_k} \big)[Y]=(-1)^n\,
\delta_{\vec{Y},\vec{Y}'}\, \eu_+(\vec{Y})\, \eu_-(\vec{Y})[\vec{Y}]\ .
\end{equation}
Denote 
\begin{equation}
\imath_n:=\bigoplus_{\vec{Y}\in\hilb{n}{X_k}^T}\, \imath_{\vec{Y}}\,
\colon \, \hilb{n}{X_k}^T\ \longrightarrow \ \hilb{n}{X_k}\ .
\end{equation}
In analogy to Equation \eqref{eq:bilinearform-C2}, define the bilinear form 
\begin{align}\label{eq:bilinearform}
\langle -,- \rangle_{\mathbb{H}_n} \, \colon \, \mathbb{H}_n\times
\mathbb{H}_n &\ \longrightarrow \  \C(\varepsilon_1, \varepsilon_2) \
,\\ \nonumber
(A, B)&\ \longmapsto \ (-1)^n\, p_n^!\, \big(\imath_n^! \big)^{-1}(A \cup B)
\end{align}
where $\mathbb{H}_n:=H_T^{\ast}(\hilb{n}{X_k})_{\mathrm{loc}}$.

As in Section \ref{sec:C2}, for any class $[\vec{Y}]\in H_T^{4n}(\hilb{n}{X_k})$ we define a distinguished class
\begin{equation}
[\alpha_{\vec{Y}}]:=\frac{[\vec{Y}]}{\eu_+(\vec{Y})} \ \in
\ H_T^{2n}\big(\hilb{n}{X_k} \big)_{\mathrm{loc}}\ .
\end{equation}
Then by the same computation as in Equation \eqref{eq:product-C2} we get
\begin{equation}\label{eq:product}
\big\langle [\alpha_{\vec{Y}}]\,,\, [\alpha_{\vec{Y}'}]
\big\rangle_{\mathbb{H}_n} =  \delta_{\vec{Y},\vec{Y}'} \,
\frac{\eu_-(\vec{Y})}{\eu_+(\vec{Y})} =
\delta_{\vec{Y},\vec{Y}'} \ \prod_{i=1}^k \ \prod_{s\in Y^i} \, \frac{
  L(s)\, \beta_i+A(s)+1 }{ \big(L(s)+1 \big)\, \beta_i+A(s) }  \ ,
\end{equation}
where analogously to \eqref{eq:beta} we defined
\begin{equation}\label{eq:betai}
\beta_i:=-\frac{\varepsilon_1^{(i)}}{\varepsilon_2^{(i)}}\ .
\end{equation}
Note that when $n=1$, $\vec{Y}$ is just a fixed point $p_i\in X_k^T$ with $i=1, \ldots, k$. Thus we have
\begin{equation}
\eu_{+}(p_i) =\varepsilon_1^{(i)}=(k-i+1)\,
\varepsilon_1-(i-1)\, \varepsilon_2 \qquad \mbox{and} \qquad 
\eu_{-}(p_i)=-\varepsilon_2^{(i)}=(k-i)\,
\varepsilon_1-i\, \varepsilon_2\ ,
\end{equation}
and therefore
\begin{equation}
\beta_i=\frac{\eu_+(p_i)}{\eu_-(p_i)}\ .
\end{equation}
If for $i=1, \ldots, k$ we define $[\alpha_i]:=[\alpha_{p_i}]$, then we get
\begin{equation}
\big\langle [\alpha_i]\,,\, [\alpha_j] \big\rangle_{\mathbb{H}_1} =
\beta^{-1}_i \, \delta_{ij} \ \in \ \C(\varepsilon_1, \varepsilon_2)\ .
\end{equation}
By the localization theorem and Equation \eqref{eq:product}, the
classes $[\alpha_{\vec{Y}}]$ with $\vert \vec{Y}\vert=n$ form a
$\C(\varepsilon_1,\varepsilon_2)$-linear basis of
$\mathbb{H}_n$. Hence the bilinear form \eqref{eq:bilinearform} is
nondegenerate; it extends to give a nondegenerate symmetric bilinear
form $\langle -,- \rangle_{\mathbb{H}}$ on the total equivariant cohomology $\mathbb{H}:=\bigoplus_{n\geq0}\, \HH_n$ of the Hilbert schemes of points on $X_k$. 

\begin{remark}\label{rem:description-Ui}
Let $i\in \{1, \ldots, k\}$. By the localization theorem, the
$\C(\varepsilon_1,\varepsilon_2)$-linear subspace of $\mathbb{H}$
generated by all classes $[\vec{Y}]$ associated to fixed points
$\vec{Y}=(Y^1, \ldots, Y^k)$ such that $Y^j=\emptyset$ for every $j\in
\{1, \ldots, k\}$ with $j\neq i$ is isomorphic to
\begin{equation}\label{eq:cohomologyUi}
\bigoplus_{m\geq 0} \, H_T^\ast\big(\hilb{m}{U_i} \big)\localizedi \ .
\end{equation}
\end{remark}
Note that $\C\big[\varepsilon_1^{(i)},\varepsilon_2^{(i)}
\big]=\C[\varepsilon_1,\varepsilon_2]$ and
$\C\big(\varepsilon_1^{(i)},\varepsilon_2^{(i)}
\big)=\C(\varepsilon_1,\varepsilon_2)$. Analogously to what we did for
$\C^2$, we can thus define
\begin{equation}
\HH_{U_i,m} := H_T^\ast\big(\hilb{m}{U_i}
\big)_{\mathrm{loc}}\qquad\mbox{and}\qquad\HH_{U_i} := \bigoplus_{m\geq 0}
\, \HH_{U_i,m}\ .
\end{equation}
By the localization theorem, there exists a $\C(\varepsilon_1,\varepsilon_2)$-linear isomorphism
\begin{equation}\label{eq:omega}
\Omega\, \colon \, \mathbb{H} \ \xrightarrow{ \ \sim \ } \
\bigotimes_{i=1}^k\, \HH_{U_i} \ .
\end{equation}
In particular, for a fixed point $\vec{Y}=(Y^1, \ldots, Y^k)$ we have
\begin{equation}
\Omega\, \colon \, [\alpha_{\vec{Y}}] \ \longmapsto \ [\alpha_{Y^1}]\otimes \cdots \otimes [\alpha_{Y^k}]\ .
\end{equation}
The isomorphism $\Omega$ interwines the bilinear forms $\langle -,-
\rangle_{\mathbb{H}}$ and $\prod_{i=1}^k\, \langle -,- \rangle_i$,
where $\langle -,- \rangle_i$ is the symmetric bilinear form on
$\HH_{U_i}$ defined analogously to \eqref{eq:bilinearform-C2}.
In a similar way, there is a $\C(\varepsilon_1,\varepsilon_2)$-linear isomorphism
\begin{equation}\label{eq:omegak}
\Omega_k\, \colon\, \HH_1 \ \xrightarrow{ \ \sim \ } \
\bigoplus_{i=1}^k \, \HH_{U_i,1}\ .
\end{equation}
In this case $\displaystyle\Omega_k\colon [\alpha_i]\mapsto (0,
\ldots, [\alpha_i], \ldots, 0)$, where the class $[\alpha_i]$ on the
left-hand side belongs to $\HH_1=H_T^\ast(X_k)_{\mathrm{loc}}$ while on the right-hand
side it belongs to $\HH_{U_i,1}$ as defined in Section
\ref{sec:C2}. The isomorphism $\Omega_k$ also intertwines the symmetric bilinear forms.

\subsubsection{Equivariant basis II: \ Torus-invariant divisors}\label{sec:divisors}

Let $[D_i]_T$ be the class in $\HH_1=H_T^\ast(X_k)_{\mathrm{loc}}$ given by the
$T$-invariant divisor $D_i$ for $i=0,1, \ldots, k$. For
$i=1,\dots,k-1$, we have
\begin{equation}\label{eq:divisor}
[D_i]_T=\frac{[p_i]}{\eu_T(T_{p_i}D_i)}+\frac{[p_{i+1}]}{\eu_T(T_{p_{i+1}}D_i)}=\frac{[p_i]}{\varepsilon_2^{(i)}}+\frac{[p_{i+1}]}{\varepsilon_1^{(i+1)}}=
-\beta_i\, [\alpha_i]+[\alpha_{i+1}]\ .
\end{equation}
Thus for $i,j =1, \ldots, k-1$ we obtain the pairings
\begin{equation}\label{eq:divisorcartan}
\big\langle[D_i]_T\,, \,[D_j]_T \big\rangle_{\mathbb{H}_1}=\left\{
\begin{array}{ll}
{2} \ ,& i=j\ ,\\
{-1} \ , & \vert i-j\vert =1\ ,\\
{0} \ , & \mbox{otherwise}\ .
\end{array}\right.
\end{equation}
By applying the localization theorem to $[D_0]_T$ and $[D_k]_T$ we
further obtain
\begin{equation}
\left[D_0\right]_T=\frac{[p_1]}{k\, \varepsilon_1}=\frac{[p_1]}{\varepsilon_1^{(1)}}=[\alpha_1]\qquad\mbox{and}\qquad
\left[D_k\right]_T=\frac{[p_k]}{k\,
  \varepsilon_2}=\frac{[p_k]}{\varepsilon_2^{(k)}}=-\beta_k\, [\alpha_k]\ .
\end{equation}
By using these expressions, one can straightforwardly obtain the pairings
\begin{equation*}
\big\langle [D_0]_T\,,\,[D_i]_T \big\rangle_{\mathbb{H}_1}=\left\{
\begin{array}{ll}
{\beta_1^{-1}} \ , & i=0\ ,\\
{-1} \ , & i=1\ ,\\
{0} \ , & \mbox{otherwise}
\end{array}\right.
\qquad\mbox{and}\qquad \big\langle [D_k]_T\,,\,[D_i]_T \big\rangle_{\mathbb{H}_1}=\left\{
\begin{array}{ll}
{\beta_k} \ , & i=k\ ,\\
{-1} \ , & i=k-1\ ,\\
{0} \ , & \mbox{otherwise}\ .
\end{array}\right.
\end{equation*}
Now we can relate the classes $[\alpha_i]$ for $i=1, \ldots, k$ to the
classes $[D_j]_T$ for $j=0,1, \ldots, k$. By using Equation
\eqref{eq:divisor}, for $i=2, \ldots, k$ one obtains
\begin{equation}\label{eq:alpha}
[\alpha_i]=\sum_{j=0}^{i-2}\, \Big(\, \prod_{s=j+1}^{i-1}\,
\beta_s\, \Big)\, [D_j]_T+[D_{i-1}]_T\ .
\end{equation}
Since $\eu_+(p_l)=\eu_-(p_{l-1})$ for $l=2,
\ldots, k$, we get 
\begin{equation}
 \prod_{s=j+1}^{i-1}\,
 \beta_s=\frac{\eu_+(p_{j+1})}{\eu_-(p_{i-1})} \
 . 
\end{equation}
By using the definition of $[\alpha_k]$ and Equation \eqref{eq:alpha} for $i=k$ we obtain
\begin{equation}
-\beta_k^{-1}\, [D_k]_T=[\alpha_k] =\sum_{j=0}^{k-1}\,
\frac{\eu_+(p_{j+1})}{\eu_-(p_{k-1})}\ [D_j]_T \
.
\end{equation}
If we formally put $\eu_+(p_{k+1}):=\eu_-(p_k)$,
we can reformulate this equation as
\begin{equation}\label{eq:sum}
\sum_{j=0}^{k}\, \eu_+(p_{j+1})\ [D_j]_T=0\ ,
\end{equation}
and in particular for all $i=0,1,\dots,k$ we have
\begin{equation}
[D_i]_T=-\sum_{\stackrel{\scriptstyle j=0} {\scriptstyle j\neq i}}^k\,
\frac{\eu_+(p_{j+1})}{\eu_+(p_{i+1})} \ [D_j]_T\ .
\end{equation}
\begin{remark}
If the action is antidiagonal, i.e., $t=t_1=t_2^{-1}$, Equation \eqref{eq:sum} implies that $\sum_{j=0}^k\, [D_j]_T=0$.
\end{remark}

As the classes $[\alpha_1], \ldots, [\alpha_k]$ form a
$\C(\varepsilon_1,\varepsilon_2)$-linear basis of $\mathbb{H}_1$, by
Equations \eqref{eq:alpha} and \eqref{eq:sum} the classes
\begin{equation}\label{eq:bases}
\big\{[D_0]_T, [D_1]_T, \ldots, [D_{k-1}]_T
\big\}\qquad\mbox{and}\qquad \big\{[D_1]_T, [D_2]_T, \ldots, [D_{k}]_T
\big\}
\end{equation}
are also $\C(\varepsilon_1,\varepsilon_2)$-linear bases for $\mathbb{H}_1$. Under the isomorphism $\Omega_k$ of Equation \eqref{eq:omegak}, we have 
\begin{equation}
\Omega_k\, \colon\, [D_i]_T \ \longmapsto \ -\beta_i \, \big(0,
\ldots,0, [\alpha_i], 0, \ldots, 0 \big)+\big(0, \ldots,0,
[\alpha_{i+1}], 0, \ldots, 0 \big)
\end{equation}
for $i=1, \ldots, k-1$, together with a similar description for $[D_0]_T$ and $[D_k]_T$.

\subsection{Heisenberg algebras}\label{sec:heisenbergalbebra}

Let $m$ be a positive integer and $Y$ a torus-invariant closed curve
in $X_k$. Define the correspondences
\begin{equation}
Y_{n,m}:=\big\{(Z, Z'\,)\in \hilb{n+m}{X_k}\times \hilb{n}{X_k}\
\big\vert\ Z'\subset Z \ , \
\mathrm{supp}(\mathcal{I}_{Z'}/\mathcal{I}_Z)=\{y\}\subset Y \big\} \ .
\end{equation}
Let $q_1$ and $q_2$ be the projections of $\hilb{n+m}{X_k}\times
\hilb{n}{X_k}$ to the two factors respectively. We define the linear
operator $\pfrak_{-m}([Y]_T)\colon \mathbb{H} \to \mathbb{H}$ which
acts on $A\in \mathbb{H}_n$ as
$\pfrak_{-m}([Y]_T)(A):=q_1^!\big(q_2^\ast(A)\cup
[Y_{n,m}]_T \big)\in \mathbb{H}_{n+m}$. This definition is well-posed
because the restriction of $q_1$ to $Y_{n,m}$ is proper. Since the
bilinear form $\langle -,- \rangle_{\mathbb{H}}$ is nondegenerate on
$\mathbb{H}$, we can define $\pfrak_m([Y]_T)$ to be the adjoint
operator of $\pfrak_{-m}([Y]_T)$. By using one of the two bases in \eqref{eq:bases}, we extend by linearity in $\alpha$ to obtain the linear operator $\pfrak_m(\alpha)$ for every $\alpha\in\mathbb{H}_1=H_T^\ast(X_k)_{\mathrm{loc}}$.

%For the following result we refer to \cite[Section 2.3]{art:qinwang2007} (antidiagonal torus action) and \cite[Section 2.2]{art:maulikoblomov2009}.
\begin{theorem}[{see \cite{art:qinwang2007,art:maulikoblomov2009}}]\label{thm:heisenberg}
The linear operators $\pfrak_m(\alpha)$, where $m\in\Z\setminus\{0\}$ and $\alpha\in{H}_T^\ast({X_k})_{\mathrm{loc}}$, satisfy the Heisenberg commutation relations
\begin{equation}
\big[\pfrak_m(\alpha),\pfrak_n(\beta) \big]= m\, \delta_{m,-n}\,
\langle\alpha,\beta\rangle_{\HH_1}\,\mathrm{id}\qquad \mbox{and}
\qquad\big[\pfrak_m(\alpha),\mathrm{id} \big]= 0\ .
\end{equation}
The vector space $\mathbb{H}$ is the Fock space of the Heisenberg algebra $\hfrak_{\mathbb{H}_1}$ modelled on $\mathbb{H}_1={H}_T^\ast({X_k})_{\mathrm{loc}}$ with highest weight vector the unit element $\vacuum$ in ${H}_T^0(\hilb{0}{X_k})_{\mathrm{loc}}$.
\end{theorem}

\subsubsection{Heisenberg algebra of rank $k$}

Let $i\in \{1, \ldots, k\}$. Consider the Heisenberg algebra $\hfrak_{i}$ over $\C(\varepsilon_1,\varepsilon_2)$ generated by the operators
\begin{equation}
\pfrak_{-m}^i:=\pfrak_{-m}([\alpha_i])\qquad\mbox{and}\qquad\mathfrak{p}_{m}^i:=\pfrak_{m}([\alpha_i])
\end{equation}
for $m\in\Z_{>0}$. By Theorem \ref{thm:heisenberg}, the commutation relations are 
\begin{equation}\label{eq:commutationalpha}
\big[\pfrak_m^i,\pfrak_n^j \big]=m\, \delta_{m,-n}\, \delta_{ij}\, 
\big\langle[\alpha_i]\,,\,[\alpha_i]
\big\rangle_{\mathbb{H}}\,\mathrm{id}=m\, \delta_{m,-n}\,
\delta_{ij}\, \beta_i^{-1}\,\mathrm{id}\ .
\end{equation}
Since $\{[\alpha_1], \ldots, [\alpha_k]\}$ is a $\C(\varepsilon_1,
\varepsilon_2)$-linear basis of $\HH_1$, the Heisenberg algebra $\hfrak_{\HH_1}$ is generated by $\pfrak_m^i$ for $i=1, \ldots, k$ and $m\in\Z\setminus\{0\}$.

Let $\HH_{U_i}$ be the $\C(\varepsilon_1,\varepsilon_2)$-linear
subspace of $\mathbb{H}$ introduced in Section
\ref{sec:basisfixedpoints}. Then by Theorem
\ref{thm:nakajimaoperators} $\HH_{U_i}$ is the Fock space for the
Heisenberg algebra $\hfrak_{i}$ for any $i\in \{1, \ldots, k\}$; therefore the
$\C(\varepsilon_1, \varepsilon_2)$-vector space $\HH_{U_i}$ is
generated by the elements $\pfrak_\lambda^i\vert 0 \rangle$ where
$\pfrak_\lambda^i:=\prod_{l\geq 1}\, (\pfrak_{-l}^i)^{m_l}$ for a
partition $\lambda=(1^{m_1}\, 2^{m_2}\, \cdots)$. One can show that
\begin{equation}
\big\langle\mathfrak{p}_\lambda^i\vacuum\,,\,\pfrak_\mu^i\vacuum
\big\rangle_{\HH_{U_i}} = \delta_{\lambda,\mu} \, z_\lambda \, \beta_i^{-\ell(\lambda)} \ .
\end{equation}

On the algebra $\Lambda_{\C(\varepsilon_1, \varepsilon_2)}$ of
symmetric functions over the field $\C(\varepsilon_1, \varepsilon_2)$
we introduce the Jack inner product \eqref{eq:jackinnerprod} with
parameter $\beta_i$.
We shall denote with $\Lambda_{\beta_i}$ the algebra $\Lambda_{\C(\varepsilon_1, \varepsilon_2)}$ endowed with the symmetric bilinear form $\langle -,- \rangle_{\beta_i}$. Thus by the isomorphism \eqref{eq:cohomologyUi} and Theorem \ref{thm:actionC2} there exists an isomorphism of $\C(\varepsilon_1, \varepsilon_2)$-vector spaces
\begin{equation}
\Phi_i\, \colon \, \HH_{U_i} \ \xrightarrow{ \ \sim \ } \
\Lambda_{\beta_i}\ , \qquad
\pfrak_\lambda^i\vert 0 \rangle \ \longmapsto \ p_\lambda\ ,
\end{equation}
which intertwines the symmetric bilinear forms $\langle -,- \rangle_i$
and  $\langle -,- \rangle_{\beta_i}$. For $m>0$ the operator
$\pfrak_{-m}^i$ acts as multiplication by $p_m$ on
$\Lambda_{\beta_i}$ while its adjoint $\pfrak_{m}^i$ with respect to
the symmetric bilinear form $\langle -,- \rangle_i $ acts as $m\,
\beta_i^{-1}\, \frac{\partial}{\partial p_m}$.

By Theorem \ref{thm:actionC2} we can also determine how $\Phi_i$ acts
on the $\C(\varepsilon_1, \varepsilon_2)$-linear basis
$\{[\alpha_{\vec{Y}}] \}$ of $\HH_{U_i}$, where
$\vec{Y}=(Y^1, \ldots, Y^k)$ is a fixed point such that $Y^j=\emptyset$ for
every $j\in \{1, \ldots, k\}$ with $j\neq i$.
\begin{proposition}
Let $\vec{Y}=(Y^1, \ldots, Y^k)$ be a fixed point such that $Y^j=\emptyset$
for every $j\in \{1, \ldots, k\}$ with $j\neq i$. Then
\begin{equation}
\Phi_i([\alpha_{\vec{Y}}])=J_{\lambda_i}(x;\beta_i^{-1}) \ ,
\end{equation}
where $Y_{\lambda_i}:= Y^i$.
\end{proposition}
Define $\Lambda_{\vec{\beta}}=\bigotimes_{i=1}^k\,
\Lambda_{\beta_i}$ endowed with the symmetric bilinear form $\langle
p,q\rangle_{\vec{\beta}} :=\prod_{i=1}^k\, \langle
p_i,q_i\rangle_{\beta_i} $ for $p=p_1\otimes \cdots\otimes p_k$ and
$q=q_1\otimes \cdots\otimes q_k$ in $\Lambda_{\vec{\beta}}\, $. 
For a $k$-tuple of Young tableaux $\vec{Y}$, define in
$\Ucal(\hfrak_{\HH_1}) $ the operators
$\pfrak_{\vec{Y}}=\prod_{i=1}^k \, \pfrak^i_{\lambda_{i}}$. We have thus
proven the following result. 
\begin{theorem}\label{thm:cohomologyALE}
There exists a $\C(\varepsilon_1, \varepsilon_2)$-linear isomorphism
\begin{equation}
\Phi:= \bigotimes_{i=1}^k\, \Phi_i \, \colon \, \HH \ \longrightarrow
\ \Lambda_{\vec{\beta}} 
\end{equation}
preserving bilinear forms such that
\begin{equation}
\Phi \left( \pfrak_{\vec{Y}} \vacuum \right) = p_{\lambda_1} \otimes
\cdots \otimes p_{\lambda_k} \qquad \mbox{and} \qquad \Phi \left(
  [\alpha_{\vec Y}] \right) = J_{\lambda_1}(x; \beta_1^{-1}) \otimes \cdots \otimes J_{\lambda_k}(x; \beta_k^{-1})\ .
\end{equation}
Via the isomorphism $\Phi$, the operators $\pfrak_m^i$ act on
$\Lambda_{\vec{\beta}}$ as multiplication by $p_{-m}$ on the $i$-th
factor for $m<0$ and as the derivation $m\, \beta_i^{-1}\, \frac{\partial}{\partial p_m}$ on the $i$-th factor for $m>0$.
\end{theorem}

\subsubsection{Lattice Heisenberg algebra of type $A_{k-1}$}

Let us now define 
\begin{equation}
\mathfrak{q}^i_{-m}:=\pfrak_{-m}([D_i]_T)\qquad\mbox{and}\qquad\mathfrak{q}^i_m:=\pfrak_m([D_i]_T)
\end{equation}
for $m\in \Z_{>0}$ and $i=1, \ldots, k-1$. By Equation \eqref{eq:divisorcartan} the operators $\mathfrak{q}^i_m$ satisfy the commutation relations
\begin{equation}
[\mathfrak{q}^i_m, \mathfrak{q}^j_n] =m\, \delta_{m,-n} \, C_{ij}\,
\mathrm{id}\qquad \mbox{for} \quad i,j=1, \ldots, k-1 \ , \ m,n\in\Z\setminus\{0\}\ ,
\end{equation}
where $C=(C_{ij})$ is the Cartan matrix of the Dynkin diagram of type $A_{k-1}$.
Let $\Lfrak\subset H_T^\ast(X_k)_{\mathrm{loc}}$ be the $\Z$-lattice generated by the classes $[D_1]_T, \ldots, [D_{k-1}]_T$ with the symmetric bilinear form given by the Cartan matrix $C$. Then the lattice Heisenberg algebra $\hfrak_{\C(\varepsilon_1, \varepsilon_2),\Lfrak}$ associated with $\Lfrak$ over $\C(\varepsilon_1, \varepsilon_2)$, which has generators $\mathfrak{q}^i_m$ for $m\in \Z\setminus\{0\}$ and $i=1, \ldots, k-1$, is isomorphic to the Heisenberg algebra $\hfrak_{\C(\varepsilon_1, \varepsilon_2), \Qfrak}$ of type $A_{k-1}$ over $\C(\varepsilon_1, \varepsilon_2)$ (cf. Example \ref{ex:latticeheisenebrg}). 

Let
\begin{equation}\label{eq:E}
E:=\sum_{i=0}^k \,a_i \, [D_i]_T
\end{equation}
where $a_i\in
\C(\varepsilon_1,\varepsilon_2)$ with $i=0,1, \ldots, k$ satisfy the relations
\begin{equation}\label{eq:relation1}
2a_j-a_{j-1}-a_{j+1}= 0\ , \quad j=1,\ldots, k-1 \qquad \mbox{and} \qquad
a_0 \, \varepsilon_2 + a_k \, \varepsilon_1 \neq 0\ .
\end{equation}
The first condition ensures that $\langle [D_i]_T, E\rangle_{\HH_1} =0$ for
$i=1, \ldots, k-1$ while the second condition implies that $\{[D_1]_T,
\ldots, [D_{k-1}]_T, E\}$ is a $\C(\varepsilon_1,
\varepsilon_2)$-linear basis of $H_T^\ast(X_k)_{\mathrm{loc}}$. By \eqref{eq:relation1} one has
\begin{equation}
\kappa:= \langle E, E\rangle_{\mathbb{H}_1} =  a_0^2\,
\beta_1^{-1}-a_0\, a_1-a_k\, a_{k-1}+a_k^2\, \beta_k \ .
\end{equation}
From now on we set $a_i=i$ in Equation \eqref{eq:E} for $i=0, 1,
\ldots, k$, which is consistent with the conditions in Equation
\eqref{eq:relation1}. This implies that $\kappa=k\, \beta$. In the
following we normalize the equivariant cohomology class $E$ such that $\langle E,
E\rangle_{\mathbb{H}_1}=1$; we denote the normalized class with the
same symbol.

Define $\pfrak_{-m}:=\pfrak_{-m}(E)$ and $\pfrak_m:=\pfrak_m(E)$ for
$m\in \Z_{>0}$. Then the operators $\mathfrak{q}^i_m$ and $\pfrak_m$ satisfy the commutation relations
\begin{equation*}
\left\{
\begin{array}{ll}
\big[\mathfrak{q}^i_m, \mathfrak{q}^j_n\big]=m\, \delta_{m,-n} \,
C_{ij}\, \mathrm{id} & \mbox{for} \quad i,j=1, \ldots, k-1 \ , \ m,n\in\Z\setminus\{0\}\ ,\\[4pt]
\left[\mathfrak{q}^i_m, \pfrak_n\right]=0 & \mbox{for} \quad i=1,
\ldots, k-1 \ , \ m,n\in\Z\setminus\{0\}\ ,\\[4pt]
\big[\pfrak_m, \pfrak_n\big]=m\, \delta_{m,-n}\, \mathrm{id} &
\mbox{for} \quad m,n\in\Z\setminus\{0\}\ .
\end{array}\right.
\end{equation*}
Let $\Lfrak'\subset H_T^\ast(X_k)_{\mathrm{loc}}$ be the $\Z$-lattice generated by
the classes $[D_1]_T, \ldots, [D_{k-1}]_T$ and $ E$. Then the operators
$\mathfrak{q}^i_m$ and $\pfrak_n$ for $m,n\in\Z\setminus\{0\}$ and
$1\leq i \leq k-1$ define the lattice Heisenberg algebra
$\hfrak_{\C(\varepsilon_1, \varepsilon_2),\Lfrak'}$ associated with
$\Lfrak'$ over $\C(\varepsilon_1, \varepsilon_2)$. In particular,
$\hfrak_{\C(\varepsilon_1, \varepsilon_2),\Lfrak'}$ is the sum
(identifying central elements) of, respectively, the Heisenberg
algebra $\hfrak_{\C(\varepsilon_1, \varepsilon_2),\Qfrak}$ of type
$A_{k-1}$ over $\C(\varepsilon_1, \varepsilon_2)$ and the Heisenberg
algebra $\hfrak_{\C(\varepsilon_1, \varepsilon_2)}$ over
$\C(\varepsilon_1, \varepsilon_2)$ generated by $\pfrak_m$ for
$m\in\Z\setminus\{0\}$.
Since $\{[D_1]_T, \ldots, [D_{k-1}]_T\}\cup \{ E\}$ is a $\C(\varepsilon_1,
\varepsilon_2)$-linear basis of $H_T^\ast(X_k)_{\mathrm{loc}}$, we have $\hfrak_{ \C(\varepsilon_1, \varepsilon_2),\Lfrak'}\cong \hfrak_{\mathbb{H}_1}$. Hence $\mathbb{H}$ is the Fock space of $\hfrak_{\C(\varepsilon_1, \varepsilon_2),\Lfrak'}$.
In what follows we omit the symbol $\C(\varepsilon_1, \varepsilon_2)$ from the notation for the Heisenberg algebras generated over the field $\C(\varepsilon_1,\varepsilon_2)$.

\begin{remark}\label{rem:divisors}
By Equation \eqref{eq:divisor}, for $l=1, \ldots, k-1$ one has
\begin{align}\label{eq:fixedpoint-divisor-1}
\mathfrak{q}^l_m &=-\beta_l\, \mathfrak{p}^l_m+\pfrak^{l+1}_m\ ,\\[4pt] \label{eq:fixedpoint-divisor-2}
 \pfrak_m &= \sqrt{-k\, \varepsilon_1\, \varepsilon_2} \ \sum_{i=1}^k\,
\frac{1}{\varepsilon_2^{(i)}}\, \pfrak_m^i=\frac{1}{\sqrt{\big\langle
    [X_k], [X_k] \big\rangle_{\HH_1}}}\,\sum_{i=1}^k\,\frac{1}{\varepsilon_2^{(i)}}\, \pfrak_m^i\ .
\end{align}
\end{remark}

\subsection{Dominant representation of $\glfrakhat_k$ on $\W_j$}

In the following we omit the dependence of the Lie algebras on the field $\C(\varepsilon_1, \varepsilon_2)$ to simplify the presentation.

\begin{proposition}\label{prop:representation}
Let $j\in\{0,1, \ldots, k-1\}$. There is an action of $\glfrakhat_k$ on $\W_j$ under which $\W_j$ is the $j$-th dominant representation of $\glfrakhat_k$ at level one, i.e., the highest weight representation of $\glfrakhat_k$ with fundamental weight $\widehat{\omega}_j$ of type $\widehat{A}_{k-1}$.
\end{proposition}
\proof
The vector space $\mathbb{H}$ is an irreducible highest weight
representation of the sum (identifying central elements)
$\hfrak \oplus
\hfrak_{\Qfrak}$ of the Heisenberg
algebra $\hfrak$ and the lattice
Heisenberg algebra $\hfrak_{\Qfrak}$
of type $A_{k-1}$ generated over the field $\C(\varepsilon_1,\varepsilon_2)$, respectively, by the operators $\pfrak_m$ and $\mathfrak{q}^i_m$ for
$m\in\Z\setminus\{0\}$ and $i=1, \ldots, k-1$. We apply the Frenkel-Kac construction (Theorem \ref{thm:frenkelkac}) to the representation $\hfrak_{\Qfrak} \to \End (\HH )$ to obtain a level one representation 
\begin{equation}
\slfrakhat_k \ \longrightarrow \ \End
\big(\HH\otimes\C(\varepsilon_1,\varepsilon_2)[\Qfrak+\omega_j] \big) \ .
\end{equation}
We can extend the representation of $\hfrak$ from $\HH$ to
$\HH\otimes\C(\varepsilon_1,\varepsilon_2)[\Qfrak+\omega_j]$ by
letting it act as the identity on the group algebra of
$\Qfrak+\omega_j$. Thus we get a level one representation of
$\glfrakhat_k$ with
\begin{equation}
\glfrakhat_k \ \longrightarrow \ \End
\big(\HH\otimes\C(\varepsilon_1,\varepsilon_2)[\Qfrak+\omega_j] \big) \ .
\end{equation}
Thanks to the isomorphism $\Psi_j$ introduced in \eqref{eq:isomWH},
this gives a level one representation of $\glfrakhat_k$ with
\begin{equation}
\glfrakhat_k \ \longrightarrow \ \End(\W_j)\ .
\end{equation}
Since $\HH$ is the Fock space of
$\hfrak \oplus\hfrak_{\Qfrak}$, the
module $\HH\otimes\C(\varepsilon_1,\varepsilon_2)[\Qfrak+\omega_j]$ is isomorphic to the $j$-th dominant representation $ \Fcal_{\C(\varepsilon_1,\varepsilon_2)}\otimes \Vcal(\,\widehat{\omega}_j\,)$ of $\glfrakhat_k$ (cf.\  Theorem \ref{thm:frenkelkac}). Hence to complete the proof it is enough to note that the class $[\emptyset, \boldsymbol e_j]$ corresponding to the fixed point $(\Rcal_j, \phi_{\Rcal_j})$ in $\Mcal(\boldsymbol e_j, 0,j)$ is sent via $\Psi_j$ to $\vacuum\otimes [\omega_j]$, which is the highest weight vector of $\HH\otimes\C(\varepsilon_1,\varepsilon_2)[\Qfrak+\omega_j]$.
\endproof
\begin{remark}
Proposition \ref{prop:representation} is analogous to a previous result derived for Nakajima quiver varieties (see e.g. \cite[Section 10]{art:nakajima:1994-3} and \cite[Section 5.1]{art:smirnov2013}).
\end{remark}
Let us introduce the Virasoro operators of $\glfrakhat_k$ by (cf.\ Sections \ref{sec:virasoro-heisenberg} and \ref{sec:virasoro-affine})
\begin{align}
L_0&:=L_0^{\hfrak}+L_0^{\widehat{\mathfrak{sl}}_{k}}=\sum_{m=
  1}^\infty\, \pfrak_{-m}\,\pfrak_m +\sum_{i=1}^{k-1}\ \sum_{m=
  1}^\infty\, \qfrak_{-m}^{\eta_i}\,\qfrak_{m}^{\eta_i}+\frac{1}{2}\,
\sum_{i=1}^{k-1} \, \big(\qfrak_0^{\eta_i}\big)^2\ ,\\[4pt]
L_n&:=L_n^{\hfrak}+L_n^{\widehat{\mathfrak{sl}}_{k}}=\frac12\,\sum_{m\in\Z}
\, \pfrak_{-m}\,\pfrak_{m+n} +\frac12\,\sum_{i=1}^{k-1}\
\sum_{m\in\Z}\, \qfrak_{-m}^{\eta_i}\,\qfrak_{m+n}^{\eta_i}\qquad
\mbox{for}\quad n\neq 0\ ,
\end{align}
where $\{\eta_i\}_{i=1}^{ k-1}$ is an orthonormal basis of the vector
space $\Qfrak\otimes_{\Z}\R$ and we set $\pfrak_0:=0$. Note that
$\{\eta_i\}_{i=1}^{k-1} \cup \{E\}$ is an orthonormal basis of the
vector space $\C(\varepsilon_1,\varepsilon_2)^k\simeq H_T^\ast(X_k)_{\mathrm{loc}}$, so after an orthonormal change of basis and a suitable normalization one can rewrite the operators $L_n$ in the form
\begin{align}
L_0&=\sum_{l=1}^k\, \beta_l \ \sum_{m= 1}^\infty \,
\pfrak^l_{-m}\,\pfrak^l_m +\frac{1}{2}\, \sum_{i=1}^{k-1} \, \big(\qfrak_0^{\eta_i}\big)^2\ ,\\[4pt]
L_n&=\frac12\,\sum_{l=1}^k\, \beta_l \ \sum_{\stackrel{\scriptstyle
    m\in\Z}{\scriptstyle m\neq 0, -n}}\, \pfrak^l_{-m}\,\pfrak^l_{m+n}
+\sum_{i=1}^{k-1}\,
\qfrak_{0}^{\eta_i}\,\qfrak_{n}^{\eta_i} \qquad
\mbox{for}\quad n\neq 0\ .
\end{align}
\begin{proposition}\label{prop:L0eigen}
Let $j\in\{0,1, \ldots, k-1\}$, $\vec{u}\in\Ufrak_j$ and $n\in\N$. Then
\begin{equation}
{L_0} \big\vert_{\W_{\vec{u}, n ,j}}=\big(n+\mbox{$\frac12$}\,
\vec{v}\cdot C\vec{v}\, \big)\,\operatorname{id}_{\W_{\vec{u}, n ,j}}\ ,
\end{equation}
where $\vec{v}:=C^{-1}\vec{u}$.
\end{proposition}
\proof
By Equation \eqref{eq:integralofmotions-fixedpoints}, we have
$\sum_{m=1}^\infty \, \pfrak^l_{-m}\,\pfrak^l_m\triangleright
[\vec{Y}]=\beta_l^{-1}\, \vert Y^l \vert \, [\vec{Y}]$ for $[\vec{Y}]=[(\emptyset,
\ldots, Y^l,\ldots, \emptyset)]$ in $\HH_{U_l}$ and $l\in\{1,\ldots,k\}$. Then by using the isomorphism $\Omega$ introduced in \eqref{eq:omega} and the isomorphism $\Psi_j$ introduced in \eqref{eq:isomWH} we get
\begin{equation}
\Big(\, \sum_{l=1}^k\, \beta_l \ \sum_{m= 1}^\infty\,
\pfrak^l_{-m}\,\pfrak^l_m\, \Big)\Big\vert_{ \W_{\vec{u}, n ,j}}=n \,\operatorname{id}_{\W_{\vec{u}, n ,j}}\ .
\end{equation}
On the other hand, since $\{\eta_i\}_{i=1}^{ k-1}$ is an orthonormal basis of $\Qfrak\otimes_{\Z}\R$ we get
\begin{equation}
\Big(\, \sum_{i=1}^{k-1}\, \big(\qfrak_0^{\eta_i}\big)^2\,
\Big)\Big\vert_{ \W_{\vec{u}, n ,j}}=\Big(\, \sum_{i=1}^{k-1}\,
\Big\langle\eta_i\,,\, \sum_{s=1}^{k-1}\, v_s \, \gamma_s
\Big\rangle_{\Qfrak\otimes_\Z\R}^2 \, \Big) \,\operatorname{id}_{\W_{\vec{u}, n ,j}}
=\big(\,\vec{v}\cdot C\vec{v} \, \big) \, \operatorname{id}_{\W_{\vec{u}, n ,j}}\ ,
\end{equation}
where $\vec{v}:=C^{-1}\vec{u}$.
\endproof

Since
\begin{equation}
L_0\triangleright [\emptyset, \vec{u}\, ]=\mbox{$\frac12$}\,
{\vec{v}\cdot C\vec{v}} \, [\emptyset, \vec{u}\, ] \qquad\mbox{and}\qquad L_n\triangleright
[\emptyset,\vec{u}\, ] = 0 \quad \mbox{for}\quad n>0\ ,
\end{equation}
we have the following result.
\begin{corollary}\label{cor:hwrep}
$\W_{\vec{u},j}$ is a highest weight representation
of the Virasoro algebra $\virfrak$ associated with $\glfrakhat_k$, which is generated by the operators $L_n$ and $\cfrak$ with highest
weight vector $[\emptyset, \vec{u}\, ]$ and conformal dimension
\begin{equation}
\Delta_{\vec{u}}:=\mbox{$\frac12$} \, {\vec{v}\cdot C\vec{v}} \ . 
\end{equation}
Moreover, the energy eigenspace decomposition of the
representation $\W_{\vec{u},j}$ is given by \eqref{eq:L0eigenspace}.
\end{corollary}

\begin{proposition}\label{prop:Wkweight}
Let $j\in\{0,1, \ldots, k-1\}$. The weight decomposition of $\W_j$ as
a $\glfrakhat_k$-module is given by
\begin{equation}
\W_j=\bigoplus_{\vec{u}\in \Ufrak_j} \, \W_{\vec{u},j}\ .
\end{equation}
\end{proposition}
\proof
For $j=0, 1, \ldots, k-1$ and for any element
$A\otimes [\gamma_{\vec{u}}+\omega_j]\in \HH\otimes \mathbb{C}(\varepsilon_1, \varepsilon_2)[\Qfrak+\omega_j]$ with $\vec{u}\in \Ufrak_j$, we have 
\begin{align}
h_0\triangleright (A\otimes [\gamma_{\vec{u}}+\omega_j ])
&=\Big(1-\sum_{l=1}^{k-1}\, u_l\Big)\, (A\otimes
[\gamma_{\vec{u}}+\omega_j] ) \ , \\[4pt]
h_i\triangleright (A\otimes [\gamma_{\vec{u}}+\omega_j ])&=u_i  \,
(A\otimes [\gamma_{\vec{u}}+\omega_j ] )\qquad
\mbox{for} \quad i=1, \ldots, k-1\ .
\end{align}
Under the $\glfrakhat_k$-action, $\W_{\vec{u},j}$ decomposes as
\begin{equation}
\W_{\vec{u},j}\simeq \Fcal_{\C(\varepsilon_1,\varepsilon_2)}\otimes\Vcal(\widehat{\omega}_j)_{\gamma_{\vec{u}}}\ ,
\end{equation}
where 
$\Vcal(\widehat{\omega}_j)_{\lambda}:=\{w \in
\Vcal(\widehat{\omega}_j)\,\vert\, h_i\triangleright
w=(\widehat{\omega}_j+\lambda)(h_i)\,w \}$ for any weight $\lambda$.
The assertion now follows by showing that for a weight $\lambda$, the
vector space $\Vcal(\widehat{\omega}_j)_\lambda$ is nonzero if and only if
$\lambda=\gamma_{\vec{u}}$ for some $\vec{u}\in\Ufrak_j$. For this,
let $\vec{\xi}^{\
  \vec{v}}:=\prod_{i=1}^{k-1} \, \xi_i^{v_i}$ for $\xi_i\in\C^*$ with
$|\xi_i|<1$ and
any vector $\vec{v}=(v_1,\dots,v_{k-1})$, and set $\vec{h}:=(h_1, \ldots, h_{k-1})$. Then
it is enough to note that the trace of the operator $\qsf^{L_0}\
\vec{\xi}^{\ C^{-1}\vec{h}}$ is given by
\begin{equation}
\Trr_{\W_j} \, \qsf^{L_0}\, \vec{\xi}^{\ C^{-1}\vec{h}} =
\sum_{\vec{u}\in\Ufrak_j} \ \sum_{n=0}^\infty\, \qsf^{n+\frac12\,
  \vec{v}\cdot C\vec{v}} \ \vec{\xi}^{\ \vec{v}}
\end{equation}
and the right-hand side is exactly the character of the $j$-th dominant
representation of $\glfrakhat_k$ by \cite[Section 5.3]{art:bruzzopedrinisalaszabo2013}.
\endproof

\subsubsection{Whittaker vectors}

Consider now the completed total equivariant cohomology
\begin{equation}
\widehat{\W}_j:=\prod_{\vec{u}\in\Ufrak_j} \ \prod_{n\geq0}\, \W_{\vec{u},n,j}\ .
\end{equation}
We can extend the isomorphism \eqref{eq:isomWH} to
\begin{equation}
\widehat{\Psi}_j\ \colon \ \widehat{\W}_j \
\xrightarrow{ \ \sim \ } \ \widehat{\HH} \otimes \Big(\,
\prod_{\gamma_{\vec{u}}\in\Qfrak}\, \C(\varepsilon_1,\varepsilon_2)\,
(\gamma_{\vec{u}}+\omega_j) \, \Big) \ ,
\end{equation}
where $\widehat{\HH}:=\prod_{n\geq0}\, H^\ast_T(\hilb{n}{X_k})_{\mathrm{loc}}$ is
the completed total equivariant cohomology of the Hilbert schemes of
points on $X_k$. In the following we drop the explicit symbols
$\widehat{\Psi}_j$ from the notation in order to simplify the
presentation, and we denote
\begin{equation}
|\omega_j\rangle:= \sum_{\vec{u}\in\Ufrak_j}
\, [\emptyset, \vec{u}\, ] \ .
\end{equation}

\begin{proposition}\label{prop:whittakerALE}
Fix $\vec{\eta}\in\C(\varepsilon_1,\varepsilon_2)^k$. In the completed total equivariant cohomology $\widehat{\W}_j$ the vector
\begin{equation}
G_j(\vec{\eta}\,):=\exp\Big(\, \sum_{i=1}^k \, \eta_i \
\pfrak_{-1}^i\, \Big) |\omega_j\rangle
\end{equation}
is a Whittaker vector of type $\chi_{\vec\eta}\, $, where
$\chi_{\vec\eta}\colon
\Ucal(\hfrak^+ \oplus
\hfrak_{\Qfrak}^+ )\to\C(\varepsilon_1,\varepsilon_2)$ is defined by
\begin{align}
\chi_{\vec\eta}(\qfrak_1^i)&= \eta_{i+1} \, \beta_{i+1}^{-1} -
\eta_i \qquad \mbox{and} \qquad \chi_{\vec\eta}(\qfrak_m^i) = 0
\quad \mbox{for} \quad m>1 \ , \ i=1,\ldots, k-1\ ,\\[4pt]
\chi_{\vec\eta}(\pfrak_1)&= \frac{1}{\sqrt{\big\langle [X_k], [X_k] \big\rangle_{\HH_1}}}\,\sum_{i=1}^k\,
\frac{\eta_i}{\varepsilon_1^{(i)}} \qquad  \mbox{and} \qquad
\chi_{\vec\eta}(\pfrak_m)= 0 \quad \mbox{for} \quad m>1 \ .
\end{align}
\end{proposition}
\proof
Let $\widehat{\HH}_{U_i}:=\prod_{n\geq0}\, H_T^\ast(\hilb{n}{U_i})_{\mathrm{loc}}$
be the completed total equivariant cohomology of the Hilbert scheme $\hilb{n}{U_i}$ for
$i=1, \ldots, k$, and define $G(\eta_i):=\exp\left(\, \eta_i \
  \pfrak_{-1}^i\, \right) \vacuum \in \widehat{\HH}_{U_i}$. By using
Theorem \ref{thm:cohomologyALE} and the completed versions of the
isomorphisms \eqref{eq:omega} and \eqref{eq:isomWH}, we can rewrite
the vector $G_j(\vec{\eta}\, )$ as
\begin{equation}
G_j(\vec{\eta}\, ) = G(\eta_1)\otimes\cdots\otimes G(\eta_k) \ \otimes
\ \sum_{\vec{u}\in\Ufrak_j}\, (\gamma_{\vec{u}}+\omega_j) \ .
\end{equation}
By Proposition \ref{prop:whittakerC2}, $G(\eta_i)$ for $i=1,\ldots, k$ is a Whittaker vector for the Heisenberg algebra $\hfrak_i$ of type $\chi_i$, where
\begin{equation}\label{eq:characterheisenberg}
\chi_i(\pfrak_1^i):= \eta_i \, \beta_i^{-1}\qquad\mbox{and}\qquad \chi_i(\pfrak_m^i):= 0\quad \mbox{for}\quad m>1\ .
\end{equation}
Again by Theorem \ref{thm:cohomologyALE}, each $\hfrak_i$ acts
trivially on $\HH_{U_l}$ for $l\neq i$ and it is easy to see that
$G_j(\vec{\eta}\, )$ is a Whittaker vector for the Heisenberg algebra
$\hfrak_{\HH_1}$ of type
$\tilde \chi_{\vec\eta}$,
where $\tilde \chi_{\vec\eta} \colon
\Ucal\big(\hfrak_{\HH_1}^+\big)\to\C(\varepsilon_1,\varepsilon_2)$
is defined by $\tilde \chi_{\vec\eta}(\pfrak_m^i):=\chi_i(\pfrak_m^i)$ for
$i=1, \ldots, k$ and $m\in\Z\setminus \{0\}$. Then by Remark
\ref{rem:divisors}, $G_j(\vec{\eta}\, )$ is a Whittaker vector for
$\glfrakhat_k$ of type $\chi_{\vec\eta}$ with $\chi_{\vec\eta}\colon
\Ucal(\hfrak^+ \oplus
\hfrak_{\Qfrak}^+)\to \C(\varepsilon_1,\varepsilon_2)$ defined for every $m>0$ by
\begin{align}
\chi_{\vec\eta}(\qfrak_m^i) & :=\tilde \chi_{\vec\eta}(\pfrak_m^{i+1})
- \beta_i \, \tilde \chi_{\vec\eta}(\pfrak_m^i) = \delta_{m,1} \,
\big(\eta_{i+1} \, \beta_{i+1}^{-1} - \eta_i \big)\ ,\\[4pt]
\chi_{\vec\eta}(\pfrak_m) &: = \sqrt{-k\, \varepsilon_1\, \varepsilon_2}\ \sum_{i=1}^k \, 
\frac{1}{\varepsilon_2^{(i)}} \, \tilde \chi_{\vec\eta}(\pfrak_m^i) = \delta_{m,1} \,
\sqrt{-k\, \varepsilon_1\, \varepsilon_2} \ \sum_{i=1}^k\,
\frac{\eta_i}{\varepsilon_1^{(i)}}\ .
\end{align}
\endproof

\bigskip \section{Chiral vertex operators for $\glfrakhat_k$\label{sec:chiralvertex}}

\subsection{Ext-bundles and bifundamental hypermultiplets\label{se:bifundXk}}

In this section we construct and study the natural generalizations of the
Ext vertex operators from Section
\ref{sec:carlsson-okounkov} for the moduli spaces
$\Mcal(\vec{u},n,j)$. 

For $n\in\N$, $j\in\{0,1,\ldots,k-1\}$ and $\vec u\in\Ufrak_j$, let $\boldsymbol{\Ecal}_{\vec u,n,j}$ denote the universal sheaf on
$\Mcal(\vec{u},n,j)\times \Xscr_k$. Define
\begin{equation}
\boldsymbol{\Ecal}_i:=p_{i 3}^\ast\big(\boldsymbol{\Ecal}_{\vec
  u_i,n_i,j_i} \big) \ \in \
K\big(\Mcal(\vec{u}_1,n_1;j_1)
\times\Mcal(\vec{u}_2,n_2;j_2)\times \Xscr_k \big) \qquad \mbox{for}
\quad i=1, 2 \ , 
\end{equation}
where $p_{ij}$ is the projection of $\Mcal(\vec{u}_1,n_1;j_1)
\times\Mcal(\vec{u}_2,n_2;j_2)\times \Xscr_k$ onto the product of the $i$-th and $
j$-th factors. Denote by $p_3$ the
projection of the same product onto $\Xscr_k$.

Let $T_\mu=\C^\ast$ and $H^\ast_{T_\mu}(\mathrm{pt};\C)=\C[\mu]$. Denote 
by $\Ocal_{\Xscr_k}(\mu)$ the trivial line bundle on $\Xscr_k$ on which $T_\mu$ acts by scaling the fibers. 
\begin{definition}[{\cite[Definition 4.17]{art:bruzzopedrinisalaszabo2013}}]
The \emph{Carlsson-Okounkov bundle} is the element
\begin{equation}
\Ebf_\mu^{\vec{u}_1,n_1,j_1;\vec{u}_2,n_2,j_2}:={p_{12}}_\ast\big(-\boldsymbol{\Ecal
}_1^\vee\cdot \boldsymbol{\Ecal}_2\cdot p_3^\ast(\Ocal_{\Xscr_k}(\mu)\otimes\Ocal_{\Xscr_k}(-\Dscr_\infty))
\big)
\end{equation}
in the K-theory $K\big(\Mcal(\vec{u}_1,n_1;j_1)\times\Mcal(\vec{u}_2,n_2;j_2)\big)$.
\end{definition}
By \cite[Section 4.5]{art:bruzzopedrinisalaszabo2013}, the fibre of
 $\Ebf_\mu^{\vec{u}_1,n_1,j_1;\vec{u}_2,n_2,j_2}$ over
 $\big([(\Ecal,\phi_\Ecal)]\,,\, [(\Ecal',\phi_{\Ecal'})]\big)$ in
 $\Mcal(\vec{u}_1,n_1;j_1)\times\Mcal(\vec{u}_2,n_2;j_2)$ is given by
\begin{equation}
\Ebf_\mu ^{\vec{u}_1
,n_1,j_1;\vec{u}_2,n_2,j_2}\big|_{\left([(\Ecal,\phi_\Ecal)]\,,\,[(\Ecal',\phi_{\Ecal'})]\right)}=\Ext^1\big(\Ecal,\Ecal'\otimes
\Ocal_{\Xscr_k}(\mu)\otimes\Ocal_{\Xscr_k}(-\Dscr_\infty) \big)\ .
\end{equation}
One can compute the dimension of this vector space by a straightforward
generalization of the dimension computations of \cite[Appendix
A]{art:bruzzopedrinisalaszabo2013} to get the rank
\begin{equation}
\rk\big(\Ebf_\mu^{\vec{u}_1,n_1,j_1;\vec{u}_2,n_2,j_2} \big) =
n_1+n_2+\mbox{$\frac12$}\, \vec v_{21}\cdot C\vec
v_{21}-\mbox{$\frac1{2k}$}\, j_{21}\, (k- j_{21} ) \ ,
\end{equation}
where $\vec v_{21}:=C^{-1}(\vec u_2-\vec u_1)$ and $j_{21}\in\{0, 1,\ldots, k-1\}$ is the equivalence class modulo $k$ of $j_2-j_1$.

Let $\W:=\bigoplus_{j=0}^{k-1}\, \W_j$ endowed with the nondegenerate $\C(\varepsilon_1,
\varepsilon_2)$-valued bilinear form $\langle -,- \rangle_{\W}$
induced by the symmetric bilinear forms $\langle
-,-\rangle_{\W_j}$. Define the operator $\V_\mu(\vec x,z)\in\mathrm{En
d}(\W\, )[[
z^{\pm\, 1},x_1^{\pm\,1},$ $\dots,x_{k-1}^{\pm\,1}]]$ by its matrix elements
\begin{multline}\label{eq:vertexoperator}
(-1)^{n_2} \, \big\langle \V_\mu(\vec x,z)A_1\,,\, A_2
\big\rangle_{\W}
:=z^{n_2-n_1 +\Delta_{\vec u_2}-\Delta_{\vec u_1} } \, \vec x\,^{\vec
v_{21}}\\ \times \ 
\int_{\Mcal(\vec{u}_1,n_1;j_1)\times\Mcal(\vec{u}_2,n_2;j_2)} \,
\eu_T\big(\Ebf_\mu^{\vec{u}_1,n_1,j_1;\vec{u}_2,n_2,j_2}
\big) \cup p_1^\ast(A_1) \cup p_2^\ast(A_2) \ ,
\end{multline} 
where $A_i\in H_T^\ast(\Mcal(\vec{u}_i,n_i;j_i))_{\mathrm{loc}}$ and $p_i$ is the projection from $\Mcal(\vec{u}_1,n_1;j_1)\times\Mcal(\vec{u}_2,n_2;j_2)$ to the $i$-th factor for
$i=1,2$. The extra isospin parameters $\vec x:=(x_1,\dots,x_{k-1})$
weigh the $\slfrakhat_k$-action, and we abbreviated $\vec x\,^{\vec v}:= \prod_{i=1}^{k-1}\, x_i^{v_i}$ for a
vector $\vec v=(v_1,\dots,v_{k-1})$. By the computations of \cite[Section 4.7]{art:bruzzopedrinisalaszabo2013}, the matrix elements \eqref{eq:vertexoperator} in the fixed point basis are given by
\begin{multline}\label{eq:vertexoperator-fixedbasis} 
\big\langle \V_\mu(\vec x,z)[\vec{Y}_1,\vec{u}_1\,] \,,\, [\vec{Y}_2,\vec{u}_2\, ]
\big\rangle_{\W}
\\ \shoveleft{=(-1)^{\vert \vec{Y}_2\vert}\,
z^{\vert \vec{Y}_2\vert-\vert \vec{Y}_1\vert+\Delta_{\vec u_2}-\Delta_{\vec u_1} } \, \vec x\,^{\vec v_{21}} \ 
\eu_T\big(\Ebf_\mu^{\vec{u}_1,n_1,j_1;\vec{u}_2,n_2,j_2} \big\vert_{\left([(\Ecal_1,\phi_{\Ecal_1})]\,,\,[(\Ecal_2,\phi_{\Ecal_2})]\right)}
\big) } \\[4pt]
\shoveleft{ = (-1)^{\vert \vec{Y}_2\vert}\,
z^{\vert \vec{Y}_2\vert-\vert \vec{Y}_1\vert+\Delta_{\vec u_2}-\Delta_{\vec u_1} } \, \vec x\,^{\vec v_{21}} \ \prod_{i=1}^k \, m_{Y_1^{i},
 {Y_2^{i}}} \big(\varepsilon_1^{(i)},\varepsilon_2^{(i)},
\mu-(\vec{v}_{21})_{i}\,
\varepsilon_1^{(i)}-(\vec{v}_{21})_{i-1}\, \varepsilon_2^{(i)} \big) }
\\ \times \ \prod_{n=
1}^{k-1}\,
\ell^{(n)}_{\vec{v}_{21}}
\big(\varepsilon_1^{(n)},\varepsilon_2^{(n)},\mu \big)\ ,
\end{multline}
where $[(\Ecal_1,\phi_{\Ecal_1})]$ and $[(\Ecal_2,\phi_{\Ecal_2})]$
are the $T$-fixed points corresponding respectively to the
combinatorial data $(\vec{Y}_1,\vec{u}_1)$ and
$(\vec{Y}_2,\vec{u}_2)$, and we use the convention
$(\vec{v}_{21})_{0}=(\vec{v}_{21})_{k}=0$. Here $m_{Y_1, Y_2}$ is
defined in \eqref{eq:m}, while $\ell^{(n)}_{\vec{v}}$ is the \emph{edge 
contribution} defined in Appendix \ref{app:edgecontributions}. This factorized expression for
the matrix elements represents the contribution of the $U(1)\times U(1)$ bifundamental
hypermultiplet for $\Ncal=2$ quiver gauge theories on the ALE space $X_k$.

\subsection{Vertex operators and primary fields}

In this subsection we factorize the operators $\V_\mu(\vec x,z)$ defined in
\eqref{eq:vertexoperator} under the decomposition $\glfrakhat_k=\hfrak\oplus\slfrakhat_k$ as tensor products of generalized bosonic
exponentials associated to the Heisenberg algebra
$\hfrak$ from Definition \ref{def:bosonops} with primary
fields of the Virasoro algebra associated to the affine Lie algebra
$\slfrakhat_k$ from Section \ref{sec:virasoro-affine}.

For $l=1,2$ fix $j_l\in\{0,1, \ldots, k-1\}$ and respectively $\vec{u}_l\in\Ufrak_{j_l}$. Set
\begin{equation}
\gamma_{21}:=\sum_{i=1}^{k-1}\,(\vec{v}_{21})_{i}\,\gamma_i=\psi_{j_2}(\vec{u}_2
)-\psi_{j_1}(\vec{u}_1) \ \in \ \Qfrak\otimes_{\Z}\Q\ .
\end{equation}
Note that $\gamma_{21}=\gamma_{\vec{u}_2}-\gamma_{\vec{u}_1}+\omega_{j_2}-\omega
_{j_1}$. We define the maps
\begin{align}
\exp\big(\gamma_{21}\,\big)\,\colon\, \HH\otimes\C(\varepsilon_1,\varepsilon_2)
[\Qfrak+\omega_{j_1}]\ &\longrightarrow \ \HH\otimes\C(\varepsilon_1,\varepsilon_2)[\Qfrak+
\omega_{j_2}]\ ,\\
v\otimes [\beta+\omega_{j_1}] \ &\longmapsto\ v\otimes [\beta+\gamma_{\vec{u}_2}-\gamma_{\vec{u}_1}+\omega_{j_2}]\ ,
\end{align}
and $\exp\big(\log z\ \cfrak+\gamma_{21}\big)\,\colon\, \HH\otimes\C(\varepsilon_1,\varepsilon_2)
[\Qfrak+\omega_{j_1}]\ \to \ \HH\otimes\C(\varepsilon_1,\varepsilon_2)[\Qfrak+
\omega_{j_2}]$, given by
\begin{equation}
\exp\big(\log z\ \cfrak+\gamma_{21}\,\big)\triangleright(v \otimes
[\beta+\omega_{j_1}]):=z^{\frac12\, \langle
  \gamma_{21},\gamma_{21}\rangle_{\Qfrak\otimes_{\Z}\Q} + \langle
  \gamma_{21}, \beta+\omega_{j_1}\rangle_{\Qfrak\otimes_{\Z}\Q} }\, (v \otimes [\beta+\gamma_{\vec{u}_2}
-\gamma_{\vec{u}_1}+\omega_{j_2}])\ .
\end{equation}
Note that the operator $\exp\big(\log z\
\cfrak-\gamma_{21}\big)\,\exp\big(\gamma_{21}\big)\in\End(\HH\otimes
\C(\varepsilon_1,\varepsilon_2)[\Qfrak+\omega_{j_1}]) [[z,z^{-1}]]$ is
given by
\begin{multline}\label{eq:expopWj}
\exp\big(\log z\
\cfrak-\gamma_{21}\,\big)\,\exp\big(\gamma_{21}\,\big)\triangleright
(v\otimes [\beta +\omega_{j_1}]) \\ =z^{-\frac12\,\langle\gamma_{21},\gamma_{21}
\rangle_{\Qfrak\otimes_{\Z}\Q} -\langle\gamma_{21},
\beta+\omega_{j_1}\rangle_{\Qfrak\otimes_{\Z}\Q} }\, (v\otimes [\beta +\omega
_{j_1}]) \ . 
\end{multline}
In the following we shall suppress the explicit isomorphism $\Psi$ from Theorem
\ref{cor:cohomologymodulihilbert} in our formulas in order to
simplify the presentation.

We will now rewrite the operator $\V_\mu(\vec x, z)$ in terms of chiral vertex operators in
$\Hom(\W_{\vec{u}_1;j_1}, \W_{\vec{u}_2;j_2})$ $[[z^{\pm\,1},x_1^{\pm\,1},\dots,x_{k-1}^{\pm\,1}]]$ between two highest
weight representations of the Virasoro algebra associated with
$\slfrakhat_k$. For this, let us define the vertex operator ${\bar \V}_\mu(\vec{v}_{21}, \vec x, z)$ of
$\Hom(\W_{j_1}, \W_{j_2})[[z^{\pm\, 1},x_1^{\pm\,1},$ $\dots,x_{k-1}^{\pm\,1}]]$ by
\begin{equation}\label{eq:barVmudef}
{\bar \V}_\mu(\vec{v}_{21}, \vec x, z):=\vec x\,^{\vec v_{21}} \, \prod_{l=1}^{k-1}\,
\ell^{(l)}_{\vec{v}_{21}}
\big(\varepsilon_1^{(l)},\varepsilon_2^{(l)},\mu \big)\
\V_{1,-1}^{\gamma_{21}}(z) \ \exp\big(\log z\ \cfrak+\gamma_{21}\,\big)\ ,
\end{equation}
where $\V_{1,-1}^{\gamma_{21}}(z)$ is the normal-ordered bosonic exponential associated with the Heisenberg algebra $\hfrak_\Qfrak$.
\begin{theorem}\label{thm:virasoroprimary}
Under the decomposition $\glfrakhat_k=\hfrak\oplus\slfrakhat_k$, the operator $\V_\mu(\vec x, z)$ is given in terms of products of vertex
operators as
\begin{multline}
\V_\mu(\vec x, z)=\V_{-\frac{\mu}{\sqrt{-k\, \varepsilon_1\, \varepsilon
_2}},\frac{\mu+\varepsilon_1+\varepsilon_2}{\sqrt{-k\, \varepsilon_1\,
\varepsilon_2}}
}(z) \\ \otimes \ \sum_{j_1,j_2=0}^{k-1} \
\sum_{\vec{u}_1\in\Ufrak_{j_1},\vec{u}_2\in\Ufrak_{j_2}}\
\bar \V_\mu(\vec{v}_{21},\vec x, z)\, z^{\Delta_{\vec u_2}-\Delta_{\vec u_1}}\ \exp\big(\log z\ \cfrak-\gamma_{21}\,\big)\,\exp\big(\gamma_{21}\,\big)\big\vert_{ \W_{\vec{u}_1, j_1}}
\ ,
\end{multline}
where $\bar \V_\mu(\vec{v}_{21},\vec x, z)$ is a primary field of the Virasoro algebra
 generated by $L_n^{\slfrakhat_k}$ and $\cfrak$ with conformal
 dimension $\Delta_{\vec u_2-\vec u_1}=\frac12\, \vec{v}_{21}\cdot C\vec{v}_{21}$, i.e., for any $n\in \Z$ we have
\begin{equation}\label{eq:primary}
\big[L^{\slfrakhat_k}_n,\bar \V_\mu(\vec{v}_{21},\vec x, z)\big]=z^n \,
\big( z\, \partial_z + \mbox{$\frac12$}\, \vec{v}_{21}\cdot
C\vec{v}_{21}\, n \big)\bar\V_\mu(\vec{v}_{21},\vec x, z)\ .
\end{equation}
\end{theorem}
\proof
By using the isomorphism $\Psi$ and Equation \eqref{eq:carlokfixedbasis} we get 
\begin{multline}\label{eq:factorization}
\V_\mu(\vec x, z)=\Psi^{-1} \circ \Big(\, \sum_{j_1,j_2=0}^{k-1} \
\sum_{\vec{u}_1\in\Ufrak_{j_1},\vec{u}_2\in\Ufrak_{j_2}} \,z^{\Delta_{\vec u_2}-\Delta_{\vec u_1}}\, \vec x\,^{\vec v_{21}} \ \prod_{n
=1}^{k-1}\,
\ell^{(n)}_{\vec{v}_{21}}
\big(\varepsilon_1^{(n)},\varepsilon_2^{(n)},\mu \big)\\
\times \ : \, \prod_{i=1}^k  \ \V\big(\Ocal_{U_i}\big(\mu-(\vec{v}_{21})_{i}\,
\varepsilon_1^{(i)}-(\vec{v}_{21})_{i-1}\,
\varepsilon_2^{(i)}\big)\,,\,z\big) \, : \
 \otimes \ \big(\psi_{j_2}(\vec{u}_2)\otimes
 \psi_{j_1}(\vec{u}_1)^\ast \big) \, \Big) \circ \Psi
\end{multline}
where $\Ocal_{U_i}(\mu)$ is the trivial line bundle on $U_i\simeq\C^2$
with an action of $T_\mu$ which rescales the fibers, and $\psi_{j_1}(\vec{u}_1)^\ast$ denotes the dual vector to
$\psi_{j_1}(\vec u_1)$ in the dual vector space $\C(\varepsilon_1,
\varepsilon_2)[\Qfrak+\omega_{j_1}]^\ast$. By Theorem \ref{thm:CO} we
get an expression determined by the operators $\pfrak_m^i$ for
$m\in\Z\setminus\{0\}$ and $i=1,\dots,k$ as
\begin{multline}
\V\big(\Ocal_{U_i}\big(\mu-(\vec{v}_{21})_{i}\,
\varepsilon_1^{(i)}-(\vec{v}_{21})_{i-1}\,
\varepsilon_2^{(i)}\big)\,,\,z\big) \\ = \exp
\Big(-\frac{\mu-(\vec{v}_{21})_{i}\,
\varepsilon_1^{(i)}-(\vec{v}_{21})_{i-1}\, \varepsilon_2^{(i)}}{\varepsilon_2^{(
i)}}\,
\sum_{m=1}^\infty\, \frac{z^m}{m}\, \pfrak_{-m}^i\Big)\\
\times \ \exp\Big(\,
\frac{\mu+ \varepsilon_1+\varepsilon_2-(\vec{v}_{21})_{i}\,
\varepsilon_1^{(i)}-(\vec{v}_{21})_{i-1}\, \varepsilon_2^{(i)}}{\varepsilon_2^{(
i)}}\,
\sum_{m=1}^\infty\, \frac{z^{-m}}{m}\, \pfrak_{m}^i\, \Big) \ .    
\end{multline}
By using Equations \eqref{eq:fixedpoint-divisor-1} and
\eqref{eq:fixedpoint-divisor-2}, we can rewrite $: \, \prod_{i=1}^k\  \V\big(\Ocal_{U_i}\big(\mu-(\vec{v}_{21})
_{i}\,
\varepsilon_1^{(i)}-(\vec{v}_{21})_{i-1}\, \varepsilon_2^{(i)}\big)
\,,\, z\big) \, :$ in terms of Heisenberg operators $\qfrak_m^l$ and $\pfrak_m$, for $l=1,\ldots, k-1$ 
and $m\in\Z\setminus\{0\}$, and the first assertion now follows. The proof of Equation \eqref{eq:primary} is somewhat lengthy and can be found in Appendix \ref{app:virasoroprimary}.
\endproof
\begin{remark}
In the following we will denote by
$\V_\mu^{j_1,j_2}(\vec x, z)$ the restriction of the vertex operator
$\V_\mu(\vec x, z)$ to $\Hom(\W_{j_1}, \W_{j_2})[[z^{\pm\,1},x_1^{\pm\,1},\dots,x_{k-1}^{\pm\,1}]]$.
\end{remark}

\subsection{Integrals of motion\label{sec:integrals}}

Let $\Vbf^{\vec{u},n,j}$ be the pushforward of
$\Ebf_0^{\vec{u},n,j;\vec{u},n,j}$ with respect to the projection of
the product $\Mcal(\vec{u},n,j)\times\Mcal(\vec{u},n,j)$ to the second
factor. It is a $T$-equivariant vector bundle on $\Mcal(\vec{u},n,j)$
of rank $n+\frac{1}{2}\, (\vec{v}\cdot C\vec{v}-\frac{1}{k}\, j\, (k-j))$,
which we shall call the natural bundle over $\Mcal(\vec{u},n,j)$.
The $T$-equivariant Chern character of $\Vbf^{\vec{u},n,j}$ at a fixed point
$[(\Ecal,\phi_\Ecal)]$ with combinatorial datum $(\vec Y,\vec u\,)$ is given by (cf.\ \cite[Section 4.7]{art:bruzzopedrinisalaszabo2013})
\begin{equation}
\ch_T\big(\Vbf^{\vec{u},n,j}\big\vert_{[(\Ecal,\phi_\Ecal)]}\big) = \sum_{i=1}^k \ \sum_{s\in
  Y^i}\, \e^{-(v_i+L'_{Y^i}(s))\, \varepsilon_1^{(i)}- (v_{i-1}+
  A_{Y^i}'(s))\, \varepsilon_2^{(i)}}+\sum_{n=1}^{k-1}\,
L_{\vec{v}}^{(n)} \big(\varepsilon_1^{(n)}, \varepsilon_2^{(n)}\big) \ ,
\end{equation}
where the edge contributions $L_{\vec{v}}^{(n)}
\big(\varepsilon_1^{(n)}, \varepsilon_2^{(n)}\big) $ are defined in
Appendix \ref{app:edgecontributions}.

Let us denote by $\Vbf^j$ the natural bundle over
$\coprod_{\vec{u}\in\Ufrak_j}\ \coprod_{n\geq 0}\,
\Mcal(\vec{u},n,j)$, and consider the operators of multiplication by $\boldsymbol{I}_1:=\rk\big(\Vbf^j \big)$ and $\boldsymbol{I}_p:=
(\cc_{p-1} )_T\big(\Vbf^j \big)$ for $p\geq 2$ on
$\prod_{\vec{u}\in\Ufrak_j}\ \prod_{n\geq
  0}\,\W_{\vec{u}, n ,j}$. For example, one has
\begin{align}
\boldsymbol{I}_1\triangleright [\vec{Y}, \vec{u}\, ]= & \ \big(\, \vert
\vec{Y}\vert+\mbox{$\frac{1}{2}$} \,
\vec{v}\cdot C\vec{v}\, \big)\, [\vec{Y}, \vec{u}\, ] -\mbox{$\frac{1}{2k}\, j\, (k-j)$}\,  [\vec{Y},
\vec{u}\, ] \ , \\[4pt]
\boldsymbol{I}_2\triangleright [\vec{Y}, \vec{u}\, ]=& \
-\sum_{i=1}^k \ \sum_{s\in
  Y^i}\, \Big( \big(v_i+L'_{Y^i}(s) \big)\, \varepsilon_1^{(i)}+
\big(v_{i-1}+ A_{Y^i}'(s) \big)\, \varepsilon_2^{(i)}\Big)\, [\vec{Y},
\vec{u}\, ]\\
  & \ +\, \sum_{n=1}^{k-1}\,
  \ell_{\vec{v}}^{(n)}\big(\varepsilon_1^{(n)},
  \varepsilon_2^{(n)}\big)_{[1]}  \, [\vec{Y}, \vec{u}\, ] \ ,
\end{align}
where the quantities $\ell_{\vec{v}}^{(n)}(\varepsilon_1^{(n)},
  \varepsilon_2^{(n)})_{[1]} $ are defined in Appendix \ref{app:edgecontributions}. Note that, by Proposition \ref{prop:L0eigen}, the operator $\boldsymbol{I}_1$ coincides (up to a constant shift) with the Virasoro operator $L_0$ for $\glfrakhat_k$.
By using the description in Section \ref{sec:integralC2}, these operators can be written partly in terms of the Heisenberg
operators $\pfrak_m^i$ of $\hfrak_{\HH_1}$ and the $\slfrakhat_k$ generators $\qfrak_0^i=h_i$; one has
\begin{align}
\boldsymbol{I}_1   =& \
\sum_{i=1}^k\,\beta_i \ 
\sum_{m=1}^\infty \, \pfrak_{-m}^i \, \pfrak_{m}^i +\frac{1}{2}\, \sum_{i=1}^{k-1} \, \big(\qfrak_0^{\eta_i}\big)^2
- \frac{1}{2k}\, j\, (k-j)\, \operatorname{id} \ , \\[4pt] 
\boldsymbol{I}_2 =& \
\sum_{i=1}^k\,\varepsilon_1^{(i)}\, \Big(\, \frac{\beta_i}{2}\, \sum_{m,n=1}^\infty \, \big(
\pfrak_{-m}^i\, \pfrak_{-n}^i\, \pfrak_{m+n}^i+\pfrak_{-m-n}^i\, \pfrak_{n}^i\,
\pfrak_{m}^i \big) \\ & \qquad \qquad -\, \frac{\beta_i-1}{2}\, \sum_{m=1}^\infty\, (m-1)\,
\pfrak_{-m}^i\, \pfrak_{m}^i\, \Big)\\
& \ +\, \sum_{i=1}^k\, \varepsilon_1^{(i)} \ \sum_{j=1}^{k-1}\, \Big( \big(C^{-1}\big)^{ij}\, \beta_i-\big(C^{-1}\big)^{i-1,j}\Big) \ \sum_{m=1}^\infty\, \pfrak_{-m}^i \, \pfrak_{m}^i\, \qfrak_0^j + \boldsymbol{L}_1 \ , 
\end{align}
where $\boldsymbol{L}_1$ is the operator defined by $\boldsymbol{L}_1\triangleright [\vec{Y}, \vec{u}\, ]:= \sum_{n=1}^{k-1}\,
  \ell_{\vec{v}}^{(n)}(\varepsilon_1^{(n)},\varepsilon_2^{(n)})_{[1]} \, [\vec{Y}, \vec{u}\, ]$ and we set $(C^{-1})^{0j}=0=(C^{-1})^{kj}$.

Following~\cite{art:belavinbershteinfeiginlitvinovtarnopolsky2011}, here we shall identify a quantum integrable system for each Heisenberg subalgebra $\hfrak$ of $\glfrakhat_k$. Then each integral of motion associated to the Heisenberg algebra $\hfrak_{\HH_1}$ is a sum of integrals of motion of $k$ non-interacting Calogero-Sutherland models from Section \ref{sec:integralC2} with prescribed couplings; in particular the Hamiltonian is given by $k$ copies of one-component bosonized Calogero-Sutherland Hamiltonians as $\sum_{i=1}^k\, \varepsilon_1^{(i)}\, \Box^{\beta_i^{-1}}$. This infinite system of commuting operators is diagonalized in the fixed-point basis $[\vec Y,\vec u\,]$. This simultaneous eigenbasis also factorizes the primary operators from Theorem \ref{thm:virasoroprimary}. The remaining $\vec v$-dependent parts of the eigenvalues are instead interpreted as particular matrix elements of our geometrically defined vertex operators $\V_\mu(\vec x,z)$ in highest weight vectors of $\glfrakhat_k$, as we discuss in Appendix \ref{app:edgecontributions}.

\begin{remark}
By Remark \ref{rem:quivervariety}, the moduli spaces $\Mcal(\vec{u},n,j)$ are Nakajima quiver varieties of type $\widehat{A}_{k-1}$. The descriptions of the integrable systems corresponding to quiver varieties are detailed in \cite{art:maulikokounkov2012,art:smirnov2013}.
\end{remark}

\bigskip \section{$\Ncal=2$ quiver gauge theories on $X_k$\label{sec:quivergaugeXk}}

\subsection{$\Ncal=2$ gauge theory\label{sec:N=2Xk}}

In this subsection we fix $j\in\{0,1, \ldots, k-1\}$ corresponding to a fixed holonomy at infinity.
The instanton partition function for the pure $\Ncal=2$ $U(1)$ gauge
theory on the ALE space $X_k$ is the generating function (cf.\ \cite[Section~5.1]{art:bruzzopedrinisalaszabo2013})
\begin{align}
\Zcal_{X_k}\big(\varepsilon_1,\varepsilon_2; \qsf, \vec{\xi} \ \big)_j
:= & \ \sum_{\stackrel{\scriptstyle \vec{v}\in\frac{1}{k}\, \Z^{k-1}
  }{\scriptstyle k\, v_{k-1}= j\bmod{k} }} \,
  \vec{\xi}^{\ \vec{v}} \ \sum_{n=0}^\infty \, \qsf^{n+\frac{1}{2}\,
  \vec{v}\cdot C\vec{v}} \ \int_{\Mcal(\vec{u},n,j)} \,
\big[\Mcal(\vec{u},n,j) \big]_T \\[4pt]
 = & \ \sum_{\vec{u}\in\Ufrak_j} \ \sum_{n=0}^\infty\,
 \qsf^{\frac{1}{2}\, \vec{u}\cdot C^{-1}\vec{u}}\ \vec{\xi}^{\
   C^{-1}\vec{u}}\, (-\qsf)^{n} \, \big\langle
 [\Mcal(\vec{u},n,j)]_T\,,\,[\Mcal(\vec{u},n,j)]_T \big\rangle_{\mathbb{W}_j}\ ,
\end{align}
where $\qsf\in\C^*$ with $|\qsf|<1$, and the fugacities
$\vec\xi:=(\xi_1,\dots,\xi_{k-1})\in(\C^*)^{k-1}$ with $|\xi_i|<1$ can
be interpreted as coordinates on the maximal torus of the Lie group $SL(k,\C)$.

In general, as described in \cite[Section
5.1]{art:bruzzopedrinisalaszabo2013}, the partition functions
factorize into products of the corresponding instanton partition
functions over the affine toric subsets $U_i\simeq \C^2$ of $X_k$ and
are weighted by the edge contributions $\ell^{(n)}_{\vec{v}}$ which
appear in the equivariant Euler classes of the Carlsson-Okounkov
bundle from Section \ref{se:bifundXk} (see Appendix \ref{app:edgecontributions}). The edge contributions for the rank one $\Ncal=2$ gauge theory
on $X_k$ are roughly speaking the equivariant Euler classes of the
vector space $H^1(\Xscr_k, \Ocal_{\Xscr_k}(-\Dscr_\infty))$, which are
zero by the computation of the rank of the natural bundle in
\cite[Appendix A]{art:bruzzopedrinisalaszabo2013}, hence the edge
contribution is always equal to one. We thus obtain a factorization in
terms of Nekrasov partition functions for the pure $\Ncal=2$ gauge
theory on $\R^4$ given by
\begin{equation}
\Zcal_{X_k}\big(\varepsilon_1,\varepsilon_2; \qsf, \vec{\xi} \ \big)_j =
\sum_{\vec{u}\in\Ufrak_j}\, \qsf^{\frac{1}{2}\, \vec{u}\cdot
  C^{-1}\vec{u}}\ \vec{\xi}^{\ C^{-1}\vec{u}} \ \prod_{i=1}^k\, 
\Zcal_{\C^2}\big(\varepsilon_1^{(i)},\varepsilon_2^{(i)}; \qsf \big)\ .
\end{equation}

Let us denote by
\begin{equation}\label{eq:characterdef}
\chi^{\widehat{\omega}_j}\big(\qsf,\vec\zeta \ \big):=
\Trr_{\Vcal(\,\widehat{\omega}_j\,)}\,
\qsf^{L_0^{\widehat{\mathfrak{sl}}_{k}}-\frac{k-1}{24} \, \id} \ \vec x^{\ \vec h}
\end{equation} 
the character
of the $j$-th dominant highest weight representation of $\slfrakhat_{k}$
at level one, with weight the $j$-th fundamental weight
$\widehat{\omega}_j$ of type $\widehat{A}_{k-1}$ for $j=0,1, \ldots,
k-1$; here $\vec\zeta:=\sum_{i=1}^{k-1}\, z_i\, H_i$ and $x_i:=\e^{2\pi\ii
  z_i}$ for $i=1,\dots,k-1$, while $\vec h:=(h_1,\dots,h_{k-1})$. Setting $\xi_i=\e^{2\pi \operatorname{i}\,
  (2z_{i}-z_{i-1}-z_{i+1})}$ for $i=1,\dots,k-1$ with $z_0=z_k=0$, by explicit computation
  of the character we get \cite[Section
  5.3]{art:bruzzopedrinisalaszabo2013}
\begin{equation}\label{eq:character}
\chi^{\widehat{\omega}_j}\big(\qsf,\vec\zeta \ \big)=\frac{1}{\eta(\qsf)^{k-1}} \,
\sum_{\vec{u}\in\Ufrak_j}\, \qsf^{\frac{1}{2}\, \vec{u}\cdot
  C^{-1}\vec{u}}\ \vec{\xi}^{\ C^{-1}\vec{u}}\ ,
\end{equation}
where
$\qsf^{\frac1{24}}\, \eta(\qsf)^{-1} = \prod_{n=1}^\infty\, (1-q^n)^{-1} =\Trr_{\Fcal_{\C(\varepsilon_1,\varepsilon_2)}}\,
\qsf^{L_0^\hfrak}$ is the character of the Fock space representation of
the Heisenberg algebra $\hfrak$.
By Equation \eqref{eq:nakajimayoshioka} and the identity
\begin{equation}\label{eq:identity}
\sum_{i=1}^k\, \frac{1}{\varepsilon_1^{(i)}\,
  \varepsilon_2^{(i)}}=\frac{1}{k\, \varepsilon_1\, \varepsilon_2}
\end{equation}
we obtain explicitly
\begin{equation}
\Zcal_{X_k}\big(\varepsilon_1,\varepsilon_2; \qsf, \vec{\xi} \ \big)_j
= \eta(\qsf)^{k-1} \ \chi^{\widehat{\omega}_j}\big(\qsf,\vec\zeta \ \big) \ \exp\Big(\, \frac \qsf{k\,\varepsilon_1\, \varepsilon_2}\, \Big) \ .
\end{equation}

\subsubsection{Gaiotto state}

Following Section \ref{sec:gaiottoC2}, we define the \emph{Gaiotto state} $G_j$ to be the sum, in the completed total equivariant cohomology $\widehat{\W}_j$,  of all fundamental classes
\begin{equation}
G_j:= \sum_{\vec{u}\in\Ufrak_j} \ \sum_{n\geq 0} \,
\big[\Mcal(\vec{u},n,j) \big]_T\ .
\end{equation}
We also define the \emph{weighted Gaiotto state}
\begin{equation}
G_j\big(\qsf,\vec{\xi} \ \big) := \sum_{\vec{u}\in\Ufrak_j} \ \sum_{
  n\geq 0} \, \qsf^{n+\frac{1}{2}\, \vec{u} \cdot C^{-1}\vec{u}}\,
\vec{\xi}^{\ C^{-1}\vec{u}}\ \big[\Mcal(\vec{u},n,j) \big]_T
\end{equation}
in the completion
\begin{equation}
\widehat{\W}_j\big(\qsf,\vec{\xi} \ \big):=\prod_{\vec{u}\in\Ufrak_j}
\ \prod_{n\geq 0}\, \qsf^{n+\frac{1}{2}\, \vec{u} \cdot
  C^{-1}\vec{u}}\ \vec{\xi}^{\ C^{-1}\vec{u}} \ \W_{\vec{u},n,j} \ .
\end{equation}
If we endow $\widehat{\W}_j\big(\qsf,\vec{\xi} \ \big)$ with the scalar product
\begin{multline}
\Big\langle \, \mbox{$\sum\limits _{\vec{u}\in\Ufrak_j}$} \ \mbox{$\sum\limits _{ n\geq0}$}\,
\qsf^{n+\frac{1}{2}\, \vec{u} \cdot C^{-1}\vec{u}}\ \vec{\xi}^{\
  C^{-1}\vec{u}} \ \eta_{\vec{u},n} \,,\,
\mbox{$\sum\limits _{\vec{u}\in\Ufrak_j} $}\ \mbox{$\sum\limits _{
  n\geq0}$}\, \qsf^{n+\frac{1}{2}\, \vec{u} \cdot C^{-1}\vec{u}}\
\vec{\xi}^{\ C^{-1}\vec{u}} \ \nu_{\vec{u},n}
\Big\rangle_{\widehat{\W}_j(\qsf,\vec{\xi}\ )}\\
:= \sum_{\vec{u}\in\Ufrak_j}\ \sum_{ n=0}^\infty\, \qsf^{\frac{1}{2}\,
  \vec{u} \cdot C^{-1}\vec{u}} \ (-\qsf)^n \ \vec{\xi}^{\
  C^{-1}\vec{u}} \ \big\langle \eta_{\vec{u},n}, \nu_{\vec{u},n} \big\rangle_{\W_{\vec{u},n,j}}\ ,
\end{multline}
\normalsize
then it is straightforward to see that the norm of the
weighted Gaiotto state is the instanton partition function for the $\Ncal=2$ $U(1)$ gauge theory on $X_k$:
\begin{equation}
\Zcal_{X_k}\big(\varepsilon_1,\varepsilon_2; \qsf, \vec{\xi} \ \big)_j
= \big\langle G_j\big(\qsf,\vec{\xi} \ \big)\,,\,
G_j\big(\qsf,\vec{\xi} \
\big) \big\rangle_{\widehat{\W}_j(\qsf,\vec{\xi} \ )} 
\ .
\end{equation}
\begin{proposition}
The Gaiotto state is a Whittaker vector for $\glfrakhat_k$ of type
$\chi$, where the algebra homomorphism
$\chi\colon \Ucal(\hfrak^+ \oplus
\hfrak_{\C(\varepsilon_1,\varepsilon_2), \Qfrak}^+) \to \C(\varepsilon_1, \varepsilon_2)$ is defined by
\begin{align}
\chi(\qfrak^i_m)  &=  0  \qquad \mbox{for} \quad m>0 \ , \ i = 1,\ldots, k-1\ ,\\[4pt]
\chi(\pfrak_m) & =  \delta_{m,1}\ 
\sqrt{\big\langle [X_k], [X_k] \big\rangle_{\HH_1}} \qquad \mbox{for}
\quad m>0\ .
\end{align}
\end{proposition}
\proof
We first note that under the isomorphism $\Psi_j$ defined in
\eqref{eq:isomWH} the fundamental class $[\Mcal(\vec{u},n,j)]_T$ is
sent to $[\hilb{n}{X_k}]_T \otimes (\gamma_{\vec{u}}+\omega_j)$. Hence
under the isomorphism \eqref{eq:omega} the Gaiotto state becomes
\begin{equation}
G_j = \bigotimes_{i=1}^k \ \sum_{n\geq0}\,  \big[\hilb{n}{U_i} \big]_T
\ \otimes \ \sum_{\vec{u}\in\Ufrak_j} \, (\gamma_{\vec{u}}+\omega_j) \
\in \ \mbox{$\bigotimes\limits_{i=1}^k$} \, \widehat{\HH}_{U_i} \ \otimes \
\mbox{$\prod\limits_{\vec{u}\in\Ufrak_j} $} \, \C(\varepsilon_1,\varepsilon_2) (\gamma_{\vec{u}}+\omega_j)\ .
\end{equation}
By Proposition \ref{prop:gaiottowhittakerC2}, $\sum_{n\geq0}\,
[\hilb{n}{U_i}]_T$ is the Whittaker vector $G(\eta_i)$ for the
Heisenberg algebra $\hfrak_i$ with $\eta_i = \big(\varepsilon_2^{(i)}
\big)^{-1}$. It follows that $G_j$ is the Whittaker vector
$G_j(\vec{\eta}\, )$ for $\glfrakhat_k$ as in Proposition
\ref{prop:whittakerALE} with $\vec{\eta}=(\eta_1, \ldots, \eta_k)$ of
type $\chi$ where
\begin{equation}
\chi(\qfrak_m^i)  =  \delta_{m,1} \, \big(\eta_{i+1}\,
\beta_{i+1}^{-1} - \eta_i \big)\qquad\mbox{and}\qquad 
\chi(\pfrak_m)  =  \delta_{m,1} \ \sqrt{-k\, \varepsilon_1\,
  \varepsilon_2} \ \sum_{i=1}^k\,
\frac{\eta_i}{\varepsilon_1^{(i)}}\ .
\end{equation}
By computing explicitly the quantities on the right-hand sides of these
equations, one gets the assertion.
\endproof

\subsection{Quiver gauge theories}

As we did in Section \ref{sec:R4quivers}, we will now add matter to
the $\Ncal=2$ gauge theory and consider $\Ncal=2$
superconformal quiver gauge theories on the ALE space $X_k$ with gauge group
$U(1)^{r+1}$ for $r\geq0$; we shall follow \cite{art:bruzzosalaszabo2015}, where superconformal quiver gauge theories on the ALE space $X_k$ are introduced. 

For a quiver $\quiv=(\quiv_0,\quiv_1)$ we
fix vectors of integers
$(\vec
u_\upsilon,n_\upsilon,j_\upsilon)_{\upsilon\in\quiv_0}$ representing
the topological numbers of the moduli spaces $\Mcal(\vec u_\upsilon,
n_\upsilon, j_\upsilon)$ at the
vertices $\quiv_0$ with $\vec u_\upsilon\in\Ufrak_{j_\upsilon}$,
$n_\upsilon\in \N$, and $j_\upsilon\in\{0,1,\dots,k-1\}$. The fundamental (resp. antifundamental) hypermultiplets of masses
$\mu^s_{\upsilon}$, $s=1,\dots,m_\upsilon$ (resp. $\bar\mu^{\bar s}_{\upsilon}$,
$\bar s=1,\dots,\bar m_\upsilon$) at the nodes $\upsilon\in\quiv_0$ correspond to the $T$-equivariant vector
bundles $\Vbf_{\mu^s_{\upsilon}}^{\vec u_\upsilon, n_\upsilon,j_\upsilon}$
(resp. $\bar\Vbf_{\bar\mu^{\bar s}_{\upsilon}}^{\vec u_\upsilon,
  n_\upsilon,j_\upsilon}$) on $\Mcal(\vec u_\upsilon,n_\upsilon,j_\upsilon)$ obtained by
pushforward of $\Ebf_{\mu^s_\upsilon}^{\vec
  u_\upsilon,n_\upsilon,j_\upsilon; \vec u_\upsilon,n_\upsilon,j_\upsilon}$
(resp. $\Ebf_{\bar\mu^{\bar s}_\upsilon}^{\vec
  u_\upsilon,n_\upsilon,j_\upsilon; \vec u_\upsilon,n_\upsilon,j_\upsilon}$) with respect to the
projection of $\Mcal(\vec u_\upsilon,n_\upsilon,j_\upsilon) \times \Mcal(\vec u_\upsilon,n_\upsilon,j_\upsilon)$ to the
second (resp. first) factor. The bifundamental hypermultiplets of
masses $\mu_e$ at the edges
$e\in \quiv_1$ correspond to
the vector bundles $\Ebf_{\mu_e}^{\vec u_{\source(e)},
  n_{\source(e)},j_{\source(e)};\vec
  u_{\tail(e)},n_{\tail(e)},j_{\tail(e)}}$ on $\Mcal(\vec
u_{\source(e)},n_{\source(e)},j_{\source(e)})\times \Mcal(\vec
u_{\tail(e)},n_{\tail(e)}, j_{\tail(e)})$; for vertex loops with $\source(e)=\tail(e)$ the restriction of
$\Ebf_{\mu_e}^{\vec u_{\source(e)},n_{\source(e)},j_{\source(e)};\vec u_{\source(e)},n_{\source(e)},j_{\source(e)}}$ to the diagonal of $\Mcal(\vec
u_{\source(e)},n_{\source(e)},j_{\source(e)}) \times \Mcal(\vec
u_{\source(e)},n_{\source(e)},j_{\source(e)})$ describes an adjoint hypermultiplet
of mass $\mu_e$. The total matter field content of
the $\Ncal=2$ quiver gauge theory on $X_k$ associated to $\quiv$
in the sector labelled by $(\vec
u_\upsilon,n_\upsilon,j_\upsilon)_{\upsilon\in\quiv_0}$
is thus described by the bundle on $\prod_{\upsilon\in\quiv_0}\,\Mcal(\vec u_\upsilon,n_\upsilon,j_\upsilon)$
 given by
\begin{align}
\boldsymbol{M}^{(\vec
  u_\upsilon,n_\upsilon,j_\upsilon)}_{(\mu^s_{\upsilon}),(\bar\mu^{\bar
    s}_\upsilon),
  (\mu_e)}:= & \ \bigoplus_{\upsilon\in\quiv_0}\, p_\upsilon^* \Big(\, \bigoplus_{s=1}^{m_\upsilon}\,
\Vbf_{\mu^s_{\upsilon}}^{\vec u_\upsilon, n_\upsilon,j_\upsilon} \
\oplus \ \bigoplus_{\bar s=1}^{\bar m_\upsilon}\,
\bar\Vbf_{\bar\mu^{\bar s}_{\upsilon}}^{\vec
  u_\upsilon,n_\upsilon,j_\upsilon}\, \Big) \\ & \ \qquad \ \oplus \ \bigoplus_{e\in\quiv_1}\,
p_e^*\Ebf_{\mu_e}^{\vec
  u_{\source(e)},n_{\source(e)},j_{\source(e)};\vec u_{\tail(e)},n_{\tail(e)},j_{\tail(e)}} \ ,
\end{align}
where $p_\upsilon$ is the projection of $\prod_{\upsilon\in\quiv_0}\,
\Mcal(\vec u_\upsilon,n_\upsilon,j_\upsilon)$ to the $\upsilon$-th
factor while $p_e$ is the projection
  to the product $\Mcal(\vec
u_{\source(e)},n_{\source(e)},j_{\source(e)})\times \Mcal(\vec
u_{\tail(e)},n_{\tail(e)}, j_{\tail(e)})$. 

The degree of the Euler class of the hypermultiplet bundle $\boldsymbol{M}^{(\vec
  u_\upsilon,n_\upsilon,j_\upsilon)}_{(\mu^s_{\upsilon}),(\bar\mu^{\bar
    s}_\upsilon),
  (\mu_e)}$ is given by
\begin{multline}
{\rm deg}\, \eu\big(\boldsymbol{M}^{(\vec
  u_\upsilon,n_\upsilon,j_\upsilon)}_{(\mu^s_{\upsilon}),(\bar\mu^{\bar
    s}_\upsilon),
  (\mu_e)}\, \big):= \sum_{\upsilon\in \quiv_0}\, \dim \Mcal(\vec
u_\upsilon,n_\upsilon,j_\upsilon)- \rk\big(\boldsymbol{M}^{(\vec
  u_\upsilon,n_\upsilon,j_\upsilon)}_{(\mu^s_{\upsilon}),(\bar\mu^{\bar
    s}_\upsilon),
  (\mu_e)}\,\big)\\[4pt] \shoveleft{=\sum_{\upsilon\in \quiv_0}\,2\, n_\upsilon -\, \sum_{\upsilon\in \quiv_0}\, \big(m_\upsilon+\bar m_\upsilon\big)\, \Big(n_\upsilon+\mbox{$\frac12$}\, \vec v_\upsilon\cdot C \vec
v_\upsilon-\mbox{$\frac1{2k}$}\,
j_\upsilon \, (k- j_\upsilon) \Big)} \\
-\sum_{e\in \quiv_1}\, \Big(\, n_{\tail(e)}+n_{\source(e)}
+\, \mbox{$\frac12$} \, \vec{v}_{\source(e)}\cdot C\vec{v}_{\source(e)}+\mbox{$\frac12$}\, \vec{v}_{\tail(e)}\cdot C\vec{v}_{\tail(e)} -\vec{v}_{\source(e)}\cdot C\vec{v}_{\tail(e)}-\mbox{$\frac1{2k}$}\,
j_{\tail(e)\source(e)} \, (k- j_{\tail(e)\source(e)})\, \Big)\ ,
\end{multline}
where $\vec{v}_\upsilon:=C^{-1}\vec{u}_\upsilon$. Using \eqref{eq:confconstr} the degree becomes
\begin{equation}
{\rm deg}\, \eu\big(\boldsymbol{M}^{(\vec
  u_\upsilon,n_\upsilon,j_\upsilon)}_{(\mu^s_{\upsilon}),(\bar\mu^{\bar
    s}_\upsilon),
  (\mu_e)}\,\big) = 
\sum_{\upsilon\in \quiv_0}\,d^{X_k}_\upsilon(\vec{\boldsymbol v}_\upsilon) \ ,
\end{equation}
where we defined
\begin{multline}
d^{X_k}_\upsilon(\vec{\boldsymbol v}_\upsilon) := \mbox{$\frac1{2k}$}\,
j_\upsilon \, (k- j_\upsilon)\,\big(2-\#\{e\in\quiv_1\, | \, \source(e)=\upsilon\}-
\#\{e\in\quiv_1\, | \, \tail(e)=\upsilon\}\big) \\
+\mbox{$\frac1{4k}$}\,\Big( \sum_{\stackrel{\scriptstyle e\in\quiv_1}{\scriptstyle \source(e)=\upsilon}} \,j_{\tail(e)\upsilon} \, (k- j_{\tail(e)\upsilon})
+\sum_{\stackrel{\scriptstyle e\in\quiv_1}{\scriptstyle \tail(e)=\upsilon}} \,j_{\upsilon\source(e)} \, (k- j_{\upsilon\source(e)}) \Big)  \\
- \vec{v}_\upsilon\cdot C \vec{v}_\upsilon+\mbox{$\frac12$}\, \sum_{\stackrel{\scriptstyle e\in\quiv_1}{\scriptstyle \source(e)=\upsilon}}\, \vec{v}_\upsilon\cdot C \vec{v}_{\tail(e)}+\mbox{$\frac12$}\, \sum_{\stackrel{\scriptstyle e\in\quiv_1}{\scriptstyle\tail(e)=\upsilon}}\, \vec{v}_\upsilon\cdot C \vec{v}_{\source(e)}
\end{multline}
for each vertex $\upsilon\in\quiv_0$; here $\vec{\boldsymbol v}_\upsilon:=\big(\vec v_\upsilon,(\vec v_{\tail(e)})_{e\in\quiv_1\,:\, \source(e)=\upsilon},(\vec v_{\source(e)})_{e\in\quiv_1\,:\, \tail(e)=\upsilon}\big)$. By analogy with the case of gauge theories on $\R^4$ (see Section \ref{sec:R4quivers}), we say that the $\Ncal=2$ quiver gauge theory on $X_k$ is conformal if $d^{X_k}_\upsilon(\vec{\boldsymbol v}_\upsilon)=0$ for all $\upsilon\in\quiv_0$; this is formally the requirement of vanishing beta-function for the running of the $\upsilon$-th gauge coupling constant. For any vertex $\upsilon\in \quiv_0$, define the set of conformal fractional instanton charges by
\begin{equation}\label{eq:confchargeset}
\Ufrak_{j_\upsilon}^{\rm conf}:=\big\{\vec u_\upsilon\in\Ufrak_{j_\upsilon} \
\big\vert \ d^{X_k}_\upsilon(\vec{\boldsymbol v}_\upsilon)=0 \big\}\ .
\end{equation}
The conformal constraint is always trivially satisfied by any $\vec u_\upsilon$ for the $\widehat{A}_0$-theory, while for the $A_0$-theory the set of conformal fractional instantons charges reduces to
\begin{equation}
\Ufrak_{j}^{\rm conf}:= \big\{\vec u\in\Ufrak_{j} \
\big\vert \ \vec u\cdot C^{-1}\vec
u= \mbox{$\frac1k$}\, j\,(k- j) \big\} \ . 
\end{equation}
Note that in this case, this is a restriction on the
conformal dimension $\Delta_{\vec u}=\frac1{2} \,\langle \omega_{j},\omega_{j}\rangle_{\Qfrak\otimes_\Z \Q}$ of the highest weight representation $\W_{\vec u,j}$ of the Virasoro algebra.

Introduce topological couplings $\qsf_\upsilon\in\C^*$ with
$|\qsf_\upsilon|<1$ and $\vec\xi_v=\big((\xi_v)_1,\dots,
(\xi_v)_{k-1}\big) \in (\C^*)^{k-1}$ with $|(\xi_v)_i|<1$ at each vertex
$\upsilon\in\quiv_0$. With notation as in Section \ref{sec:R4quivers},
the quiver gauge theory partition function on $X_k$ is then defined by
the generating function
\begin{multline}
\Zcal_{X_k}^{\quiv}
\big(\varepsilon_1,\varepsilon_2,\mubf;\qbf, \vec\xibf \ \big)_{\jbf} :=
\sum_{(\vec u_\upsilon\in\Ufrak_{j_\upsilon}^{\rm conf})} \,
\vec\xibf{}^{\ C^{-1}\vec\ubf} \ \sum_{(n_\upsilon)\in\N^{\quiv_0}}\,
\qbf^{\boldsymbol{n} + \frac12\, \vec\ubf\cdot C^{-1}\vec\ubf} \\ \shoveright{
\times \ 
  \int_{\prod\limits_{\upsilon\in\quiv_0}\, \Mcal(\vec
    u_\upsilon,n_\upsilon,j_\upsilon)} \, \eu_{T\times
  T_{\mubf}}\big(\boldsymbol{M}^{(\vec
  u_\upsilon,n_\upsilon,j_\upsilon)}_{(\mu^s_{\upsilon}),(\bar\mu^{\bar
    s}_\upsilon),
  (\mu_e)} \big) } \\[4pt] \shoveleft{ = \sum_{(\vec u_\upsilon\in\Ufrak_{j_\upsilon}^{\rm conf})} \
\sum_{(n_\upsilon)\in\N^{\quiv_0}}\,
\qbf^{\boldsymbol{n} + \frac12\, \vec\ubf\cdot C^{-1}\vec\ubf} \
\vec\xibf{}^{\ C^{-1}\vec\ubf} } \\ \times \
  \int_{\prod\limits_{\upsilon\in\quiv_0}\, \Mcal(\vec
    u_\upsilon,n_\upsilon,j_\upsilon)} \ 
  \prod_{\upsilon\in\quiv_0}\, p_\upsilon^*\Big(\, \prod_{s=1}^{m_\upsilon}\,
  \eu_T\big(\Vbf_{\mu^s_{\upsilon}}^{\vec u_\upsilon, n_\upsilon,j_\upsilon}\big) \
  \prod_{\bar s=1}^{\bar m_\upsilon}\,
  \eu_T\big(\bar\Vbf_{\bar\mu^{\bar s}_{\upsilon}}^{\vec u_\upsilon,
  n_\upsilon,j_\upsilon}\big)\, \Big) \\
  \times \ \prod_{e\in\quiv_1}\,
  p_e^*\, \eu_T\big(\Ebf_{\mu_e}^{\vec u_{\source(e)},
  n_{\source(e)},j_{\source(e)};\vec
  u_{\tail(e)},n_{\tail(e)},j_{\tail(e)}}
  \big) \ ,
\end{multline}
where $\qbf^{\nbf + \frac12\, \vec\ubf\cdot C^{-1}\vec\ubf}:=
\prod_{\upsilon\in\quiv_0}\, \qsf_\upsilon^{n_\upsilon+\frac12\, \vec
  u_\upsilon\cdot C^{-1}\vec u_\upsilon}$ and $\vec\xibf{}^{\
  C^{-1}\vec\ubf} := \prod_{\upsilon\in\quiv_0}\, \vec\xi_\upsilon^{\
  C^{-1}\vec u_\upsilon}$. By applying 
the localization theorem, and using Equations \eqref{eq:ZcalC2quiver}
and \eqref{eq:vertexoperator-fixedbasis}, we obtain a factorization in
terms of $\Ncal=2$ quiver gauge theory partition functions on $\R^4$
weighted by edge contributions. For the fundamental and
antifundamental matter fields, the relevant edge contributions
$\ell_{\vec v_\upsilon}^{(n)}$ are the equivariant Euler classes of
$H^1(\Xscr_k,\Rcal^{\vec{u}_\upsilon}\otimes
\Ocal_{\Xscr_k}(-\Dscr_\infty))$ for $\upsilon\in\quiv_0$; by
\cite[Section 5]{art:bruzzopedrinisalaszabo2013} this vector space is
zero if and only if $\vec u_\upsilon\in\Ufrak_{j_\upsilon}^{\rm conf}$ and the
corresponding edge contribution is equal to one. Thus only the arrows
of the quiver yield edge contributions and the partition function is
given by
\begin{multline}
\Zcal_{X_k}^{\quiv}\big(\varepsilon_1,\varepsilon_2,\mubf;\qbf, \vec\xibf
\ \big)_{\jbf} = \prod_{\upsilon\in\quiv_0}\, \qsf_\upsilon^{\frac12\,
  \langle \omega_{j_\upsilon},\omega_{j_\upsilon}\rangle_{\Qfrak\otimes_\Z \Q}} \ \sum_{(\vec u_\upsilon\in\Ufrak_{j_\upsilon}^{\rm conf})} \,
\vec\xibf{}^{\ C^{-1}\vec\ubf} \ \prod_{n=1}^{k-1} \ \prod_{e\in\quiv_1}\, \ell^{(n)}_{\vec
  v_e}\big(\varepsilon_1^{(n)},\varepsilon_2^{(n)},\mu_e \big) \\
\times \ \prod_{i=1}^k\,
\Zcal_{\C^2}^{\quiv}\big(\varepsilon_1^{(i)},\varepsilon_2^{(i)},\mubf^{(i)}
;\qbf \big) \ , 
\end{multline}
where $\vec v_e:= \vec v_{\tail(e)}-\vec v_{\source(e)}$; the shifted
masses are
\begin{equation}
\big(\mu_v^s\big)^{(i)}:= \mu_v^s-(\vec v_v)_i\, \varepsilon_1^{(i)} -
(\vec v_v)_{i-1} \, \varepsilon_2^{(i)} 
\end{equation}
for $v\in\quiv_0$, $s=1,\dots,m_v$ and $i=1,\dots,k$, and similarly
for $\big(\bar\mu_v^{\bar s}\big)^{(i)}$, whereas
\begin{equation}
\mu_e^{(i)}:= \mu_e-(\vec v_e)_i\, \varepsilon_1^{(i)} -
(\vec v_e)_{i-1} \, \varepsilon_2^{(i)}
\end{equation}
for $e\in\quiv_1$ and $i=1,\dots,k$.

In the remainder of this section we consider in detail each of the
admissible quivers in turn.

\subsection{$\widehat{A}_{r}$ theories}

With the conventions of Section
\ref{sec:hatArtheories}, the instanton partition function for the
$\Ncal=2$ $U(1)^{r+1}$ quiver gauge theory of type $\widehat{A}_{r}$ on the ALE
space $X_k$ reads
as
\begin{multline}
\Zcal_{X_k}^{\widehat{A}_r}\big(\varepsilon_1,\varepsilon_2,\mubf;\qbf, \vec\xibf
\ \big)_{\jbf} = \prod_{\upsilon=0}^r \, \qsf_\upsilon^{\frac12\,
  \langle \omega_{j_\upsilon},\omega_{j_\upsilon}\rangle_{\Qfrak\otimes_\Z \Q}} \ \sum_{(\vec u_\upsilon\in\Ufrak_{j_\upsilon}^{\rm conf})} \,
\vec\xibf{}^{\
  C^{-1}\vec\ubf} \ \prod_{n=1}^{k-1} \ \prod_{\upsilon=0}^r \, \ell^{(n)}_{\vec
  v_{\upsilon,\upsilon+1}} \big(\varepsilon_1^{(n)},\varepsilon_2^{(n)},\mu_\upsilon \big) \\
\times \ \prod_{i=1}^k\,
\Zcal_{\C^2}^{\widehat{A}_r} \big(\varepsilon_1^{(i)},\varepsilon_2^{(i)},\mubf^{(i)}
;\qbf \big) \ , \label{eq:ZXkhatAr} 
\end{multline}
where $\vec\xibf{}^{\
  C^{-1}\vec\ubf} := \prod_{\upsilon=0}^r\, \vec\xi_\upsilon^{\
  C^{-1}\vec u_\upsilon}$ with $\vec u_{r+1}:= \vec u_0$, while
$\mu_\upsilon^{(i)}:= \mu_\upsilon-(\vec v_{\upsilon,\upsilon+1})_i\, \varepsilon_1^{(i)} -
(\vec v_{\upsilon,\upsilon+1})_{i-1} \, \varepsilon_2^{(i)}$ with
$\vec v_{\upsilon,\upsilon+1}:=\vec
v_{\upsilon+1}-\vec v_{\upsilon}$.

\subsubsection{Conformal blocks}

We will relate the partition function \eqref{eq:ZXkhatAr} to the trace of vertex
operators $\V_{\mu_\upsilon}^{j_\upsilon, j_{\upsilon+1}}(\vec x_\upsilon,z_\upsilon)$ from Section \ref{sec:chiralvertex},
analogously to what we did in
Section \ref{sec:hatArtheories}, and hence interpret it as a torus
$(r+1)$-point conformal block. For this, we fix vertices
$\upsilon,\upsilon'\in\{0,1,\dots, r\}$ and introduce the \emph{conformal
  restriction operators} $\delta^{\rm
  conf}_{\upsilon,{\upsilon'}}\colon\W\to \W$ which are defined by
their matrix elements in the fixed point basis of the vector space $\W$ by
\begin{equation}
\big\langle \delta^{\rm
  conf}_{\upsilon,{\upsilon'}}\triangleright [\vec Y,\vec u\, ]\,,\, [\vec
Y',\vec u\,'\,] \big\rangle_{\W}:= 
\left\{
\begin{array}{ll}
1 & \mbox{if} \quad \vec u\in\Ufrak_{j_\upsilon}^{\rm conf} \ ,\ \vec u\,'\in
\Ufrak_{j_{\upsilon'}}^{\rm conf} \ , \\[4pt]
0 & \mbox{otherwise} \ .
\end{array} \right. \label{eq:confrestrop}
\end{equation}
Suitable insertions of this operator restrict the first Chern classes $\vec
u_\upsilon\in\Ufrak_{j_\upsilon}$ in the way required by the
superconformal constraints of the quiver gauge theory and the
constrained conformal dimensions of the associated Virasoro algebras
at the nodes of the $\widehat{A}_r$-type quivers. Using Proposition
\ref{prop:L0eigen} and Proposition \ref{prop:Wkweight}, by performing
analogous manipulations to those used in the proof of Proposition
\ref{prop:hatArtrace} we arrive at the following result.

\begin{proposition}\label{prop:ZcalXkArTr}
The partition function of the $\widehat{A}_r$-theory on $X_k$ is given
by
\begin{equation}
\Zcal_{X_k}^{\widehat{A}_r}\big(\varepsilon_1,\varepsilon_2,\mubf;\qbf, \vec\xibf
\ \big)_{\jbf} = \Trr_{\W_{j_0}} \, \qsf^{L_0}\ \vec\xi{}^{\ C^{-1}\vec h} \ \prod_{\upsilon=0}^r\,
\V_{\mu_\upsilon}^{j_\upsilon,j_{\upsilon+1}}(\vec x_\upsilon,z_\upsilon)\ \delta^{\rm
  conf}_{\upsilon,{\upsilon+1}}
\end{equation}
independently of $z_0\in\C^*$ and $\vec x_0\in(\C^*)^{k-1}$, where
$\qsf:= \qsf_0\, \qsf_1\cdots \qsf_r$, $( \, \vec\xi\ )_i:=
(\vec\xi_0)_i\, (\vec\xi_1)_i\, \cdots (\vec\xi_r)_i$, $z_\upsilon:=
z_0\, \qsf_1\cdots \qsf_\upsilon$, and $(\vec x_\upsilon)_i:= (\vec
x_0)_i\, (\vec\xi_1)_i\cdots (\vec\xi_\upsilon)_i$ for
$\upsilon=1,\dots,r$ and $i=1,\dots,k-1$.
\end{proposition}

By combining Theorem \ref{thm:virasoroprimary} and Proposition \ref{prop:ZcalXkArTr}, it follows that the quiver gauge theory partition function completely factorizes under the isomorphism of Proposition \ref{prop:representation} into partition functions associated to the affine algebras $\hfrak$ and $\slfrakhat_k$.
\begin{corollary}\label{cor:trace}
Let $\V_\mu(\vec{v}_{21},\vec x,
z)$ be the vertex operator in
$\Hom(\W_{j_1},
\W_{j_2})[[
z^{\pm\,1},x_1^{\pm\,1},\dots,x_{k-1}^{\pm\,1}]]$ given by
\begin{equation}\label{eq:Vmudef}
\V_\mu(\vec{v}_{21},\vec x,
z):=  \, z^{\Delta_{\vec u_2}-\Delta_{\vec u_1}} {\bar \V}_\mu(\vec{v}_{21}, \vec x, z)\ \exp\big(\log z\ \cfrak-\gamma_{21} \big)\, \exp\big(\gamma_{21}\,\big)\ .
\end{equation}
Then
\begin{multline}
\Zcal_{X_k}^{\widehat{A}_r}\big(\varepsilon_1,\varepsilon_2,\mubf;\qbf, \vec\xibf
\ \big)_{\jbf}
= \Zcal^{\widehat{A}_{r}}_{\C^2}(\varepsilon_1,\varepsilon_2,\mubf;\qsf)^{\frac{1}{k}}\ \qsf^{\frac{1}{24}\,(1-\frac{1}{k})}\,\eta(\qsf)^{\frac{1}{k}-1} \\
\times\ \Trr_{\Vcal(\,\widehat{\omega}_{j_0}\,)} \,
\qsf^{L_0^{\slfrakhat_k}}\ \vec\xi{}^{\ C^{-1}\vec h}\ \prod_{\upsilon=0}^r\ 
\sum_{(\vec{u}_\upsilon\in\Ufrak^{\rm conf}_{j_\upsilon})} \,
\V_{\mu_\upsilon}(\vec{v}_{\upsilon,\upsilon+1},\vec x_\upsilon, z_\upsilon) \, \big\vert
_{\W_{\vec{u}_\upsilon,j_\upsilon}} \ .
\end{multline}
\end{corollary}
\proof
By Theorem \ref{thm:virasoroprimary} and Proposition \ref{prop:representation} we get
\begin{multline}
\Trr_{\W_{j_0}} \, \qsf^{L_0}\ \vec\xi{}^{\ C^{-1}\vec h} \ \prod_{\upsilon=0}^r\,
\V_{\mu_\upsilon}^{j_\upsilon,j_{\upsilon+1}}(\vec x_\upsilon,z_\upsilon)\ \delta^{\rm
  conf}_{\upsilon,{\upsilon+1}}
  =\Trr_{\Fcal_{\C(\varepsilon_1,\varepsilon_2)}} \, \qsf^{L_0^{\hfrak}} \ \prod_{\upsilon=0}^r\, 
\V_{-\frac{\mu_\upsilon}{\sqrt{-k\,\varepsilon_1\, \varepsilon_2}}, \frac{\mu_\upsilon+\varepsilon_1+\varepsilon_2}{\sqrt{-k\, \varepsilon_1\, \varepsilon_2}}}(z_\upsilon)\\
\times\ \Trr_{\Vcal(\,\widehat{\omega}_{j_0}\,)} \, \qsf^{L_0^{\slfrakhat_k}}\ \vec\xi{}^{\ C^{-1}\vec h} \ \prod_{\upsilon=0}^r \
\sum_{\vec{u}_1\in\Ufrak_{j_\upsilon},\vec{u}_2\in\Ufrak_{j_{\upsilon+1}}} \,
\V_\mu(\vec{v}_{21},z_\upsilon)\, \big\vert
_{\W_{\vec{u}_1,j_\upsilon}}\ \delta^{\rm
  conf}_{\upsilon,{\upsilon+1}} \ .
\end{multline}
Then by using the same arguments as in the proof of \cite[Corollary 1]{art:carlssonokounkov2012} one gets
\begin{equation}
\Trr_{\Fcal_{\C(\varepsilon_1,\varepsilon_2)}} \, \qsf^{L_0^{\hfrak}} \ \prod_{\upsilon=0}^r\,
\V_{-\frac{\mu_\upsilon}{\sqrt{-k\, \varepsilon_1\, \varepsilon_2}}, \frac{\mu_\upsilon+\varepsilon_1+\varepsilon_2}{\sqrt{-k\, \varepsilon_1\, \varepsilon_2}}}(z_\upsilon)=\prod_{\upsilon=0}^r\,\big(\qsf_\upsilon^{-\frac{1}{24}}\, \eta(\qsf_\upsilon)\big)^{-\frac{\mu_\upsilon\,
    (\mu_\upsilon+\varepsilon_1+\varepsilon_2)}{k\, \varepsilon_1\, \varepsilon_2}} \ \qsf^{\frac{1}{24}}\,\eta(\qsf)^{-1} \ .
\end{equation}
The result now follows from Proposition \ref{prop:tracededekind}.
\endproof

\subsubsection{$\widehat{A}_{0}$ theory}

For the $\Ncal=2^*$ gauge theory on $X_k$, similar arguments to those
of Section \ref{sec:N=2Xk} show that the edge contributions are also
equal to
one in this case. In this instance the gauge theory is automatically
conformal without further restriction of the first Chern classes $\vec
u\in\Ufrak_j$. Then the instanton partition function for $U(1)$
gauge theory on the ALE space $X_k$ with a single adjoint
hypermultiplet of mass $\mu$ can be written in a factorized form in
terms of the Nekrasov partition function for $\Ncal=2^\ast$ gauge
theory on $\R^4$ given by
\begin{equation}
\Zcal_{X_k}^{\widehat{A}_0} \big(\varepsilon_1,\varepsilon_2,\mu;
\qsf, \vec{\xi} \ \big)_j =\eta(\qsf)^{k-1} \ \chi^{\widehat{\omega}_j}\big(\qsf,\vec\zeta \ \big) \
\prod_{i=1}^k \,
\Zcal_{\C^2}^{\widehat{A}_0} \big(\varepsilon_1^{(i)},\varepsilon_2^{(i)}, \mu ;
\qsf \big)\ .
\end{equation}
In this case Proposition
\ref{prop:ZcalXkArTr} may be stated in a factorized form under the decomposition of Theorem \ref{thm:virasoroprimary} in terms of characters of
$\hfrak \subset
\glfrakhat_k$ and $\glfrakhat_k$ as
\begin{equation}
\Zcal_{X_k}^{\widehat{A}_0} \big(\varepsilon_1,\varepsilon_2,\mu;
\qsf, \vec{\xi} \ \big)_j = \eta(\qsf)^{k-1} \ \chi^{\widehat{\omega}_j}\big(\qsf,\vec\zeta \ \big) \
\Trr_{\HH} \, \qsf^{L_0^{\hfrak}} \ \V_{-\frac{\mu}{\sqrt{-k\, \varepsilon_1\, \varepsilon
_2}},\frac{\mu+\varepsilon_1+\varepsilon_2}{\sqrt{-k\, \varepsilon_1\,
\varepsilon_2}}}(1) \ .
\end{equation}
By using the
identities \eqref{eq:Ahat0dedekind} and \eqref{eq:identity}, we obtain explicitly
\begin{equation}
\Zcal_{X_k}^{\widehat{A}_0} \big(\varepsilon_1,\varepsilon_2, \mu;
\qsf, \vec{\xi} \ \big)_j \ = \qsf^{\frac{k}{24}}\ \eta(\qsf)^{-1} \ \chi^{\widehat{\omega}_j}\big(\qsf,\vec\zeta \ \big) \
\big(\qsf^{-\frac{1}{24}}\,
\eta(\qsf)\big)^{-\frac{\mu\, (\mu+\varepsilon_1+\varepsilon_2)}{k\,
    \varepsilon_1\, \varepsilon_2}}\ .
\end{equation}
\begin{remark}
Note that $\HH$ is not the Fock space of $\hfrak$, as we have
\begin{equation}
\Trr_{\HH} \, \qsf^{L_0^{\hfrak}} \ \V_{-\frac{\mu}{\sqrt{-k\, \varepsilon_1\, \varepsilon
_2}},\frac{\mu+\varepsilon_1+\varepsilon_2}{\sqrt{-k\, \varepsilon_1\,
\varepsilon_2}}}(1) =\big(\qsf^{\frac{1}{24}}\,
\eta(\qsf)^{-1} \big)^{k-1}\ \Trr_{\Fcal_{\C(\varepsilon_1,\varepsilon_2)}} \, \qsf^{L_0^{\hfrak}} \ \V_{-\frac{\mu}{\sqrt{-k\, \varepsilon_1\, \varepsilon
_2}},\frac{\mu+\varepsilon_1+\varepsilon_2}{\sqrt{-k\, \varepsilon_1\,
\varepsilon_2}}}(1)\ .
\end{equation}
\end{remark}

\subsection{$A_r$ theories}

With the conventions of Section
\ref{sec:Artheories}, the instanton partition function for the
$\Ncal=2$ $U(1)^{r+1}$ quiver gauge theory of type ${A}_{r}$ on the ALE
space $X_k$ reads
as
\begin{multline}
\Zcal_{X_k}^{{A}_r}\big(\varepsilon_1,\varepsilon_2,\mubf;\qbf, \vec\xibf
\ \big)_{\jbf} = \prod_{\upsilon=0}^r \, \qsf_\upsilon^{\frac12\,
  \langle \omega_{j_\upsilon},\omega_{j_\upsilon}\rangle_{\Qfrak\otimes_\Z \Q}} \ \sum_{(\vec u_\upsilon\in\Ufrak_{j_\upsilon}^{\rm conf})} \,
\vec\xibf{}^{\
  C^{-1}\vec\ubf} \ \prod_{n=1}^{k-1} \ \prod_{\upsilon=0}^{r-1} \, \ell^{(n)}_{\vec
  v_{\upsilon,\upsilon+1}} \big(\varepsilon_1^{(n)},\varepsilon_2^{(n)},\mu_{\upsilon+1} \big) \\
\times \ \prod_{i=1}^k\,
\Zcal_{\C^2}^{{A}_r} \big(\varepsilon_1^{(i)},\varepsilon_2^{(i)},\mubf^{(i)}
;\qbf \big) \ . \label{eq:ZXkAr}
\end{multline}

\subsubsection{Conformal blocks}

By performing analogous manipulations to those used in the proof of
Proposition \ref{prop:fundamental-matter}, we can express the
partition function \eqref{eq:ZXkAr} as a particular matrix element of vertex operators and hence interpret it as an $(r+4)$-point conformal block on the sphere. For this, let
\begin{equation}
|0\rangle_{\rm conf}:= \prod_{\upsilon=0}^r\,\delta^{\rm
  conf}_{0,\upsilon}\triangleright [\emptyset, \vec{0}\,] \ .
\end{equation}
\begin{proposition}\label{prop:ZcalXkArcorr}
The partition function of the $A_r$-theory on $X_k$ is given
by
\begin{multline}
\Zcal_{X_k}^{{A}_r}\big(\varepsilon_1,\varepsilon_2,\mubf;\qbf, \vec\xibf
\ \big)_{\jbf} \\
= \Big\langle |0\rangle_{\rm conf} \, , \, V_{\mu_0}(\vec x_0,z_0) \, \Big(\, \prod_{\upsilon=1}^{r}\,
\V_{\mu_{\upsilon}}^{j_{\upsilon-1}, j_{\upsilon}}(\vec x_\upsilon,z_\upsilon)\ \delta^{\rm
  conf}_{{\upsilon-1},{\upsilon}} \, \Big)\, \V_{\mu_{r+1}}(\vec x_{r+1},z_{r+1}) |0\rangle_{\rm conf} \Big\rangle_{\W}
\end{multline}
independently of $z_0\in\C^*$ and $\vec x_0\in(\C^*)^{k-1}$, where
$z_\upsilon:=
z_0\, \qsf_0\, \qsf_1\cdots \qsf_\upsilon$ and $(\vec x_\upsilon)_i:= (\vec
x_0)_i\, (\vec\xi_0)_i\, (\vec\xi_1)_i$ $\cdots (\vec\xi_{\upsilon-1})_i$ for
$\upsilon=1,\dots,r+1$ and $i=1,\dots,k-1$.
\end{proposition}

Again, combining Theorem \ref{thm:virasoroprimary} and Proposition \ref{prop:ZcalXkArcorr} yields a completely factorized form for the quiver gauge theory partition function under the isomorphism of Proposition \ref{prop:representation}. In the following we denote $\Vcal:=\bigoplus_{j=0}^{k-1}\, \Vcal(\,\widehat{\omega}_j\,)$.
\begin{corollary}
Let $\V_\mu(\vec{v}_{21},\vec x,
z)$ be the vertex operator in
$\Hom(\W_{j_1},
\W_{j_2})[[
z^{\pm\,1},x_1^{\pm\,1},\dots,x_{k-1}^{\pm\,1}]]$ given by \eqref{eq:Vmudef}. Then
\begin{multline}
\Zcal_{X_k}^{{A}_r}\big(\varepsilon_1,\varepsilon_2,\mubf;\qbf, \vec\xibf
\ \big)_{\jbf} 
= \Zcal_{\C^2}^{{A}_r}(\varepsilon_1,\varepsilon_2,\mubf;\qbf)^{\frac1k} \\ \shoveleft{ \times \
\Big\langle |0\rangle_{\rm conf} \, , \, \Big(\, \sum_{j_0,j_0'=0}^{k-1} \
\sum_{\vec{u}_0\in\Ufrak_{j_0},\vec{u}_0'\in\Ufrak_{j_0'}} \
\V_{\mu_0}(\vec{v}_{0',0},\vec x_0, z_0)\, \big\vert
_{\W_{\vec{u}_0,j_0}} \, \Big) } \\ \times \
\prod_{\upsilon=1}^{r} \ \sum_{(\vec u_\upsilon\in \Ufrak^{\rm conf}_{j_\upsilon})}\, 
\V_{\mu_{\upsilon}}(\vec v_{\upsilon-1,\upsilon},\vec x_\upsilon,z_\upsilon)\, \big|_{\W_{\vec u_\upsilon,j_\upsilon}} \\ \times \ \Big(\, 
\sum_{j_{r+1},j_{r+1}'=0}^{k-1} \
\sum_{\vec{u}_{r+1}\in\Ufrak_{j_{r+1}},\vec{u}_{r+1}'\in\Ufrak_{j_{r+1}'}} \
\V_{\mu_{r+1}}(\vec{v}_{r+1',r+1},\vec x_{r+1}, z_{r+1})\, \big\vert
_{\W_{\vec{u}_{r+1},j_{r+1}}}\, \Big)
|0\rangle_{\rm conf} \Big\rangle_{\Vcal} \ .
\end{multline}
\end{corollary}
\proof
The proof follows that of Corollary \ref{cor:trace}, and by repeating the proof of Proposition \ref{prop:ZC2Arexpl} to compute
\begin{equation}
\Big\langle |0\rangle\,,\, \prod_{\upsilon=0}^{r+1}\, \V_{-\frac{\mu_\upsilon}{\sqrt{-k\, \varepsilon_1\, \varepsilon_2}}, \frac{\mu_\upsilon+\varepsilon_1+\varepsilon_2}{\sqrt{-k\, \varepsilon_1\, \varepsilon_2}}}(z_\upsilon) |0\rangle\Big\rangle_{\Fcal_{\C(\varepsilon_1,\varepsilon_2)}} = \prod_{0\leq \upsilon<\upsilon'\leq r+1}\,
\big(1-\qsf_{\upsilon+1}\cdots
\qsf_{\upsilon'}\big)^{-\frac{\mu_{\upsilon'}\,
    (\mu_\upsilon+\varepsilon_1+\varepsilon_2)}{k\, \varepsilon_1\,
    \varepsilon_2}}
\ .
\end{equation}
\endproof

\subsubsection{$A_0$ theory}

For the $\Ncal=2$ superconformal gauge theory on $X_k$ with two
fundamental hypermultiplets of masses $\mu_0,\mu_1$, the set
$\Ufrak_j^{\rm conf}$ coincides with the rank one limit of the more
general conformal charge sets obtained in \cite[Section
5.4]{art:bruzzopedrinisalaszabo2013}. Analogously to Equation
\eqref{eq:characterdef}, let us define the restricted
$\slfrakhat_k$ characters
\begin{align}
\chi_{\mathrm{conf}}^{\widehat{\omega}_j}\big(\qsf,\vec\zeta \ \big):=& \ \Trr_{\Vcal(\,\widehat{\omega}_j\,)}\,
\qsf^{L_0^{\widehat{\mathfrak{sl}}_{k}}-\frac{k-1}{24}\, \id}\ \vec x\,^{
  \vec h} \ \delta^{\rm conf}_{j,j} \\[4pt] =& \ \frac1{\eta(\qsf)^{k-1}} \
\sum_{\vec{u}\in\Ufrak_j^{\mathrm{conf}}}\,
\qsf^{\frac{1}{2}\, \vec{u}\cdot C^{-1}\vec{u}} \ \vec{\xi}^{\
  C^{-1}\vec{u}}= \frac{\qsf^{\frac{1}{2}\, \langle
    \omega_j,\omega_j\rangle_{\Qfrak\otimes_\Z \Q}}}{\eta(\qsf)^{k-1}} \
\sum_{\vec{u}\in\Ufrak_j^{\mathrm{conf}}}\, \vec{\xi}^{\ C^{-1}\vec{u}}\ .
\end{align}
Then the instanton partition function is given by the factorization
\begin{equation}
\Zcal_{X_k}^{A_0}\big(\varepsilon_1,\varepsilon_2, \mu_0,\mu_1; \qsf,
\vec{\xi} \ \big)_j
=\eta(\qsf)^{k-1} \ \chi_{\mathrm{conf}}^{\widehat{\omega}_j}\big(\qsf,\vec\zeta \ \big) \
\prod_{i=1}^k \,
\Zcal_{\C^2}^{A_0}\big(\varepsilon_1^{(i)},\varepsilon_2^{(i)},
\mu_0,\mu_1 ; \qsf \big) \ .
\end{equation}
In this instance Proposition \ref{prop:ZcalXkArcorr} factorizes under the decomposition of Theorem \ref{thm:virasoroprimary} as
\begin{multline}
\Zcal_{X_k}^{A_0} \big(\varepsilon_1,\varepsilon_2, \mu_1,\mu_2; \qsf,
\vec{\xi} \ \big)_j \\
=\eta(\qsf)^{k-1}\ \chi_{\mathrm{conf}}^{\widehat{\omega}_j}\big(\qsf,\vec\zeta \ \big) \ \Big\langle \vacuum\,,\, \V_{-\frac{\mu_0}{\sqrt{-k\, \varepsilon_1\, \varepsilon
_2}},\frac{\mu_0+\varepsilon_1+\varepsilon_2}{\sqrt{-k\, \varepsilon_1\,
\varepsilon_2}}}(1) \ \V_{-\frac{\mu_1}{\sqrt{-k\, \varepsilon_1\, \varepsilon
_2}},\frac{\mu_1+\varepsilon_1+\varepsilon_2}{\sqrt{-k\, \varepsilon_1\,
\varepsilon_2}}}(\qsf) \vacuum\Big\rangle_{\mathbb{H}}\ .
\end{multline}
By Equation \eqref{eq:fundamental} we then obtain explicitly
\begin{equation}
\Zcal_{X_k}^{A_0}\big(\varepsilon_1,\varepsilon_2, \mu_0,\mu_1; \qsf,
\vec{\xi} \ \big)_j
=\eta(\qsf)^{k-1} \ \chi_{\mathrm{conf}}^{\widehat{\omega}_j}\big(\qsf,\vec\zeta \ \big) \
(1-\qsf)^{-\frac{\mu_1\,(\mu_0+\varepsilon_1+\varepsilon_2)}{k\, \varepsilon_1\, \varepsilon_2}} \ .
\end{equation}

\appendix

\bigskip \section{Virasoro primary fields}\label{app:virasoroprimary}

In this appendix we prove Theorem \ref{thm:virasoroprimary}. We need
to show that the vertex operator $\bar \V_\mu(\vec{v}_{21},\vec x, z)$ is a primary field of the Virasoro algebra
 generated by $L_n^{\slfrakhat_k}$ and $\cfrak$. We begin with the
 following result establishing the commutation relations between the
 Virasoro operators $L_n^{\widehat{\mathfrak{sl}}_k}$ introduced in
 Section \ref{sec:virasoro-affine} and the
 normal-ordered bosonic exponentials $\V_{1,-1}^{\gamma_{21}}(z)$ associated with the Heisenberg algebra $\hfrak_\Qfrak$.
\begin{lemma} For $n\neq0$ we have
\begin{align}
\big[L_n^{\slfrakhat_k}\,,\,\V_{1,-1}^{\gamma_{21}}(z)\big]&=z^n
\, \Big( z\, \partial_
z + \frac12\, \vec{v}_{21}\cdot C\vec{v}_{21} \, (n+1) +
\sum_{i=1}^{k-1} \, \big(\vec{v}_{21}\big)_i\,
 \big(\qfrak_0^i - z^{-n}\, \qfrak_n^i\big) \Big) \, \V_{1,-1}^{\gamma_{21}}(z)
 \ , \\[4pt]
\big[L_0^{\slfrakhat_k}\,,\,
\V_{1,-1}^{\gamma_{21}}(z)\big]&=z\, \partial_z \, \V_{1,-1}^{\gamma_{21}}(z) \ .
\end{align}
\end{lemma}
\proof
Let $\{\eta_i\}_{i=1}^{k-1}$ be an orthonormal basis of the vector
space $\Qfrak \otimes_\Z \R$. 
By the commutation relations \eqref{eq:commutationrelations}, we easily get
\begin{equation}\label{eq:commutationvirheis}
\big[L^{\slfrakhat_k}_n,\qfrak_m^{\eta_j}\big] = -m\, \qfrak_{n+m}^{\eta_j}
\end{equation}
for $n,m\in\Z$ and $ j=1,\ldots,k-1$. Fix a vector $\vec{v}\in\R^{k-1}$ and for $j=1,\ldots,k-1$ define
\begin{equation}
A_{\vec{v}}^j(z)_- := v_j \, \varphi_-^{\eta_j}(z) = v_j \, \sum_{m=
  1}^\infty\, \frac{z^m}{m} \,
 \qfrak_{-m}^{\eta_j} \qquad\mbox{and}\qquad
A_{\vec{v}}^j(z)_+ := - v_j\, \varphi_+^{\eta_j}(z) = - v_j\, \sum_{m= 1}^\infty\, \frac
{z^{-m}}{m}\, \qfrak_{m}^{\eta_j}\ .
\end{equation}
Using the commutation relations \eqref{eq:commutationvirheis}, we get
\begin{align}
\big[L^{\slfrakhat_k}_n, A_{\vec{v}}^j(z)_-\big] &= v_j\, \sum_{m=
  1}^\infty\, z^m \, \qfrak_{n-m}^{\eta_j} =: v_j\, \varphi_{-,n}^{\eta_j}(z)\ ,\\[4pt]
\big[L^{\slfrakhat_k}_n, A_{\vec{v}}^j(z)_+\big] &= v_j\, \sum_{m=
  1}^\infty\, z^{-m} \, \qfrak_{n+m}^{
\eta_j} =: v_j\, \varphi_{+,n}^{\eta_j}(z)\ .
\end{align}
For $n\leq 0$ the operator $\varphi_{-,n}^{\eta_j}(z)$ is a series in
the Heisenberg operators $\qfrak_l^{\eta_j}$ with $l<0$, thus it
commutes with $A_{\vec{v}}^j(z)_-$. For $n>0$, using again the relations \eqref{eq:commutationrelations} we get
\begin{equation}\label{eq:commAphin}
\big[ A_{\vec{v}}^j(z)_-, \varphi_{-,n}^{\eta_j}(z)\big] = v_j \, \sum_{m,l=1}^{n-
1} \, \frac{z^{m+l}}{m} \, \left[ \qfrak_{-m}^{\eta_j}, \qfrak_{n-l}^{\eta_j}\right] =
 -(n-1)\, v_j \, z^n \, \cfrak\ .
\end{equation}
Analogously, for $n\geq 0$ the operators $A_{\vec{v}}^j(z)_+$ and $\varphi_{+,n}^{\eta_j}(z)$ commute, while for $n<0$ we have
\begin{equation}
\big[A_{\vec{v}}^j(z)_+,\varphi_{+,n}^{\eta_j}(z)\big] = (n+1)\, v_j\, z^n\, \cfrak\ .
\end{equation}
Now we compute the commutator
\begin{align}
\big[L^{\slfrakhat_k}_n,\exp\big(A_{\vec{v}}^j(z)_- \big)\big] &= \sum_{l=
  1}^\infty\, \frac{1}{l!} \ \sum_{i=0
}^{l-1} \, A_{\vec{v}}^j(z)_-^i \,
\big[L^{\slfrakhat_k}_n,A_{\vec{v}}^j(z)_-\big] \, A_{\vec{v}}^j(z)_-^{l
-i-1} \\[4pt]
&= v_j \, \sum_{l= 1}^\infty\, \frac{1}{l!} \ \sum_{i=0}^{l-1} \,
A_{\vec{v}}^j(z)_-^i \, \varphi_{
-,n}^{\eta_j}(z) \, A_{\vec{v}}^j(z)_-^{l-i-1}\ .
\end{align}
For $n\leq 0$ we get simply
$\big[L^{\slfrakhat_k}_n,\exp\big(A_{\vec{v}}^j(z)_- \big)\big]=v_j\,
\varphi_{-,n}^{\eta_j}(z)\, \exp\big(A_{\vec{v}}^j(z)_- \big)$, while
for $n>0$ we can apply \eqref{eq:commAphin} iteratively to obtain
\begin{equation}
\big[L^{\slfrakhat_k}_n,\exp\big(A_{\vec{v}}^j(z)_- \big)\big]=\big(
v_j\, \varphi_{-,n}^{\eta_j}(z) - \mbox{$\frac12$}\, v_j^2\, (n-1) \,
z^n \,\big)\, \exp\big(A_{\vec{v}}^j(z)_- \big)\ .
\end{equation}
Noting that
\begin{equation}
z^{n+1}\, \partial_z \exp\big(A_{\vec{v}}^j(z)_- \big) = \left\{
\begin{array}{ll}
v_j\, \Big( \varphi_{-,n}^{\eta_j}(z) -\sum\limits_{t=1}^n \, z^t \, \qfrak_{n-t}^{\eta_
j} \Big)\, \exp\big(A_{\vec{v}}^j(z)_- \big) & n>0\\
v_j\, \Big( \varphi_{-,n}^{\eta_j}(z) -\sum\limits_{t=0}^{-n-1}\, z^{-t}\, \qfrak_{n+t
}^{\eta_j} \Big) \, \exp\big(A_{\vec{v}}^j(z)_- \big) \quad & n\leq 0
\end{array}
\right.\ ,
\end{equation}
and substituting in the previous expressions we finally get
\begin{equation}
\big[L^{\slfrakhat_k}_n,\exp\big(A_{\vec{v}}^j(z)_- \big)\big] = \left\{
\begin{array}{ll}
z^n \, \Big( z\, \partial_z + v_j \, \sum\limits_{s=0}^{n-1}\, z^{-s}
\, \qfrak_s^{\eta_j} - 
\frac12\, v_j^2\, (n-1) \Big)\exp\big(A_{\vec{v}}^j(z)_- \big) \quad & n>0\\
z^n \, \Big( z\, \partial_z - v_j \, \sum\limits_{s=1}^{-n}\, z^{s} \,
\qfrak_{-s}^{\eta_j} \Big)\exp\big(A_{\vec{v}}^j(z)_- \big) & n\leq 0
\end{array}
\right.\ .
\end{equation}
Repeating these computations for the operator $A_{\vec{v}}^j(z)_+$, we get analogous relations
\begin{equation}
\big[L^{\slfrakhat_k}_n,\exp\big(A_{\vec{v}}^j(z)_+ \big)\big] = \left\{
\begin{array}{ll}
z^n \, \Big( z\, \partial_z - v_j \, \sum\limits_{s=1}^{n}\, z^{-s} \,
\qfrak_{s}^{\eta_j} \Big) \exp\big(A_{\vec{v}}^j(z)_+ \big) & n\geq 0\\
z^n \, \Big( z\, \partial_z + v_j \, \sum\limits_{s=0}^{-n-1}\,  z^{s}
\, \qfrak_{-s}^{\eta_j}
 + \frac12\, v_j^2\, (n+1) \Big) \exp\big(A_{\vec{v}}^j(z)_+ \big) \quad & n<0
\end{array}
\right.\ .
\end{equation}
Now set
\begin{equation}
A_{\vec{v}}(z)_-:=\sum_{j=1}^{k-1} \, A_{\vec{v}}^j(z)_- \qquad\mbox{and}\qquad A_{\vec{v}}(z)_+:=\sum_{j=1}^{k-1} \, A_{\vec{v}}^j(z)_+\ .
\end{equation}
We are ready to compute
\begin{multline}
\big[L^{\slfrakhat_k}_n,\exp\big(A_{\vec{v}}(z)_- \big)\,
\exp\big(A_{\vec{v}}(z)_+ \big)\big]\\
=\, \sum_{j=1}^{k-1}\, \Big( \exp\big(A_{\vec{v}}^1(z)_-\big) \cdots
  \big[L^{\slfrakhat_k}_n,\exp(A_{\vec{v}}^j(z)_-)\big] \cdots
  \exp\big(A_{\vec{v}}^{k-1}(z)_-\big)\, \exp\big(A_{\vec{v}}(z)_+\big)\\
 + \, \exp\big(A_{\vec{v}}(z)_-\big)\, \exp\big(A_{\vec{v}}^1(z)_+\big) \cdots \big[L^{\slfrakhat_k}_n,
\exp\big(A_{\vec{v}}^j(z)_+\big)\big] \cdots
\exp\big(A_{\vec{v}}^{k-1}(z)_+ \big) \Big) \ .
\end{multline}
Fix $n<0$. Using the commutation relations computed before and noting
that the Heisenberg operator $\qfrak_{m}^{\eta_j}$ commutes with the
vertex operators $\exp\big(A_{\vec{v}}^l(z)_- \big)$ and
$\exp\big(A_{\vec{v}}^l(z)_+ \big)$ for $l\neq j$ and any $m\in\Z\setminus\{0\}$, we obtain
\begin{multline}
\big[L^{\slfrakhat_k}_n,\exp\big(A_{\vec{v}}(z)_-\big)\,
\exp\big(A_{\vec{v}}(z)_+ \big)\big]\\
=\, z^n\, \Big(z\, \partial_z + 
\frac12 \, (n+1) \, \sum_{j=1}^{k-1}\, v_j^2 - \sum_{j=1}^{k-1}\, v_j
\ \sum_
{t=1}^{-n} \, z^t \, \qfrak_{-t}^{\eta_j}\Big)
\Big(\exp\big(A_{\vec{v}}(z)_-\big)\, \exp \big(A_{\vec
{v}}(z)_+\big)\Big)\\
+\,  \sum_{j=1}^{k-1} v_j \, \exp\big(A_{\vec{v}}(z)_- \big)\, \Big(
\, \sum_{t=0}^{-n-1} \, z^t \, \qfrak_{
-t}^{\eta_j} \, \Big)\, \exp\big(A_{\vec{v}}(z)_+ \big)\ .
\end{multline}
Since the Heisenberg operators $\qfrak_{-t}^{\eta_j}$ for $t\geq0$
also commute with $\exp\big(A_{\vec{v}}^j(z)_- \big)$, we thus find
\begin{multline}\label{eq:commVirexpexp}
\big[L^{\slfrakhat_k}_n,\exp\big(A_{\vec{v}}(z)_- \big)\,
\exp\big(A_{\vec{v}}(z)_+ \big)\big]\\ =\, z^n\, \Big(z\, \partial
_z + \frac12\, (n+1)\, \sum_{j=1}^{k-1}\, v_j^2 - \sum_{j=1}^{k-1} \,
v_j \, \big(\qfrak_{0}^{\eta_j} - z
^{-n}\, \qfrak_{n}^{\eta_j} \big)\Big) \Big(\exp\big(A_{\vec{v}}(z)_-
\big)\, \exp\big(A_{\vec{v}}(z)_+\big) \Big) \ .
\end{multline}
For $n>0$ we arrive at the similar expression
\begin{multline}
\big[L^{\slfrakhat_k}_n,\exp\big(A_{\vec{v}}(z)_- \big)\,
\exp\big(A_{\vec{v}}(z)_+ \big)\big] \\ =\, z^n\, \Big(z\, \partial_z
- \frac12\, (n-1)\, \sum_{j=1}^{k-1}\,  v_j^2 + \sum_{j=1}^{k-1} \,
v_j \ \sum_{
t=0}^{n-1} \, z^{-t} \,
\qfrak_t^{\eta_j}\Big)\Big(\exp\big(A_{\vec{v}}(z)_- \big)\, \exp\big(A_{\vec
{v}}(z)_+\big) \Big)\\
- \, \sum_{j=1}^{k-1}\,  v_j \, \exp\big(A_{\vec{v}}(z)_-\big)\,
\Big(\, \sum_{t=1}^{n} \, z^{-t} \, \qfrak_t
^{\eta_j} \, \Big)\, \exp\big(A_{\vec{v}}(z)_+ \big)\ ,
\end{multline}
but this time the operators $\qfrak_t^{\eta_j}$ for $t>0$ do not
commute with $\exp\big(A_{\vec{v}}^j(z)_- \big)$. Since
\begin{equation}
\Big[A_{\vec{v}}^j(z)_-\,,\,\sum_{t=1}^{n} \, z^{-t} \, \qfrak_t^{\eta_j}\Big] = -n\,v_j\,
\cfrak\ ,
\end{equation}
we get
\begin{equation}
\Big[\exp\big(A_{\vec{v}}^j(z)_- \big)\,,\,\sum_{t=1}^{n} \, z^{-t} \,
\qfrak_t^{\eta_j}\Big] = -n\, v
_j\, \exp\big(A_{\vec{v}}^j(z)_- \big)
\end{equation}
and thus we arrive again at Equation \eqref{eq:commVirexpexp}. For $n=0$ we obtain
\begin{equation}
\big[L^{\slfrakhat_k}_0,\exp\big(A_{\vec{v}}(z)_-\big)\,
\exp\big(A_{\vec{v}}(z)_+ \big)\big] = z\, \partial_z \Big(
\exp\big(A_{\vec{v}}(z)_- \big)\, \exp\big(A_{\vec{v}}(z)_+\big) \Big)\ .
\end{equation}
Finally, to get the assertion it is sufficient to note that if
$D=(d_{ij})$ is the change of basis matrix such that
$\gamma_i=\sum_{j=1}^{k-1}\, d_{ij}\, \eta_j$ with $d_{ij}\in\R$, 
then $\V_{1,-1}^{\gamma_{21}}(z) =
\exp\big(A_{D\vec{v}_{21}}(z)_-\big)\, \exp\big(A_{D\vec{v}_{2
1}}(z)_+ \big)$, and moreover
\begin{align}
\sum_{j=1}^{k-1}\, \big( D\vec{v}_{21} \big)_j\, \eta_j&=
\sum_{i=1}^{k-1} \, (\vec
{v}_{21})_i\, \gamma_i \ , \\[4pt]
\sum_{j=1}^{k-1}\, \big( D\vec{v}_{21} \big)_j^2&=\Big\langle \, \sum_{j=1}^{k-1}\, \big(D\vec{v}_{21}\big)_j\,\eta_j\,,\, \sum_{j=1}^{k-1}\, \big(D\vec{v}_{21}\big)_
j\,\eta_j\, \Big\rangle_{\Qfrak\otimes_\Z\R}=\vec{v}_{21}\cdot C\vec{v}_{21}\ .
\end{align}
\endproof

\begin{remark}
A similar (but simpler) calculation shows that the vertex operators $\V_{\alpha,\beta}(z)$ are primary fields in the sense stated in Remark \ref{rem:genbosprimary}.
\end{remark}

The proof of Theorem \ref{thm:virasoroprimary} is now completed once
we establish the following commutation relations.
\begin{lemma}
For any $n\in \Z$ we have
\begin{equation}
\big[L^{\slfrakhat_k}_n,\bar \V_\mu(\vec{v}_{21},\vec x, z)\big]=z^n \,
\big( z\, \partial_z + \mbox{$\frac12$}\, \vec{v}_{21}\cdot
C\vec{v}_{21}\, n \big)\bar\V_\mu(\vec{v}_{21},\vec x, z)\ .
\end{equation}
\end{lemma}
\proof
To get the assertion, it is enough to derive the commutation relations
\begin{align}
\big[L^{\slfrakhat_k}_n,\exp\big(\log z\
\cfrak+\gamma_{21}\,\big)\big]&=z^n\, \Big( z\, \partial_z + \frac12\,
\vec{v}_{21}\cdot C\vec{v}_{21} \\ & \qquad \qquad \qquad \qquad -\, \sum_{i=1}^{k-1} \, \big(\vec
{v}_{21} \big)_i \, \big(\qfrak_0^i - z^{-n}\, \qfrak_n^i\big) \Big)
\exp\big(\log z\ \cfrak+\gamma_{21}\,\big) \ , \\[4pt]
\big[L^{\slfrakhat_k}_0, \exp\big(\log z\
\cfrak+\gamma_{21}\,\big)\big]&=z\, \partial_z \exp\big(\log z\ \cfrak+\gamma_{21}\,\big)\ .
\end{align}
\endproof

\bigskip \section{Edge contributions\label{app:edgecontributions}}

In this appendix we begin by listing the edge contributions
\begin{equation}
\sum_{n=
1}^{k-1}\,
L^{(n)}_{\vec{v}}
\big(\varepsilon_1^{(n)},\varepsilon_2^{(n)}\big)
\end{equation}
to the $T$-equivariant Chern character of the natural bundle
$\Vbf^{\vec{u},n,j} $ on $\Mcal(\vec{u},n;j)$
which were derived in \cite[Appendix C]{art:bruzzopedrinisalaszabo2013}.  For this, we
first introduce some notation. Let
$j\in\{0,1,\dots,k-1\}$ be the equivalence class of $k\,
v_{k-1}$ modulo $k$. Set $(C^{-1})^{n 0}=0$ for
$n\in\{1, \ldots, k-1\}$ and $(C^{-1})^{k,j}=0$. We also set $\vec
s:=C^{-1}(\vec u-\boldsymbol{e}_j)$ if $j>0$ and $\vec s:= \vec v$ if
$j=0$; then $\vec s\in\Z^{k-1}$. We denote by $\lfloor x\rfloor\in\Z$
  the integer part and by $\{x\}:=x-\lfloor
  x\rfloor\in[0,1)$ the fractional part  of a rational 
  number $x$.

If $s_ n\geq 0$ for every $ n=1, \ldots, k-1$, consider the equation
\begin{equation}\label{eq:condition-index+}
\frac{C_{ n  n}}{2}\, i^2- i\, \Big(\vec{v}-\sum_{p=1}^{ n-1}\, s_p\,
\boldsymbol{e}_p\Big)\cdot C\boldsymbol{e}_ n+\frac{1}{2}\,
\bigg(\Big(\vec{v}-\sum_{p=1}^{ n-1}\, s_p\, \boldsymbol{e}_p\Big)\cdot
C\Big(\vec{v}-\sum_{p=1}^{ n-1}\, s_p\, \boldsymbol{e}_p\Big)-\big(C^{-1} \big)^{c c}\bigg)=0\ ,
\end{equation}
and define the set 
\begin{equation}
S_ n^+:=\{i\in\N\, \vert\, i\leq s_ n \ \mbox{ is a solution of Equation \eqref{eq:condition-index+}}\}\ .
\end{equation}
Let $d_ n^+:= \min(S_ n^+)$ if
$S_ n^+\neq \emptyset$ and $d_ n^+:=s_ n$ otherwise.

When $s_ n<0$ for $ n=1,\ldots, k-1$ consider the equation
\begin{equation}\label{eq:condition-index-}
\frac{C_{ n  n}}{2}\, i^2+i\, \Big(\vec{v}-\sum_{p=1}^{ n-1}\, s_p\,
\boldsymbol{e}_p\Big)\cdot C\boldsymbol{e}_ n +\frac{1}{2}\,
\bigg(\Big(\vec{v}-\sum_{p=1}^{ n-1}\, s_p\, \boldsymbol{e}_p\Big)\cdot
C\Big(\vec{v}-\sum_{p=1}^{ n-1}\, s_p\, \boldsymbol{e}_p\Big)-\big(C^{-1} \big)^{c c}\bigg)=0\ ,
\end{equation}
and define the set 
\begin{equation}
S_ n^-:=\{i\in\N\, \vert\, i\leq -s_ n \ \mbox{ is a solution of Equation \eqref{eq:condition-index-}}\}\ .
\end{equation}
Let $d_ n^-:=
\min(S_ n^-)$ if $S_ l^-\neq \emptyset$ and $d_ n^-:=-s_ n$
otherwise. Let $m$ be the smallest integer $n\in\{1, \ldots, k-1\}$
such that $S_ n^+$ or $S_ n^-$ is nonempty; if all of these sets are
empty,  let $m:=k-1$.

Then for fixed $n=1, \ldots, m$ we set:
\begin{itemize} \setlength{\itemsep}{4mm}
\item[\scriptsize$\blacksquare$] For $v_ n-(C^{-1})^{ n j}> 0$:
\smallskip
\begin{itemize} \setlength{\itemsep}{0.8cm}
\item[$\bullet$]For $\delta_{ n,j} -v_{ n+1}+(C^{-1})^{ n+1, j}+2(v_{
    n}-(C^{-1})^{ n j}-d_ n^+)\geq 0$:
\small
\begin{equation}
L_{\vec v}^{(n)}\big(\varepsilon_1^{(n)},\varepsilon_2^{(n)} \big) =- \sum_{i=v_{
    n}-(C^{-1})^{ n j}-d_ n^+}^{v_ n-(C^{-1})^{ n j}-1}\hspace{0.2cm}
\sum_{j=0}^{2i+\delta_{ n,j} -v_{ n+1}+(C^{-1})^{ n+1, j}}\,
\big(\chi_1^ n \big)^{i+\big\lfloor\frac{\delta_{ n,j} -v_{
      n+1}+(C^{-1})^{n+1, j}}{2} \big\rfloor}\, \big(\chi_2^ n \big)^j\ .
\end{equation}
\normalsize
\item[$\bullet$]For $2\leq \delta_{ n,j}-v_{ n+1}+(C^{-1})^{ n+1,
    j}+2(v_ n-(C^{-1})^{ n j})<2 d_ n^+$:
\small
\begin{multline}
L_{\vec v}^{(n)}\big(\varepsilon_1^{(n)},\varepsilon_2^{(n)} \big) \\
\shoveleft{=\sum_{i=v_ n-(C^{-1})^{ n j}-d_ n^+}^{-\big\lfloor
    \frac{\delta_{ n,j} -v_{ n+1}+(C^{-1})^{ n+1, j}}{2}\big\rfloor -
    1}\hspace{6mm} \sum_{j=1}^{2i- (\delta_{ n,j} -v_{
      n+1}+(C^{-1})^{ n+1, j})-1} \, \big(\chi_1^ n
  \big)^{i-\big\lfloor -\frac{\delta_{ n,j} -v_{ n+1}+(C^{-1})^{ n+1,
        j}}{2}\big\rfloor}\, \big(\chi_2^ n \big)^{-j}} \\[4pt]
 - \, \sum_{i=-\big\lfloor \frac{\delta_{ n,j} -v_{ n+1}+(C^{-1})^{
       n+1, j}}{2}\big\rfloor}^{2(v_ n-(C^{-1})^{ n j})+\delta_{ n,j}
   -v_{ n+1}+(C^{-1})^{ n+1, j} -2}
 \hspace{6mm}\sum_{j=0}^{2i+\delta_{ n,j}-v_{ n+1}+(C^{-1})^{ n+1,
     j}}\, \big(\chi_1^ n \big)^{i+\big\lfloor \frac{\delta_{ n,j}
     -v_{ n+1}+(C^{-1})^{ n+1, j}}{2}\big\rfloor} \, \big(\chi_2^ n \big)^j\ .
\end{multline}
\normalsize
\item[$\bullet$]For $\delta_{ n,j} -v_{ n+1}+(C^{-1})^{n+1, j}<2-2(v_
  n-(C^{-1})^{ n j})$:
\small
\begin{multline}
L_{\vec v}^{(n)}\big(\varepsilon_1^{(n)},\varepsilon_2^{(n)} \big) \\
=\sum_{i=v_{
    n}-(C^{-1})^{ n j}-d_ n^+}^{v_ n-(C^{-1})^{ n j}-1}\hspace{0.2cm}
\sum_{j=1}^{-2i-\delta_{ n,j} +v_{ n+1}-(C^{-1})^{n+1,j}-1} \,
\big(\chi_1^ n \big)^{i-\big\lfloor - \frac{\delta_{ n,j} -v_{
      n+1}+(C^{-1})^{n+1,j}}{2} \big\rfloor} \, \big(\chi_2^ n \big)^{-j} \ .
\end{multline}
\normalsize
\end{itemize}
\item[\scriptsize$\blacksquare$] For $v_ n-(C^{-1})^{ n j}= 0$:
  $\qquad L_{\vec v}^{(n)} \big(\varepsilon_1^{(n)},\varepsilon_2^{(n)} \big)=0$.
\item[\scriptsize$\blacksquare$] For $v_ n-(C^{-1})^{ n j}< 0$:
\smallskip
\begin{itemize}\setlength{\itemsep}{0.5cm}
\item[$\bullet$]For $\delta_{ n,j} -v_{ n+1}+(C^{-1})^{n+1,j} +2v_{
    n}-2(C^{-1})^{ n j} < 2-2 d_ n^{-}$:
\small
\begin{multline}
L_{\vec v}^{(n)}\big(\varepsilon_1^{(n)},\varepsilon_2^{(n)} \big) \\
=-\sum_{i=1-v_{
    n}+(C^{-1})^{ n j}-d_ n^-}^{-v_{ n}+(C^{-1})^{ n j}}\hspace{0.2cm}
\sum_{j=1}^{2i-(\delta_{ n,j} -v_{ n+1}+(C^{-1})^{n+1,j})-1} \,
\big(\chi_1^ n \big)^{-i-\big\lfloor -\frac{\delta_{ n,j} -v_{
      n+1}+(C^{-1})^{n+1,j}}{2}\big\rfloor} \, \big(\chi_2^ n \big)^{-j}\ .
\end{multline}
\normalsize
\item[$\bullet$]For $2-2 d_ n^{-}\leq \delta_{ n,j} -v_{
    n+1}+(C^{-1})^{n+1,j}+2v_{ n}-2(C^{-1})^{ n j}<0$:
\small
\begin{multline}
L_{\vec v}^{(n)}\big(\varepsilon_1^{(n)},\varepsilon_2^{(n)} \big) \\
\shoveleft{=\sum_{i=1-v_{ n}+(C^{-1})^{ n j}-d_ n^-}^{\big\lfloor
    \frac{\delta_{ n,j} -v_{ n+1}+(C^{-1})^{n+1,j}}{2}
    \big\rfloor}\hspace{0.2cm} \sum_{j=0}^{-2i+\delta_{ n,j} -v_{
      n+1}+(C^{-1})^{n+1,j}} \, \big(\chi_1^ n \big)^{-i+\big\lfloor
    \frac{\delta_{ n,j} -v_{ n+1}+(C^{-1})^{n+1,j}}{2} \big\rfloor} \,
  \big(\chi_2^ n \big)^j}\\
-\, \sum_{i=\big\lfloor \frac{\delta_{ n,j} -v_{
      n+1}+(C^{-1})^{n+1,j}}{2} \big\rfloor+1}^{-v_{ n}+(C^{-1})^{ n
    j}}\hspace{0.2cm} \sum_{j=1}^{2i-(\delta_{ n,j} -v_{
    n+1}+(C^{-1})^{n+1,j})-1 } \, \big(\chi_1^ n \big)^{-i-
  \big\lfloor -\frac{\delta_{ n,j} -v_{
      n+1}+(C^{-1})^{n+1,j}}{2}\big\rfloor} \, \big(\chi_2^ n \big)^{-j}\ .
\end{multline}
\normalsize
\item[$\bullet$]For $\delta_{ n,j} -v_{ n+1}+(C^{-1})^{n+1,j}\geq
  -2v_{ n}+2(C^{-1})^{ n j}$:
\small
\begin{equation}
L_{\vec v}^{(n)}\big(\varepsilon_1^{(n)},\varepsilon_2^{(n)} \big) =\sum_{i=1-v_{
    n}+(C^{-1})^{ n j}-d_ n^-}^{-v_{ n}+(C^{-1})^{ n j}}\hspace{0.2cm}
\sum_{j=0}^{-2i+\delta_{ n,j}-v_{ n+1}+(C^{-1})^{n+1,j}} \,
\big(\chi_1^ n \big)^{-i+\big\lfloor \frac{\delta_{ n,j} -v_{
      n+1}+(C^{-1})^{n+1,j}}{2}\big\rfloor}\, \big(\chi_2^ n \big)^{j}\ .
\end{equation}
\normalsize
\end{itemize}
\end{itemize}
For $n=m+1, \ldots, k-1$ we set $
L^{(n)}_{\vec{v}}\big(\varepsilon_1^{(n)},
\varepsilon_2^{(n)} \big)=0$. Note that for any fixed $n\in\{1,\dots,k-1\}$,
$d_n^\pm=0$ implies $L^{(n)}_{\vec{v}}
\big(\varepsilon_1^{(n)},\varepsilon_2^{(n)}
\big)=0$.

The edge factors $\ell^{(n)}_{\vec{v}_{21}}
\big(\varepsilon_1^{(n)},\varepsilon_2^{(n)},\mu \big)$ which
contribute to the $T$-equivariant Euler class of the Carlsson-Okounkov bundle
$\Ebf_\mu^{\vec{u}_1,n_1,j_1;\vec{u}_2,n_2,j_2} $ on
$\Mcal(\vec{u}_1,n_1;j_1)\times\Mcal(\vec{u}_2,n_2;j_2)$
are then obtained in the following way. We replace $\vec v$ in the
above by $\vec v_{21}$ and $j$ by $j_{21}$. If
\begin{equation}
L^{(n)}_{\vec{v}_{21}}\big(\varepsilon_1^{(n)},
\varepsilon_2^{(n)} \big) = \sum_{i=1}^{D}\, \eta_i\, \e^{\sigma_i}
\end{equation}
with $\eta_i=0,\pm\,1$, then
\begin{equation}
\ell^{(n)}_{\vec{v}_{21}}
\big(\varepsilon_1^{(n)},\varepsilon_2^{(n)},\mu \big) =
\prod_{i=1}^D\, (\mu+\sigma_i)^{\eta_i} \ .
\end{equation}
Explicit formulas are written in \cite[Section
4.7]{art:bruzzopedrinisalaszabo2013}.

The contribution of $L^{(n)}_{\vec{v}_{21}}\big(\varepsilon_1^{(n)},
\varepsilon_2^{(n)} \big)$ to the $p$-th equivariant Chern class $(\cc_{p} )_T\big(\Vbf^j \big)$ is gotten by extracting the monomial terms of total degree $p$ in $\varepsilon_1,\varepsilon_2$. In particular, the contribution to the first Chern class is given by
\begin{equation}
\ell_{\vec{v}}^{(n)}\big(\varepsilon_1^{(n)},
  \varepsilon_2^{(n)}\big)_{[1]} = \Big(\varepsilon_1\, \left.\frac\partial{\partial\varepsilon_1}\,\right|_{\varepsilon_1=0}+ \varepsilon_2\, \left.\frac\partial{\partial\varepsilon_2}\,\right|_{\varepsilon_2=0} \Big)L^{(n)}_{\vec{v}_{21}}\big(\varepsilon_1^{(n)},
\varepsilon_2^{(n)} \big) \ .
\end{equation}

\begin{example}
Let $k=2$. Then $j_{21}\in\{0,1\}$, while
$\{v_{21}\}=\frac12\, \delta_{1,j_{21}}$ and
$\lfloor v_{21}\rfloor = v_{21}-(C^{-1})^{1,
  j_{21}}$. Since $m=1$, $d_1^+=\lfloor
v_{21}\rfloor$ and $d_1^-=-\lfloor v_{21}\rfloor$,
and we get
\begin{equation*}
\ell_{v_{21}}\big(\varepsilon_1,\varepsilon_2,
\mu \big)=\left\{
\begin{array}{cl}
\prod_{i=0}^{\lfloor v_{21}\rfloor-1}\limits\hspace{0.2cm}
\prod_{j=0}^{2i+2\{v_{21}\}}\limits\big(\mu+i\,
\varepsilon_1+j\, \varepsilon_2\big) & \mbox{ for } \lfloor v_{21}\rfloor> 0\ ,\\[8pt]
1 & \mbox{ for } \lfloor v_{21}\rfloor= 0\ ,\\[8pt]
\prod_{i=1}^{-\lfloor v_{21}\rfloor}\limits\hspace{0.2cm} \prod_{j=1}^{2i-2\{v_{21}\}-1}\limits
\big(\mu+(2\{v_{21}\}-i)\,
\varepsilon_1-j\, \varepsilon_2\big)  & \mbox{ for } \lfloor v_{21}\rfloor< 0 \ .
\end{array}
\right.
\end{equation*}
For $\{v_{21}\}=0$ these formulas coincide with the blowup factors obtained in \cite{art:gasparimliu2010} up to a redefinition of the equivariant parameters (see also \cite{art:ciraficiszabo2012}). Moreover, for $\lfloor v_{21}\rfloor>0$ they can be easily written in the form
\begin{equation}
\ell_{v_{21}}\big(\varepsilon_1,\varepsilon_2,
\mu \big) = \prod_{\stackrel{\scriptstyle i,j\geq1\,,\,i+j\leq 2\lfloor v_{21}\rfloor}{\scriptstyle i+j\equiv 0\ {\rm mod}\,2}}\, \big(\mu+(i-1)\, \tilde\varepsilon_1+(j-1)\, \tilde\varepsilon_2\big)
\end{equation}
with $\tilde\varepsilon_1=\frac{\varepsilon_1}2$ and $\tilde\varepsilon_2=\frac{\varepsilon_1}2+\varepsilon_2$, which coincide with the blowup factors of \cite{art:bonellimaruyoshitanzini2011, art:bonellimaruyoshitanzini2012, art:belavinbershteinfeiginlitvinovtarnopolsky2011} (similarly for $\lfloor v_{21}\rfloor<0$ and/or $\{v_{21}\}=\frac12$). In \cite{art:belavinbershteinfeiginlitvinovtarnopolsky2011} it is stated that these edge factors can be represented as suitable matrix elements of primary fields from Theorem \ref{thm:virasoroprimary} in highest weight states of $\glfrakhat_2=\hfrak\oplus\slfrakhat_2$ at level one; the proof makes use of the Frenkel-Kac construction and the Dotsenko-Fateev integrals of \cite{art:felder1989}.
\end{example}

%\bibliography{article}
%\bibliographystyle{siam}

\end{document}